\documentclass[11pt]{article}
\usepackage{amsthm,amsfonts,amsmath,amssymb}
\usepackage{tabularx}
\usepackage{multirow,topcapt}
\newtheorem{thm}{Theorem}[section]

\newtheorem{lem}[thm]{Lemma}

\def\di{\bigm|} \def\lg{\langle} \def\rg{\rangle}
\def\nd{\mathrel{\bigm|\kern-.7em/}}

\def\f{\noindent}

\def\mod{\hbox{\rm mod }}

\def\demo{\f{\bf Proof}\hskip10pt}

\def\diag{\hbox{\rm diag}}

\def\qed{\hfill $\Box$}

\def\lg{\langle}
\def\rg{\rangle}

\def\rr#1{\item[{\rm (#1)}]}

\numberwithin{equation}{section}
\numberwithin{table}{section}

\begin{document}
\title{Finite $p$-groups with a minimal non-abelian
subgroup of index $p$ (IV)
\thanks{This work was supported by NSFC (no. 11371232),
by NSF of Shanxi Province (no. 2012011001-3 and 2013011001-1). }}
\author{Lijian An\thanks{Corresponding author. e-mail: anlj@sxnu.edu.cn}, Ruifang Hu and Qinhai Zhang
\\
Department of Mathematics,
Shanxi Normal University\\
Linfen, Shanxi 041004, P. R. China }

\maketitle

\begin{abstract}
In this paper, we completely classify the finite $p$-groups $G$ such
that $\Phi(G')G_3\le C_p^2$, $\Phi(G')G_3\le Z(G)$ and
$G/\Phi(G')G_3$ is minimal non-abelian. This paper is a part of the
classification of finite $p$-groups with a minimal non-abelian
subgroup of index $p$. Together with other four papers, we solve a problem proposed by
Y. Berkovich.

\medskip


\end{abstract}

\baselineskip=16pt

\section{Introduction}

Groups in this paper are finite $p$-groups. Notation and
terminology are consistent with that in \cite{ALQZ,QXA,QYXA}. For a positive integer $t$, we use $G\in\mathcal{A}_t$ to denote that $G$ is a non-abelian $p$-group and the maximal index of minimal non-abelian subgroups of $G$ is $p^{t-1}$. So $\mathcal{A}_1$-groups just mean the minimal non-abelian $p$-groups.
 In \cite{Yak}, Y. Berkovich proposed the following

\medskip

\noindent{\bf Problem 239:} Classify the $p$-groups containing an
$\mathcal{A}_1$-subgroup of index $p$.

\medskip

 This problem may be divided into two parts:

\medskip

\noindent{\bf Part 1.} Classify the finite $p$-groups with at least two
$\mathcal{A}_1$-subgroups of index $p$;

\noindent{\bf Part 2.} Classify the finite $p$-groups with a unique
$\mathcal{A}_1$-subgroup  of index $p$.

In \cite{QYXA} and \cite{QZGA}, we solved Part 2 for $p>2$ and $p=2$
respectively. In \cite{ALQZ} and \cite{QXA}, we
solved Part 1 for $d(G)=3$. So, for problem 239, it remains the groups of Part 1 with $d(G)=2$. For such groups, we
can prove that $c(G)\le 3$, $G/\Phi(G')G_3\in \mathcal{A}_1$ and $\Phi(G')G_3\le
C_p^2$ (Lemma \ref{part1}). Hence $G$ is a central extension of an elementary abelian
$p$-group by an $\mathcal{A}_1$-group.
If $G/\Phi(G')G_3$ is metacyclic, then, by \cite[Theorem 2.3]{Bla}, $G$ is metacyclic and hence
$G'\cong C_p^2$. Furthermore, by \cite[Lemma J(i)]{Yak0}, $G$ is an
$\mathcal{A}_2$-group. Hence we only need to deal with the case that
$G/\Phi(G')G_3$ is non-metacyclic. By \cite{R}, $G/\Phi(G')G_3$ is isomorphic to $M_p(n,m,1):=\lg a,b,c\di a^{p^n}=b^{p^m}=c^p=1,[a,b]=c,[c,a]=[c,b]=1\rg$ where $n>1$ for $p=2$ and $n\ge m$.
In this case, it is
independently interesting to classify the groups satisfying the
following property:

\medskip

\noindent{\bf Property $\mathcal{P}$.} $\Phi(G')G_3\le C_p^2$,
$\Phi(G')G_3\le Z(G)$ and $G/\Phi(G')G_3\cong M_p(n,m,1)$.

In this paper, we classify the groups satisfying Property $\mathcal{P}$.
As a direct application, we solve Part 1 of problem 239 for
$d(G)=2$, and completely solve problem 239 together with
\cite{ALQZ,QXA,QYXA,QZGA}.

By the way, we also look into some properties for those groups with
Property $\mathcal{P}$. For the case $\Phi(G')G_3=C_p$, we give the minimal and the
maximal index of $\mathcal{A}_1$-subgroups. For the case $\Phi(G')G_3=C_p^2$, we pick up the groups having an $\mathcal{A}_1$-subgroup of index $p$ and $\mathcal{A}_3$-groups respectively. These results are
useful in the classification of $\mathcal{A}_3$-groups \cite{zhang}
.

\section{Preliminaries}

In this paper, $p$ is always a prime. We use $F_p$ to denote the
finite field containing $p$ elements. $F_p^*$ is the multiplicative
group of $F_p$. $(F_p^*)^2=\{ a^2\di a\in F_p^*\}$ is a subgroup of
$F_p^*$. $F_p^2=(F_p^*)^2\cup\{0\}$. For a finite non-abelian $p$-group $G$, we use $p^{I_{\min}}$ and $p^{I_{\max}}$ to denote the minimal index and the maximal index of $\mathcal{A}_1$-subgroups of $G$ respectively. 
For a square matrix $A$, $|A|$ denotes the determinant of $A$.  We need the following lemmas.

\begin{lem}{\rm (\cite[Lemma 2.2]{Alj})}
\label{minimal non-abelian equivalent conditions} Suppose that $G$
is a finite non-abelian $p$-group. Then the following conditions are
equivalent:

{\rm (1)} $G$ is an $\mathcal{A}_1$-group;%anlijian

{\rm (2)} $d(G)=2$ and $|G'|=p$;

{\rm (3)} $d(G)=2$ and $\Phi(G)=Z(G)$.
\end{lem}

%\begin{lem}{\rm (\cite[Theorm 2.3]{Bla})}
%\label{metacyclic equivalent conditions} A $p$-group $G$ is metacyclic if and only if $G/\Phi(G')G_3$ is metacyclic.
%\end{lem}

\begin{lem}\label{part1}{\rm (\cite[Theorem 2.7]{ALQZ})} Suppose that $G$ has at least two
$\mathcal{A}_1$-subgroups of index $p$. The the following conclusions hold:

{\rm (1)} If $d(G)=3$, then $\Phi(G)\le Z(G)$;

{\rm (2)} If $d(G)=2$, then $c(G)\le 3$, $\Phi(G')G_3\le C_p^2$ and
$G/\Phi(G')G_3\in\mathcal{A}_1$.

\end{lem}

\begin{lem}
\label{cor 2.4}{\rm (\cite[Corollary 2.4]{ALQZ})} {\rm (1)} Let $M$ be an $\mathcal {A}_t$-group, and
$A$ be an abelian group of order $p^k$. Then $G=M\times A$ is an
$\mathcal {A}_{t+k}$-group.

{\rm (2)} Let $M$ be an $\mathcal {A}_t$-group with $|M'|=p$,
$G=M\ast A$, where $A$ is abelian with order $p^{k+1}$ and $M\cap
A=M'$. Then $G$ is an $\mathcal {A}_{t+k}$-group.
\end{lem}

\begin{lem}\label{hetong-p}{\rm (\cite[Lemma 2.1]{QXA})} Suppose that $p$ is odd, $\{1,\eta\}$ is a transversal for
$(F_p^*)^2$ in $F_p^*$. Then the following matrices form a
transversal for the congruence classes of
invertible matrices of order $2$ over $F_p$:

\medskip
\ \ \ {\rm (1)} $\left(
 \begin{array}{cc}
 0& 1\\
 -1& 0
 \end{array}
 \right)$, \ \ \
{\rm (2)} $\left(
 \begin{array}{cc}
 \nu  & 1\\
 -1& 0
 \end{array}
 \right)$,\ \ \
{\rm (3)} $\left(
 \begin{array}{cc}
 1& 0\\
 0& \nu
 \end{array}
 \right)$,\ \ \
{\rm (4)} $\left(
 \begin{array}{cc}
 1& 1\\
 -1& r
 \end{array}
 \right)$,

\medskip

\noindent where $\nu=1$ or $\eta$, $r=1,2,\dots,p-2$.
 \end{lem}

\begin{lem}{\rm (\cite[Lemma 4.3]{ALQZ})}\label{hetong-2}
The following matrices form a transversal for the congruence classes of
invertible matrices of order $2$ over $F_2$.

\medskip
 \ \ \ {\rm (1)} $\left(
 \begin{array}{cc}
 1& 0\\
 0& 1
 \end{array}
 \right)$, \ \ \
{\rm (2)} $\left(
 \begin{array}{cc}
 0& 1\\
 1& 0
 \end{array}
 \right)$, \ \ \
{\rm (3)} $\left(
 \begin{array}{cc}
 1& 0\\
 1& 1
 \end{array}
 \right)$.
\end{lem}

\begin{lem}\label{hetong-p-non-invertible}{\rm (\cite[Lemma 2.3]{QXA})} Suppose that $p$ is a prime $(p=2$ is possible$)$. For odd $p$, $\{1,\eta\}$ is a transversal for
$(F_p^*)^2$ in $F_p^*$. Then the following matrices form a
transversal for the congruence classes of
non-invertible matrices of order $2$ over $F_p$:

\medskip
\ \ \ {\rm (1)} $\left(
 \begin{array}{cc}
 0& 1\\
 0& 0
 \end{array}
 \right)$, \ \ \
{\rm (2)} $\left(
 \begin{array}{cc}
 0& 0\\
 0& \nu
 \end{array}
 \right)$,\ \ \
{\rm (3)} $\left(
 \begin{array}{cc}
 0& 0\\
 0& 0
 \end{array}
 \right)$,
where $\nu=1$ or $\eta$.
 \end{lem}

\section{The case $G_3\le \Phi(G')\cong C_p$}

Suppose that $G$ is a finite $p$-group with $G_3\le \Phi(G')\cong C_p$ and $G/\Phi(G')\cong M_p(n,m,1)$, where $n>1$ for $p=2$ and
$n\ge m$. Let $$G/\Phi(G')=\lg \bar{a},\bar{b},\bar{c}\mid \bar{a}^{p^n}=\bar{b}^{p^m}=\bar{c}^{p}=1,[\bar{c},\bar{a}]=[\bar{b},\bar{c}]=1\rg.$$
Then, without loss of generality, we may assume that $G=\lg a,b,c\rg$ such that $[a,b]=c$. Since $G_3\le \Phi(G')\cong C_p$, $c$ is of order $p^2$ and $\Phi(G')=\lg c^p\rg$.
Since $a^{p^{n}}\in \Phi(G')$, we may assume that $a^{p^{n}}=c^{w_{11}p}$. By similar reasons, we may assume that $b^{p^m}=c^{w_{21}p}$, $[c,a]=c^{w_{12}p}$ and $[c,b]=c^{w_{22}p}$.
Then we get a $2\times 2$ matrix over $F_p$. $w(G)=(w_{ij})$ is called a characteristic matrix of $G$. Notice that $w(G)$ will be changed if we change the generators $a,b$. We also call $a,b$ a set of characteristic generators of $w(G)$.

In this paper, $w(G)$ is always a characteristic matrix of $G$, and $a,b$ is always a set of characteristic generators of $w(G)$.

\begin{thm}
\label{property-3.1}Let $G$ be a finite $p$-group such that
$G_3\le \Phi(G')\cong C_p$ and $G/\Phi(G')\cong M_p(n,m,1)$ where
$n>1$ for $p=2$ and $n\ge m$. Then the following conclusions hold:

{\rm (1)} If $m=1$, then $p=2$, $[c,b]=c^2$ and $G$ has a unique $\mathcal{A}_1$-subgroup of index $p$;

{\rm (2)} $I_{\min}=2$ for $m\ge 2$.
\end{thm}
\demo (1) If $m=1$, since $1=[a,b^p]=c^p[c,b]^{p\choose 2}$ and $c^p\neq 1$, then $p=2$ and $c^2=[c,b]$. Since $[c,ab][c,a]=c^2$, we have $[c,a]=c^2$ or $[c,ab]=c^2$. Without loss of generality, we may assume that $[c,a]=c^2$. Hence $[c,ab]=1$. Notice that the maximal subgroups of $G$ are $\lg a,\Phi(G)\rg=\lg c,a\rg$, $\lg b,\Phi(G)\rg=\lg c,b,a^2\rg$ and $\lg ab,\Phi(G)\rg=\lg c,ab\rg$. Since $\lg c,a\rg\in\mathcal{A}_1$, $\lg c,b,a^2\rg\not\in\mathcal{A}_1$ and $\lg c,ab\rg$ is abelian, $\lg c,a\rg$ is the unique $\mathcal{A}_1$-subgroup of index $p$.

(2) If $m\ge 2$, then we claim that $I_{\min}>1$. Otherwise, we may assume that $D$ is an $\mathcal{A}_1$-subgroup of index $p$. It is easy to see that $D'=\Phi(G')$ and $d(D/\Phi(G'))=2$. Since $\Phi(G)\le D$, $d(\Phi(G)/\Phi(G'))\le 2$. On the other hand, $\Phi(G)/\Phi(G')=\lg \bar{a}^p,\bar{b}^p,\bar{c}\rg$ is of type $(p^{n-1},p^{m-1},p)$, a contradiction. Hence $I_{\min}>1$.

Since $|G_3|\le p$, without loss of generality, we may assume that $[c,a]=1$ or $[c,b]=1$. If $[c,a]=1$, then $\lg a^p,b\rg$ is an $\mathcal{A}_1$-subgroup of index $p^2$. If $[c,b]=1$, then $\lg a,b^p\rg$ is an $\mathcal{A}_1$-subgroup of index $p^2$. To sum up, $I_{\min}=2$.\qed

\medskip

\begin{thm}
\label{property-3.1/2}Suppose that $G$ is a finite $p$-group such that
$G_3\le \Phi(G')\cong C_p$ and $G/\Phi(G')\cong M_p(n,m,1)$ where
$n>1$ for $p=2$ and $n\ge m$. Let $w(G)=(w_{ij})$ be a characteristic matrix of $G$. Then the following conclusions hold:

{\rm (1)} If $w_{22}=w_{12}=0$, then $G\in \mathcal{A}_{3}$;

{\rm (2)} If $w_{22}=0$ and $w_{12}\neq 0$, then $G\in \mathcal{A}_{m+1}$;

{\rm (3)} If $w_{22}\neq 0$ and $w_{12}=0$, then $G\in \mathcal{A}_{n+1}$.
\end{thm}
\demo We get $I_{\max}$ by investigating all maximal subgroups of $G$. Let $N=\lg b,a^p,c\rg$ and $M_i=\lg ab^i, b^p,c\rg$. Then all maximal subgroups of $G$ are $N$ and $M_i$ where $0\le i\le p-1$. Note that $m>1$ for $p>2$ by Theorem \ref{property-3.1} (1).

(1) In this case, $[c,a]=[c,b]=1$. By calculation, $M_i=\lg b^p,ab^i\rg\ast \lg c\rg$ and $N=\lg b,a^p\rg\ast \lg c\rg$ where $\lg b^p,ab^i\rg\in\mathcal{A}_1$ and $\lg b,a^p\rg\in\mathcal{A}_1$. By Lemma \ref{cor 2.4}, $M_i\in\mathcal{A}_2$ and $N\in\mathcal{A}_2$. Hence $G\in \mathcal{A}_{3}$.

(2) In this case, $[c,b]=1$. By Theorem \ref{property-3.1} (1), $m>1$. By calculation, $M_i=\lg c,ab^i\rg\ast\lg cb^{w_{12}p}\rg$ where $\lg c, ab^i\rg\in \mathcal{A}_1$. By Lemma \ref{cor 2.4}, $M_i\in\mathcal{A}_m$. If $p=2$, then $N=\lg c,b,a^2\rg$ is abelian. If $p>2$, then $N=\lg b,a^p\rg\ast \lg c\rg\in \mathcal{A}_2$. To sum up, $G\in \mathcal{A}_{m+1}$.

(3) In this case, $[c,b]\neq 1$ and $[c,a]=1$. By calculation, $N=\lg c,b\rg\ast \lg ca^{-w_{22}p}\rg$ where $\lg c,b\rg\in \mathcal{A}_1$. By Lemma \ref{cor 2.4}, $N\in \mathcal{A}_{n}$. If $p=2$, then $M_0=\lg c,a, b^2\rg$ is abelian and $M_1=\lg c,ab\rg\ast\lg b^2\rg\in \mathcal{A}_{m}$. If $p>2$, then $M_0=\lg c,a, b^p\rg=\lg a,b^p\rg\ast \lg c\rg\in \mathcal{A}_2$ and $M_i=\lg c,ab^i\rg\ast\lg cb^{iw_{22}p}\rg\in \mathcal{A}_{m}$  for $i=1,2,\dots,p-1$. To sum up, $G\in\mathcal{A}_{n+1}$.
 \qed
\begin{thm}\label{isomorphic-3.1}
Suppose that $G$ and $\bar{G}$ are finite $p$-groups such that $G_3\le \Phi(G')\cong C_p$ and $G/\Phi(G')\cong M_p(n,m,1)$, where $p>2$ and $n\ge m\ge 2$. Let two characteristic matrices of $G$ and $\bar{G}$ be $w(G)=(w_{ij})$ and $w(\bar{G})=(\bar{w}_{ij})$ respectively.
Then $G\cong \bar{G}$ if and only if there exists $X={\left(
 \begin{array}{cc}
 x_{11}& x_{12}\\
 x_{21}p^{n-m}& x_{22}
 \end{array}
 \right)}$, an invertible matrix over $F_p$, such that $$
\left(
 \begin{array}{c}
 \bar{w}_{11}\\
 \bar{w}_{21}
  \end{array}
 \right)=
|X|^{-1}
{\left(
 \begin{array}{cc}
 x_{11}& x_{12}p^{n-m}\\
 x_{21}& x_{22}
 \end{array}
 \right)}\left(
 \begin{array}{c}
 w_{11}\\
 w_{21}
 \end{array}
 \right)$$ and $$
\left(
 \begin{array}{c}
 \bar{w}_{12}\\
 \bar{w}_{22}
  \end{array}
 \right)=X\left(
 \begin{array}{cc}
 w_{12}\\
  w_{22}
 \end{array}
 \right).$$
 \end{thm}
\demo Suppose that $a,b$ and $\bar{a},\bar{b}$ are two set of characteristic generators of $w(G)$ and $w(\bar{G})$
respectively.
Let $\theta$ be an isomorphism from $\bar{G}$ onto $G$. We have $\Phi(\bar{G})^\theta=\Phi(G)$ and $\Omega_m(\bar{G})^\theta=\Omega_m(G)$ since these four groups are characteristic in $G$ or $\bar{G}$. So we may let
$$\bar{a}^\theta= a^{x_{11}}b^{x_{12}}\phi_1,\
\bar{b}^\theta=a^{x_{21}p^{n-m}}b^{x_{22}}\phi_2$$
where $\phi_1\in \Phi(G)$, $\phi_2\in \Phi(G)\cap \Omega_m(G)$ and $X:=\left(
 \begin{array}{cc}
 x_{11}& x_{12}\\
 x_{21}p^{n-m}& x_{22}
 \end{array}
 \right)$ is an invertible matrix over $F_p$.
By calculations,
$$\bar{c}^\theta=[\bar{a},\bar{b}]^\theta=[\bar{a}^\theta,\bar{b}^\theta]\equiv [a^{x_{11}}b^{x_{12}},a^{x_{21}p^{n-m}}b^{x_{22}}]\equiv
c^{|X|}\ (\mod \Phi(G')).$$
Transforming $\bar{c}^{\bar{w}_{11}p}=\bar{a}^{p^{n}}$ by $\theta$, we have 
$c^{|X|\bar{w}_{11}p}=a^{x_{11}p^n}b^{x_{12}p^n}$. Hence
\begin{equation}\label{eq 3.1}
|X|\bar{w}_{11}=
 (x_{11},x_{12}p^{n-m})
 \left(
 \begin{array}{c}
 w_{11}\\
 w_{21}
 \end{array}
 \right)
\end{equation}
Transforming $\bar{c}^{\bar{w}_{21}p}=\bar{b}^{p^{m}}$ by $\theta$, we have
 $c^{|X|\bar{w}_{21}p}=a^{x_{21}p^n}b^{x_{22}p^m}$. Hence
\begin{equation}\label{eq 3.2}
|X|\bar{w}_{21}=
 (x_{21},x_{22})
 \left(
 \begin{array}{c}
 w_{11}\\
 w_{21}
 \end{array}
 \right)
\end{equation}
By Equation \ref{eq 3.1} and \ref{eq 3.2},
\begin{equation}\label{eq 3.3}
\left(
 \begin{array}{c}
 \bar{w}_{11}\\
 \bar{w}_{21}
  \end{array}
 \right)=
 |X|^{-1}\left(
 \begin{array}{cc}
 x_{11}& x_{12}p^{n-m}\\
 x_{21}& x_{22}
 \end{array}
 \right)
 \left(
 \begin{array}{c}
 w_{11}\\
 w_{21}
 \end{array}
 \right)
\end{equation}
 Transforming $\bar{c}^{\bar{w}_{12}p}=[\bar{c},\bar{a}]$ by $\theta$, we have $c^{|X|\bar{w}_{12} p}=[c^{|X|},a^{x_{11}}b^{x_{12}}]$. It follows that $c^{\bar{w}_{12}p}=[c,a]^{x_{11}}[c,b]^{x_{12}}=c^{w_{12}x_{11}p+w_{22}x_{12}p}$. Hence
\begin{equation}\label{eq 3.4}
\bar{w}_{12}=
 (x_{11},x_{12})
 \left(
 \begin{array}{c}
 w_{12}\\
 w_{22}
 \end{array}
 \right)
\end{equation}
 Transforming $\bar{c}^{\bar{w}_{22}p}=[\bar{c},\bar{b}]$ by $\theta$, we have $c^{|X|\bar{w}_{22} p}=[c^{|X|},a^{x_{21}p^{n-m}}b^{x_{22}}]$. It follows that $c^{\bar{w}_{22}p}=[c,a]^{x_{21}p^{n-m}}[c,b]^{x_{22}}=c^{w_{12}x_{21}p^{n-m+1}+w_{22}x_{22}p}$. Hence
\begin{equation}\label{eq 3.5}
\bar{w}_{22}=
 (x_{21}p^{n-m},x_{22})
 \left(
 \begin{array}{c}
 w_{12}\\
 w_{22}
 \end{array}
 \right)
\end{equation}
By Equation \ref{eq 3.4} and \ref{eq 3.5},
\begin{equation}\label{eq 3.6}
\left(
 \begin{array}{c}
 \bar{w}_{12}\\
 \bar{w}_{22}
  \end{array}
 \right)=
 \left(
 \begin{array}{cc}
 x_{11}& x_{12}\\
 x_{21}p^{n-m}& x_{22}
 \end{array}
 \right)
 \left(
 \begin{array}{c}
 w_{12}\\
 w_{22}
 \end{array}
 \right)
\end{equation}

Conversely, if there exists an invertible matrix $X=\left(
 \begin{array}{cc}
 x_{11}& x_{12}\\
 x_{21}p^{n-m}& x_{22}
 \end{array}
 \right)$ over
$F_p$ such that Equation \ref{eq 3.3} and \ref{eq 3.6}, then, by using the above argument, it is easy to check that the map $\theta: \bar{a}\mapsto a^{x_{11}}b^{x_{12}},\bar{b}\mapsto a^{x_{21}p^{n-m}}b^{x_{22}}$ is an isomorphism from
$\bar{G}$ onto $G$.
 \qed

\medskip

 If $p=2$, $m\ge 2$ and $n\ge 3$, then we also have Equation \ref{eq 3.3} and \ref{eq 3.6}.
If $p=2$, $m=1$ and $n\ge 3$, then $b^2\in Z(G)$. It follows that $1=[a,b^2]=c^2[c,b]$. Hence $[c,b]=c^2$. That is, $w_{22}=1$. In this case, we may let $\bar{b}^\theta=a^{x_{21}2^{n-1}}bc^{x_{23}}$. By calculations, $(\bar{b}^{p^{m}})^\theta=(\bar{b}^{2})^\theta=(a^{x_{21}2^{n-1}}bc^{x_{23}})^{2}=a^{x_{21}2^n}b^{2}c^{2x_{23}}[c,b]^{x_{23}}=a^{x_{21}2^n}b^{2}$. Hence we also have Equation \ref{eq 3.2}. It follows that Equation \ref{eq 3.3} and \ref{eq 3.6} hold.

Therefore we have the following theorem.

\begin{thm}\label{isomorphic-3.2}
Suppose that $G$ and $\bar{G}$ are finite $2$-groups such that $G_3\le \Phi(G')\cong C_2$ and $G/\Phi(G')\cong M_2(n,m,1)$, where $n\ge 3$. Let two characteristic matrices of $G$ and $\bar{G}$ be $w(G)=(w_{ij})$ and $w(\bar{G})=(\bar{w}_{ij})$ respectively.
Then $G\cong \bar{G}$ if and only if there exists $X={\left(
 \begin{array}{cc}
 x_{11}& x_{12}\\
 x_{21}2^{n-m}& x_{22}
 \end{array}
 \right)}$, an invertible matrix over $F_2$, such that $$
\left(
 \begin{array}{c}
 \bar{w}_{11}\\
 \bar{w}_{21}
  \end{array}
 \right)=
{\left(
 \begin{array}{cc}
 x_{11}& x_{12}2^{n-m}\\
 x_{21}& x_{22}
 \end{array}
 \right)}\left(
 \begin{array}{c}
 w_{11}\\
 w_{21}
 \end{array}
 \right)$$ and $$
\left(
 \begin{array}{c}
 \bar{w}_{12}\\
 \bar{w}_{22}
  \end{array}
 \right)=X\left(
 \begin{array}{cc}
 w_{12}\\
  w_{22}
 \end{array}
 \right).$$ In addition, if $m=1$, then $w_{22}=1$.
 \end{thm}

\begin{thm}
\label{c_{p^2}}Let $G$ be a finite $p$-group such that
$G_3\le \Phi(G')\cong C_p$ and $G/\Phi(G')\cong M_p(n,m,1)$ where
$n>1$ for $p=2$ and $n\ge m$. Then $G$ is one of the following non-isomorphic groups:

{\rm (A1)}  $\langle a, b, c \di a^{4}=b^2=c^4=1,
[a,b]=c,[c,a]=[c,b]=c^2 \rangle$; %(C1)

{\rm (A2)} $\langle a, b, c \di a^{4}=b^{4}=1, c^2=a^{2},
[a,b]=c,[c,a]=[c,b]=c^2 \rangle$; %(C2)

{\rm (A3)}  $\langle a, b, c \di a^{8}=b^2=1,c^{2}=a^4,
[a,b]=c,[c,a]=[c,b]=c^2 \rangle$; %(C3)

{\rm (B1)}  $\langle a, b, c \di a^{2^{n+1}}=b^{2}=1, c^2=a^{2^n},
[a,b]=c,[c,a]=1,[c,b]=c^2 \rangle$, where $n\ge 3$; %(D1)

{\rm (B2)}  $\langle a, b, c \di a^{2^{n}}=b^{2}=c^4=1,
[a,b]=c,[c,a]=1,[c,b]=c^2 \rangle$, where $n\ge 3$; %(D2)

{\rm (B3)}  $\langle a, b, c \di a^{2^{n}}=b^{4}=1, c^2=b^{2},
[a,b]=c,[c,a]=1,[c,b]=c^2 \rangle$, where $n\ge 3$; %(D3)

{\rm (C1)}  $\langle a, b, c \di a^{8}=1, c^2=a^{4}=b^4,
[a,b]=c,[c,a]=1,[c,b]=1 \rangle$; %(D4)

{\rm (C2)} $\langle a, b, c \di a^{8}=b^{4}=1, c^2=a^{4},
[a,b]=c,[c,a]=[c,b]=1 \rangle$; %(D5) for $p=n=m=2$

{\rm (C3)}  $\langle a, b, c \di a^{8}=1, c^2=a^{4}=b^4,
[a,b]=c,[c,a]=1,[c,b]=c^2 \rangle$; %(D6)

{\rm (C4)} $\langle a, b, c \di a^{8}=b^{4}=1, c^2=a^{4},
[a,b]=c,[c,a]=1,[c,b]=c^{2} \rangle$; %(D7) for $p=n=m=2$

{\rm (C5)} $\langle a, b, c \di a^{8}=b^{4}=1, c^2=a^{4},
[a,b]=c,[c,a]=c^2,[c,b]=1 \rangle$;%(D8) for $p=n=m=2$

{\rm (D1)} $\langle a, b, c \di a^{p^{n+1}}=b^{p^n}=1, c^p=a^{p^n},
[a,b]=c,[c,a]=1,[c,b]=c^{tp} \rangle$, where
\mbox{\hskip0.65in}$n\ge 3$ for $p=2$, $n\ge 2$ and $t\in F_p^*$;%(D7) for $p>2$ and $n=m\ge 2$

{\rm (D2)} $\langle a, b, c \di a^{p^{n+1}}=b^{p^n}=1, c^p=a^{p^n},
[a,b]=c,[c,a]=c^p,[c,b]=1 \rangle$, where $n\ge 3$ \mbox{\hskip0.65in}for $p=2$ and $n\ge 2$;%(D8) for $p>2$ and $n=m\ge 2$

{\rm (D3)} $\langle a, b, c \di a^{p^{n+1}}=b^{p^n}=1, c^p=a^{p^n},
[a,b]=c,[c,a]=1,[c,b]=1 \rangle$, where $n\ge 3$ \mbox{\hskip0.65in}for $p=2$ and $n\ge 2$; %(D5) for $p>2$ and $n=m\ge 2$

{\rm (D4)} $\langle a, b, c \di a^{p^{n}}=b^{p^{n}}=c^{p^2}=1,
[a,b]=c,[c,a]=c^p,[c,b]=1 \rangle$, where $n\ge 3$ \mbox{\hskip0.65in}for $p=2$ and $n\ge 2$;
%(D13) for $p>2$ and $n=m\ge 2$

{\rm (D5)} $\langle a, b, c \di a^{p^{n}}=b^{p^{n}}=c^{p^2}=1,
[a,b]=c,[c,a]=1,[c,b]=1 \rangle$, where $n\ge 3$ \mbox{\hskip0.65in}for $p=2$ and $n\ge 2$;
%(D12) for $p>2$ and $n=m\ge 2$

{\rm (E1)} $\langle a, b, c \di a^{p^{n+1}}=b^{p^m}=1, c^p=a^{p^n},
[a,b]=c,[c,a]=1,[c,b]=c^{tp} \rangle$, where
\mbox{\hskip0.65in}$n>m\ge 2$ and $t\in F_p^*$;%(D7) for $n>m\ge 2$

{\rm (E2)} $\langle a, b, c \di a^{p^{n+1}}=b^{p^m}=1, c^p=a^{p^n},
[a,b]=c,[c,a]=c^p,[c,b]=1 \rangle$, where \mbox{\hskip0.65in}$n>m\ge 2$;%(D8) for $n>m\ge 2$

{\rm (E3)} $\langle a, b, c \di a^{p^{n+1}}=b^{p^m}=1, c^p=a^{p^n},
[a,b]=c,[c,a]=1,[c,b]=1 \rangle$, where \mbox{\hskip0.65in}$n>m\ge 2$; %(D5) for $n>m\ge 2$

{\rm (E4)} $\langle a, b, c \di a^{p^{n}}=b^{p^{m+1}}=1,
c^p=b^{p^m}, [a,b]=c,[c,a]=1,[c,b]=c^{p} \rangle$, where
\mbox{\hskip0.65in}$n>m\ge 2$; %(D10)

{\rm (E5)} $\langle a, b, c \di a^{p^{n}}=b^{p^{m+1}}=1,
c^p=b^{p^m}, [a,b]=c,[c,a]=c^{tp},[c,b]=1 \rangle$, where
\mbox{\hskip0.65in}$n>m\ge 2$ and $t\in F_p^*$; %(D11)

{\rm (E6)} $\langle a, b, c \di a^{p^{n}}=b^{p^{m+1}}=1,
c^p=b^{p^m}, [a,b]=c,[c,a]=1,[c,b]=1 \rangle$, where \mbox{\hskip0.65in}$n>m\ge 2$; %(D9)

{\rm (E7)} $\langle a, b, c \di a^{p^{n}}=b^{p^{m}}=c^{p^2}=1,
[a,b]=c,[c,a]=1,[c,b]=c^p \rangle$, where $n>m\ge 2$;%(D13) for $n>m\ge 2$

{\rm (E8)} $\langle a, b, c \di a^{p^{n}}=b^{p^{m}}=c^{p^2}=1,
[a,b]=c,[c,a]=c^p,[c,b]=1 \rangle$, where $n>m\ge 2$; %(D14)

{\rm (E9)} $\langle a, b, c \di a^{p^{n}}=b^{p^{m}}=c^{p^2}=1,
[a,b]=c,[c,a]=1,[c,b]=1 \rangle$, where $n>m\ge 2$.%(D12) for $n>m\ge 2$
\end{thm}
\demo
{\bf Case 1:} $m=1$.

If $p>2$, then $c^p=[a,b^p]=1$, a contradiction. Hence $p=2$. It follows that $n\ge 2$. If
$n=2$, then $|G|=2^5$. By checking the list of groups of order
$2^5$, we get the groups of Type (A1)--(A3). In the following, we may assume that $n\ge 3$. Suppose that $G$ and $\bar{G}$ are two groups described in the theorem. By Theorem \ref{isomorphic-3.2}, $\bar{w}_{22}=w_{22}=1$, and $\bar{G}\cong G$ if and only if there exists
$X={\left(
 \begin{array}{cc}
 1& x_{12}\\
 x_{21}2^{n-1}& 1
 \end{array}
 \right)}$, an invertible matrix over $F_2$, such that $$
\left(
 \begin{array}{c}
 \bar{w}_{11}\\
 \bar{w}_{21}
  \end{array}
 \right)=
{\left(
 \begin{array}{cc}
 1& 0\\
 x_{21}& 1
 \end{array}
 \right)}\left(
 \begin{array}{c}
 w_{11}\\
 w_{21}
 \end{array}
 \right)$$ and $$
\left(
 \begin{array}{c}
 \bar{w}_{12}\\
 1
  \end{array}
 \right)={\left(
 \begin{array}{cc}
 1& x_{12}\\
 0& 1
 \end{array}
 \right)}\left(
 \begin{array}{cc}
 w_{12}\\
  1
 \end{array}
 \right).$$ Let $x_{12}=w_{12}$. Then we have $\bar{w}_{12}=0$. If $w_{11}=0$, then $\bar{w}_{11}=0$ and $\bar{w}_{21}={w}_{21}$. If $w_{11}=1$, then, letting $x_{21}=w_{21}$, we have $\bar{w}_{21}=0$. To sum up, characteristic matrices of non-isomorphic groups are:
\begin{center}
(b1)  $\left(
 \begin{array}{cc}
 1& 0\\
 0&1
 \end{array}
 \right)$\ \
(b2)  $\left(
 \begin{array}{cc}
 0& 0\\
 0&1
 \end{array}
 \right)$\ \
(b3)  $\left(
 \begin{array}{cc}
 0& 0\\
 1&1
 \end{array}
 \right).$
\end{center}

\noindent Hence we get the groups of Type (B1)--(B3).

\medskip

{\bf Case 2:} $p=n=m=2$.

In this case, $|G|=2^6$. By checking the list of groups of order
$2^6$, we get the groups of Type (C1)--(C5).

\medskip

 {\bf Case 3:} $n\ge 3$ for $p=2$ and $n=m\ge 2$.

Suppose that $G$ and $\bar{G}$ are two groups described in the theorem. By Theorem \ref{isomorphic-3.1} and \ref{isomorphic-3.2},
$G\cong \bar{G}$ if and only if there exists $X={\left(
 \begin{array}{cc}
 x_{11}& x_{12}\\
 x_{21}& x_{22}
 \end{array}
 \right)}$, an invertible matrix over $F_p$, such that
$$\left(
 \begin{array}{c}
 \bar{w}_{11}\\
 \bar{w}_{21}
  \end{array}
 \right)=
{|X|^{-1}}
X\left(
 \begin{array}{c}
 w_{11}\\
 w_{21}
 \end{array}
 \right)$$ and
 $$\left(
 \begin{array}{c}
 \bar{w}_{12}\\
 \bar{w}_{22}
  \end{array}
 \right)=X\left(
 \begin{array}{cc}
 w_{12}\\
  w_{22}
 \end{array}
 \right).$$
That is,
\begin{equation}\label{eq 3.9}
\left(
 \begin{array}{cc}
 \bar{w}_{11}&\bar{w}_{12}\\
 \bar{w}_{21}&\bar{w}_{22}
  \end{array}
 \right)=
X\left(
 \begin{array}{cc}
 w_{11}&w_{12}\\
 w_{21}&w_{22}
 \end{array}
 \right)
\left(
 \begin{array}{cc}
 {|X|^{-1}}&0\\
 0&1
 \end{array}
 \right)
  \end{equation}

If $w(G)$ is invertible, then,
letting $X=\diag(1,|w(G)|)(w(G))^{-1}$, we have $w(\bar{G})=\diag(1,|w(G)|)$. Hence we get the group of Type (D1) where $t=|w(G)|$. By
Equation \ref{eq 3.9}, $|w(\bar{G})|=|w(G)|$. Hence different $t$ give non-isomorphic groups.

If $w(G)=0$, then $G$ is the group of Type (D5). In the following, we assume that $w(G)$ is of rank 1.

By suitably choosing $X$ such that $|X|=1$, that is, using some row operations, we can simplify $w(G)$ to be $\left(
 \begin{array}{cc}
 w_{11}&w_{12}\\
 0&0
 \end{array}
 \right)$
 where $(w_{11},w_{12})\neq (0,0)$.
Let $w(\bar{G})$ and $w(G)$ be such matrices. Then, by Equation \ref{eq 3.9}, it is easy to check that
$G\cong \bar{G}$ if and only if there exists $X=\diag(x_{11},x_{22})$, an invertible matrix over $F_p$, such that
$w(\bar{G})=Xw(G)\diag(|X|^{-1},1)$.\begin{table}[h]
  \centering
\small
%\topcaption{}
%Properties of the groups of Type (A) in Theorem
%\ref{main}}

{
\begin{tabular}{cccc}
Case &$X$ & $w(\bar{G})$&Group\\
   \hline
  $w_{11}\neq 0$ and $w_{12}\neq 0$&$\diag(w_{12}^{-1},w_{11})$&$\left(
 \begin{array}{cc}
 1& 1\\
 0& 0
 \end{array}
 \right)$& (D2)\\
 % \hline
  $w_{11}\neq 0$ and $w_{12}=0$& $\diag(1,w_{11})$&$\left(
 \begin{array}{cc}
 1& 0\\
 0& 0
 \end{array}
 \right)$& (D3)\\
 %\hline
  $w_{11}=0$ and $w_{12}\neq 0$& $\diag(w_{12}^{-1},w_{12})$&$\left(
 \begin{array}{cc}
0& 1\\
 0& 0
 \end{array}
 \right)$& (D4)\\

   \hline
\end{tabular}
\caption{Case 3 in Theorem \ref{c_{p^2}} where $w(G)$ is of rank 1}
\label{Table 1}
}\end{table}
By Table \ref{Table 1}, we get the groups of (D2)--(D4).
It is easy to check that different types give non-isomorphic groups.

{\bf Case 4:} $n>m\ge 2$.

\medskip

Suppose that $G$ and $\bar{G}$ are two groups described in the theorem. By Theorem \ref{isomorphic-3.1} and \ref{isomorphic-3.2},
$G\cong \bar{G}$ if and only if there exists $X={\left(
 \begin{array}{cc}
 x_{11}& x_{12}\\
 x_{21}p^{n-m}& x_{22}
 \end{array}
 \right)}$, an invertible matrix over $F_p$, such that
 \begin{equation}\label{eq: 3.10}
\left(
 \begin{array}{c}
 \bar{w}_{11}\\
 \bar{w}_{21}
  \end{array}
 \right)=
 |X|^{-1}
{\left(
 \begin{array}{cc}
 x_{11}& 0\\
 x_{21}& x_{22}
 \end{array}
 \right)}\left(
 \begin{array}{c}
 w_{11}\\
 w_{21}
 \end{array}
 \right)
 \end{equation} and
 \begin{equation}\label{eq: 3.11}
\left(
 \begin{array}{c}
 \bar{w}_{12}\\
 \bar{w}_{22}
  \end{array}
 \right)={\left(
 \begin{array}{cc}
 x_{11}& x_{12}\\
 0 & x_{22}
 \end{array}
 \right)}\left(
 \begin{array}{cc}
 w_{12}\\
  w_{22}
 \end{array}
 \right)
 \end{equation}
 By suitably choosing $x_{21}$, that is, using an elementary row operation, we can simplify $\left(
 \begin{array}{cc}
 w_{11}\\
  w_{21}
 \end{array}
 \right)$ to be one of the following three types:

\begin{center}(a1)  $\left(
 \begin{array}{c}
 w_{11}\\
 0
 \end{array}
 \right)$ where $w_{11}\neq 0$,
(b1)  $\left(
 \begin{array}{c}
 0\\
 w_{21}
 \end{array}
 \right)$ where $w_{21}\neq 0$,
(c1)  $\left(
 \begin{array}{c}
 0\\
 0
 \end{array}
 \right).$
\end{center}
In the following, we assume that both $\left(
 \begin{array}{cc}
 w_{11}\\
  w_{21}
 \end{array}
 \right)$ and $\left(
 \begin{array}{cc}
 \bar{w}_{11}\\
  \bar{w}_{21}
 \end{array}
 \right)$ are such matrices. By Equation \ref{eq: 3.10} and \ref{eq: 3.11}, it is easy to check that (i) different types give non-isomorphic groups, and (ii)
$G\cong \bar{G}$ if and only if there exists $X={\left(
 \begin{array}{cc}
 x_{11}& x_{12}\\
 0 & x_{22}
 \end{array}
 \right)}$, an invertible matrix over $F_p$, such that
 \begin{equation}\label{eq: 3.10'}
\left(
 \begin{array}{c}
 \bar{w}_{11}\\
 \bar{w}_{21}
  \end{array}
 \right)=
 |X|^{-1}
{\left(
 \begin{array}{cc}
 x_{11}& 0\\
 0& x_{22}
 \end{array}
 \right)}\left(
 \begin{array}{c}
 w_{11}\\
 w_{21}
 \end{array}
 \right)
 \end{equation} and
 \begin{equation}\label{eq: 3.11'}
\left(
 \begin{array}{c}
 \bar{w}_{12}\\
 \bar{w}_{22}
  \end{array}
 \right)={\left(
 \begin{array}{cc}
 x_{11}& x_{12}\\
 0 & x_{22}
 \end{array}
 \right)}\left(
 \begin{array}{cc}
 w_{12}\\
  w_{22}
 \end{array}
 \right).
 \end{equation}
 By suitably choosing $x_{12}$, that is, using an elementary row operation, we can simplify $\left(
 \begin{array}{cc}
 w_{12}\\
  w_{22}
 \end{array}
 \right)$ to be one of the following three types:
\begin{center}
(a2)  $\left(
 \begin{array}{c}
 0\\
 w_{22}
 \end{array}
 \right)$ where $w_{22}\neq 0$,
(b2)  $\left(
 \begin{array}{c}
 w_{12}\\
 0
 \end{array}
 \right)$ where $w_{21}\neq 0$,
(c2)  $\left(
 \begin{array}{c}
 0\\
 0
 \end{array}
 \right).$
\end{center}
In the following, we assume that both $\left(
 \begin{array}{cc}
 w_{12}\\
  w_{22}
 \end{array}
 \right)$ and $\left(
 \begin{array}{cc}
 \bar{w}_{12}\\
  \bar{w}_{22}
 \end{array}
 \right)$ are such matrices. By Equation \ref{eq: 3.10'} and \ref{eq: 3.11'}, it is easy to check that (i) different types give non-isomorphic groups, and (ii)
$G\cong \bar{G}$ if and only if there exists $X=\diag(x_{11},x_{22})$, an invertible matrix over $F_p$, such that
 \begin{equation}\label{eq: 3.10''}
\left(
 \begin{array}{c}
 \bar{w}_{11}\\
 \bar{w}_{21}
  \end{array}
 \right)=
 |X|^{-1}
X\left(
 \begin{array}{c}
 w_{11}\\
 w_{21}
 \end{array}
 \right)
 \ \rm{ and } \
\left(
 \begin{array}{c}
 \bar{w}_{12}\\
 \bar{w}_{22}
  \end{array}
 \right)=X\left(
 \begin{array}{cc}
 w_{12}\\
  w_{22}
 \end{array}
 \right).
 \end{equation}
By Equation \ref{eq: 3.10''}, $G\cong \bar{G}$ if and only if there exists $X=\diag(x_{11},x_{22})$, an invertible matrix over $F_p$, such that $w(\bar{G})=Xw(G)\diag(|X|^{-1},1)$.
\begin{table}[h]
  \centering
{\scriptsize
\begin{tabular}{ccccccc}
\hline
$\left(
 \begin{array}{cc}
 w_{11}\\
  w_{21}
 \end{array}
 \right)$ &$\left(
 \begin{array}{cc}
 w_{12}\\
  w_{22}
 \end{array}
 \right)$ &$w(G)$&$X$ & $w(\bar{G})$&Group & Remark\\
   \hline
(a1)&  (a2) & $\left(
 \begin{array}{cc}
 w_{11}& 0\\
 0& w_{22}
 \end{array}
 \right)$ &$\diag(1,w_{11})$&$\left(
 \begin{array}{cc}
 1& 0\\
 0& w_{11}w_{22}
 \end{array}
 \right)$& (E1) & $t=w_{11}w_{22}$\\
 \hline
(a1)&  (b2) &$\left(
 \begin{array}{cc}
 w_{11}& w_{12}\\
 0&0
 \end{array}
 \right)$ &$\diag(w_{12}^{-1},w_{11})$&$\left(
 \begin{array}{cc}
 1& 1\\
 0& 0
 \end{array}
 \right)$& (E2) & \\
\hline
(a1)&  (c2) &$\left(
 \begin{array}{cc}
 w_{11}&0\\
 0&0
 \end{array}
 \right)$

 &$\diag(1,w_{11})$&$\left(
 \begin{array}{cc}
 1& 0\\
 0& 0
 \end{array}
 \right)$& (E3) & \\
 \hline
(b1)&  (a2) &$\left(
 \begin{array}{cc}
 0& 0\\
 w_{21}& w_{22}
 \end{array}
 \right)$& $\diag(w_{21},w_{22}^{-1})$&$\left(
 \begin{array}{cc}
 0& 0\\
 1& 1
 \end{array}
 \right)$& (E4) & \\
 \hline
(b1)&  (b2) &$\left(
 \begin{array}{cc}
 0& w_{12}\\
 w_{21}&0
 \end{array}
 \right)$ &$\diag(w_{21},1)$&$\left(
 \begin{array}{cc}
 0& w_{12}w_{21}\\
 1& 0
 \end{array}
 \right)$& (E5) & $t=w_{12}w_{21}$\\
\hline
(b1)&  (c2)&$\left(
 \begin{array}{cc}
 0&0\\
 w_{21}&0
 \end{array}
 \right)$
 & $\diag(w_{21},1)$&$\left(
 \begin{array}{cc}
 0& 0\\
 1& 0
 \end{array}
 \right)$& (E6) & \\
\hline
(c1)&  (a2) &$\left(
 \begin{array}{cc}
 0& 0\\
 0& w_{22}
 \end{array}
 \right)$& $\diag(1,w_{22}^{-1})$&$\left(
 \begin{array}{cc}
 0& 0\\
 0& 1
 \end{array}
 \right)$& (E7) & \\
\hline
(c1)&  (b2) &$\left(
 \begin{array}{cc}
 0& w_{12}\\
 0&0
 \end{array}
 \right)$& $\diag(w_{12}^{-1},1)$&$\left(
 \begin{array}{cc}
 0& 1\\
 0& 0
 \end{array}
 \right)$& (E8) & \\
\hline
(c1)&  (c2)&$\left(
 \begin{array}{cc}
 0&0\\
 0&0
 \end{array}
 \right)$
 & &$\left(
 \begin{array}{cc}
 0& 0\\
 0& 0
 \end{array}
 \right)$& (E9) & \\
   \hline
\end{tabular}
\caption{Case 4 in Theorem \ref{c_{p^2}}}
\label{table 2}
}\end{table}
By Table \ref{table 2}, we get the groups of Type (E1)--(E9).\qed

\medskip

 By Theorem \ref{property-3.1/2}, we have the following theorem.

\begin{thm}
\label{property-3.2} Let $G$ be a finite $p$-group listed in Theorem $\ref{c_{p^2}}$. Then the following conclusions hold:

{\rm (1)} $I_{\max}=2$ for types {\rm (A1)--(A3), (C1)--(C5), (D3), (D5), (E3), (E6), (E9)};

{\rm (2)} $I_{\max}=n$  for types {\rm (B1)--(B3), (D1)--(D2), (D4), (E1), (E4), (E7)};

{\rm (3)} $I_{\max}=m$ for types {\rm (E2), (E5), (E8)}.
\end{thm}

\section{The case $G_3\cong C_p$ and $\Phi(G')=1$}

Suppose that $G$ is a finite $p$-group with $G_3\cong C_p$, $\Phi(G')=1$ and $G/G_3\cong M_p(n,m,1)$, where $n>1$ for $p=2$ and
$n\ge m$. Let $$G/G_3=\lg \bar{a},\bar{b},\bar{c}\mid \bar{a}^{p^n}=\bar{b}^{p^m}=\bar{c}^{p}=1,[\bar{c},\bar{a}]=[\bar{b},\bar{c}]=1\rg.$$
Then, without loss of generality, we may assume that $G=\lg a,b,c\rg$ such that $[a,b]=c$. Since $\Phi(G')=1$, $c^p=1$. Let $G_3=\lg z\rg$.
Since $a^{p^{n}}\in G_3$, we may assume that $a^{p^{n}}=z^{w_{11}}$. By similar reasons, we may assume that $b^{p^m}=z^{w_{21}}$, $[c,a]=z^{w_{12}}$ and $[c,b]=z^{w_{22}}$.
Then we get a $2\times 2$ matrix over $F_p$. $w(G)=(w_{ij})$ is called a characteristic matrix of $G$. Notice that $w(G)$ will be changed if we change the generators $a,b$. We also call $a,b$ a set of characteristic generators of $w(G)$.

\begin{thm}
\label{property-4.1}Let $G$ be a finite $p$-group such that
$G_3\cong C_p$, $\Phi(G')=1$ and $G/G_3\cong M_p(n,m,1)$, where $n>1$ for $p=2$ and
$n\ge m$. Then the following conclusion hold:

{\rm (1)} If $m=1$, then $I_{\min}=1$;

{\rm (2)} If $m\ge 2$, then $I_{\min}=2$.
\end{thm}
\demo (1) Since $G_3=\lg [c,ab],[c,a]\rg\cong C_p$, $[c,a]\neq 1$ or $[c,ab]\neq 1$. Hence either $\lg a,c\rg$ or $\lg ab,c\rg$ is an $\mathcal{A}_1$-subgroup of index $p$. Thus $I_{\min}=1$.

(2) Firstly, we claim that $I_{\min}>1$. Otherwise, $I_{\min}=1$. Assume that $D$ is an $\mathcal{A}_1$-subgroup of index $p$. It is easy to see that $D'=G_3$ and $d(D/G_3)=2$. Since $\Phi(G)\le D$, $d(\Phi(G)/G_3)\le 2$. On the other hand, $\Phi(G)/G_3=\lg \bar{a}^p,\bar{b}^p,\bar{c}\rg$ is of type $(p^{n-1},p^{m-1},p)$, a contradiction.

If $[c,b]=1$, then $[c,a]\neq 1$. In this case, $\lg cb^p,a\rg$ is an $\mathcal{A}_1$-subgroup of index $p^2$. If $[c,b]\neq 1$, then, without loss of generality, we may assume that $[c,a]=1$. In this case, $\lg b,ca^p\rg$ is an $\mathcal{A}_1$-subgroup of index $p^2$. To sum up, $I_{\min}=2$.\qed

\begin{thm}
\label{property-4.1/2} Suppose that $G$ is a finite $p$-group such that
$G_3\cong C_p$, $\Phi(G')=1$ and $G/G_3\cong M_p(n,m,1)$, where $n>1$ for $p=2$ and
$n\ge m$. Let $w(G)=(w_{ij})$ be a characteristic matrix of $G$. Then the following conclusions hold:

{\rm (1)} If $p=2$ and $m=1$, then $G\in \mathcal{A}_{3}$;

{\rm (2)} If $m>1$ for $p=2$, $w_{22}=0$ and $w_{12}\neq 0$, then $G\in \mathcal{A}_{m+1}$;

{\rm (3)} If $m>1$ for $p=2$, $w_{22}\neq 0$ and $w_{12}=0$, then $G\in \mathcal{A}_{n+1}$.
\end{thm}
\demo We get $I_{\max}$ by investigating all maximal subgroups of $G$. Let $G_3=\lg z\rg$, $N=\lg b,a^p,c,z\rg$ and $M_i=\lg ab^i, b^p,c,z\rg$. Then all maximal subgroups of $G$ are $N$ and $M_i$ where $0\le i\le p-1$.

(1) In this case, $[c,b]=[a,b^2]=1$. Hence $[c,a]\neq 1$. It follows that $M_i=\lg c,ab^i\rg\in \mathcal{A}_1$. Since $\lg a^2,b\rg\in\mathcal{A}_1$, we have $N=\lg a^2,b\rg\times \lg c\rg\in\mathcal{A}_2$ by Lemma \ref{cor 2.4}. Hence $G\in \mathcal{A}_{3}$.

(2) In this case, $[c,b]=1$ and $[c,a]\neq 1$. Hence $\lg c,ab^i\rg\in \mathcal{A}_1$. If $m=1$, then $M_i=\lg c,ab^i\rg\in \mathcal{A}_m$. If $m>1$, then $M_i=\lg c,ab^i\rg\ast\lg b^p\rg$. By Lemma \ref{cor 2.4}, we also have $M_i\in\mathcal{A}_m$. If $p=2$, then $m>1$ and $N=\lg a^2,b\rg\times \lg c\rg\in\mathcal{A}_2$. If $p>2$, then $N=\lg a^p,b, c,z\rg$ is abelian. To sum up, $G\in \mathcal{A}_{m+1}$.

(3) In this case, $[c,a]=1$ and $[c,b]\neq 1$. By Lemma \ref{cor 2.4}, $N=\lg c,b\rg\ast \lg a^{p}\rg\in\mathcal{A}_n$. If $p=2$, then $M_0=\lg a, b^2\rg\times \lg c\rg\in \mathcal{A}_2$ and $M_1=\lg c,ab\rg\ast\lg cb^2\rg\in \mathcal{A}_{m}$. If $p>2$, then $M_0=\lg c,a,b^p,z\rg$ is abelian and $M_i=\lg c,ab^i\rg\ast\lg b^{p}\rg\in \mathcal{A}_{m}$  for $i=1,2,\dots,p-1$. To sum up, $G\in\mathcal{A}_{n+1}$.
 \qed
\begin{thm}\label{isomorphic-4.1}
Suppose that $G$ and $\bar{G}$ are finite $p$-groups such that $G_3\cong C_p$, $\Phi(G')=1$ and $G/G_3\cong M_p(n,m,1)$, where  $n\ge m\ge 2$. Let two characteristic matrices of $G$ and $\bar{G}$ be $w(G)=(w_{ij})$ and $w(\bar{G})=(\bar{w}_{ij})$ respectively.
Then $G\cong \bar{G}$ if and only if there exists $X={\left(
 \begin{array}{cc}
 x_{11}& x_{12}\\
 x_{21}p^{n-m}& x_{22}
 \end{array}
 \right)}$, an invertible matrix over $F_p$, and $\lambda\in F_p^*$, such that $$
\left(
 \begin{array}{c}
 \bar{w}_{11}\\
 \bar{w}_{21}
  \end{array}
 \right)=
{\lambda}^{-1}
{\left(
 \begin{array}{cc}
 x_{11}& x_{12}p^{n-m}\\
 x_{21}& x_{22}
 \end{array}
 \right)}\left(
 \begin{array}{c}
 w_{11}\\
 w_{21}
 \end{array}
 \right)$$ and $$
\left(
 \begin{array}{c}
 \bar{w}_{12}\\
 \bar{w}_{22}
  \end{array}
 \right)={\lambda}^{-1}|X|X\left(
 \begin{array}{cc}
 w_{12}\\
  w_{22}
 \end{array}
 \right).$$
 \end{thm}
\demo Suppose that $a,b$ and
$\bar{a},\bar{b}$ are two set of characteristic generators of $w(G)$ and $w(\bar{G})$ respectively.
Let $\theta$ be an isomorphism from $\bar{G}$ onto $G$. 
We have $\Phi(\bar{G})^\theta=\Phi(G)$ and $\Omega_m(\bar{G})^\theta=\Omega_m(G)$ since these four groups are characteristic in $G$ or $\bar{G}$. So we may let
$$\bar{a}^\theta= a^{x_{11}}b^{x_{12}}\phi_1,\
\bar{b}^\theta=a^{x_{21}p^{n-m}}b^{x_{22}}\phi_2$$
where $\phi_1\in \Phi(G)$, $\phi_2\in \Phi(G)\cap \Omega_m(G)$ and $X:=\left(
 \begin{array}{cc}
 x_{11}& x_{12}\\
 x_{21}p^{n-m}& x_{22}
 \end{array}
 \right)$ is an invertible matrix over $F_p$.
By calculations,
$$\bar{c}^\theta=[\bar{a},\bar{b}]^\theta=[\bar{a}^\theta,\bar{b}^\theta]\equiv [a^{x_{11}}b^{x_{12}},a^{x_{21}p^{n-m}}b^{x_{22}}]\equiv
c^{|X|}\ (\mod G_3).$$ Let $\bar{z}^\theta=z^{\lambda}$. Then $\lambda\neq 0$.
Transforming $\bar{z}^{\bar{w}_{11}}=\bar{a}^{p^{n}}$ by $\theta$, we have $z^{\lambda \bar{w}_{11}}=a^{x_{11}p^n}b^{x_{12}p^n}$. Hence
\begin{equation}\label{eq 4.1}
\lambda\bar{w}_{11}=
 (x_{11},x_{12}p^{n-m})
 \left(
 \begin{array}{c}
 w_{11}\\
 w_{21}
 \end{array}
 \right)
\end{equation}
 Transforming $\bar{z}^{\bar{w}_{21}}=\bar{b}^{p^{m}}$ by $\theta$, we have $z^{\lambda \bar{w}_{21}}=a^{x_{21}p^n}b^{x_{22}p^m}$. Hence
\begin{equation}\label{eq 4.2}
\lambda\bar{w}_{21}=
 (x_{21},x_{22})
 \left(
 \begin{array}{c}
 w_{11}\\
 w_{21}
 \end{array}
 \right)
\end{equation}
By Equation \ref{eq 4.1} and \ref{eq 4.2},
\begin{equation}\label{eq 4.3}
\left(
 \begin{array}{c}
 \bar{w}_{11}\\
 \bar{w}_{21}
  \end{array}
 \right)=
 {\lambda}^{-1}\left(
 \begin{array}{cc}
 x_{11}& x_{12}p^{n-m}\\
 x_{21}& x_{22}
 \end{array}
 \right)
 \left(
 \begin{array}{c}
 w_{11}\\
 w_{21}
 \end{array}
 \right).
\end{equation}
Transforming $\bar{z}^{\bar{w}_{12}}=[\bar{c},\bar{a}]$ by $\theta$, we have $z^{\lambda \bar{w}_{12}}=[c^{|X|},a^{x_{11}}b^{x_{12}}]$. It follows that
\begin{equation}\label{eq 4.4}
\lambda\bar{w}_{12}=|X|
 (x_{11},x_{12})
 \left(
 \begin{array}{c}
 w_{12}\\
 w_{22}
 \end{array}
 \right).
\end{equation}
 Transforming $\bar{z}^{\bar{w}_{22}}=[\bar{c},\bar{b}]$ by $\theta$, we have $z^{\lambda \bar{w}_{22}}=[c^{|X|},a^{x_{21}p^{n-m}}b^{x_{22}}]$. It follows that
\begin{equation}\label{eq 4.5}
\lambda\bar{w}_{22}=
 |X|(x_{21}p^{n-m},x_{22})
 \left(
 \begin{array}{c}
 w_{12}\\
 w_{22}
 \end{array}
 \right).
\end{equation}
By Equation \ref{eq 4.4} and \ref{eq 4.5},
\begin{equation}\label{eq 4.6}
\left(
 \begin{array}{c}
 \bar{w}_{12}\\
 \bar{w}_{22}
  \end{array}
 \right)={\lambda}^{-1}|X|
 \left(
 \begin{array}{cc}
 x_{11}& x_{12}\\
 x_{21}p^{n-m}& x_{22}
 \end{array}
 \right)
 \left(
 \begin{array}{c}
 w_{12}\\
 w_{22}
 \end{array}
 \right).
\end{equation}

Conversely, if there exists an invertible matrix $X=\left(
 \begin{array}{cc}
 x_{11}& x_{12}\\
 x_{21}p^{n-m}& x_{22}
 \end{array}
 \right)$ over
$F_p$ and $\lambda\in F_p^*$ such that Equation \ref{eq 4.3} and \ref{eq 4.6}, then, by using the above argument, it is easy to check that the map $\theta: \bar{a}\mapsto a^{x_{11}}b^{x_{12}},\bar{b}\mapsto a^{x_{21}p^{n-m}}b^{x_{22}}$ is an isomorphism from
$\bar{G}$ onto $G$.
 \qed

\medskip

For odd $p$ and $m=1$, if $p>3$ or $n\ge 2$, then we also have
Equation \ref{eq 4.3} and \ref{eq 4.6}. Hence we get the following theorem.

\begin{thm}\label{isomorphic-4.2}
Suppose that $G$ and $\bar{G}$ are finite $p$-groups such that $G_3\cong C_p$, $\Phi(G')=1$ and $G/G_3\cong M_p(n,1,1)$, where $p>2$ and $n>1$ for $p=3$. Let two characteristic matrices of $G$ and $\bar{G}$ be $w(G)=(w_{ij})$ and $w(\bar{G})=(\bar{w}_{ij})$ respectively.
Then $G\cong \bar{G}$ if and only if there exists $X={\left(
 \begin{array}{cc}
 x_{11}& x_{12}\\
 x_{21}p^{n-m}& x_{22}
 \end{array}
 \right)}$, an invertible matrix over $F_p$, and $\lambda\in F_p^*$ such that $$
\left(
 \begin{array}{c}
 \bar{w}_{11}\\
 \bar{w}_{21}
  \end{array}
 \right)=
{\lambda}^{-1}
{\left(
 \begin{array}{cc}
 x_{11}& x_{12}p^{n-m}\\
 x_{21}& x_{22}
 \end{array}
 \right)}\left(
 \begin{array}{c}
 w_{11}\\
 w_{21}
 \end{array}
 \right)$$ and $$
\left(
 \begin{array}{c}
 \bar{w}_{12}\\
 \bar{w}_{22}
  \end{array}
 \right)={\lambda}^{-1}{|X|}X\left(
 \begin{array}{cc}
 w_{12}\\
  w_{22}
 \end{array}
 \right).$$
 \end{thm}
\qed
\medskip

If $p=2$ and $m=1$, then $1=[a,b^2]=[c,b]$. Since $G_3\neq 1$, $[c,a]\neq 1$. That is, $(w_{12},w_{22})=(1,0)$. In this case, we may let $\bar{b}^\theta=a^{x_{21}2^{n-1}}bc^{x_{23}}$. If $n=2$, then $(\bar{b}^{2})^\theta=(a^{x_{21}2}bc^{x_{23}})^{2}=a^{x_{21}2^2}b^{2}[c,a]^{x_{21}}$ and Equation \ref{eq 4.2} is changed to be
$$\bar{w}_{21}=
 (x_{21},1)
 \left(
 \begin{array}{c}
 w_{11}\\
 w_{21}
 \end{array}
 \right)+x_{21}.
\eqno(4.2')
$$
If $n\ge 3$, then we also have
Equation \ref{eq 4.3} and \ref{eq 4.6}. Hence we get the following theorem.

\begin{thm}\label{isomorphic-4.3}
Suppose that $G$ and $\bar{G}$ are finite $2$-groups such that $G_3\cong C_2$, $\Phi(G')=1$ and $G/G_3\cong M_2(n,1,1)$, where $n\ge 3$. Let two characteristic matrices of $G$ and $\bar{G}$ be $w(G)=(w_{ij})$ and $w(\bar{G})=(\bar{w}_{ij})$ respectively.
Then $G\cong \bar{G}$ if and only if there exists $X=
{\left(
 \begin{array}{cc}
 1& 0\\
 x_{21}& 1
 \end{array}
 \right)}$, a matrix over $F_2$ such that $
\left(
 \begin{array}{c}
 \bar{w}_{11}\\
 \bar{w}_{21}
  \end{array}
 \right)=X\left(
 \begin{array}{c}
 w_{11}\\
 w_{21}
 \end{array}
 \right)$.
 \end{thm}

\begin{thm}
\label{c_p^2}Let $G$ be a finite $p$-group such that
$G_3\cong C_p$, $\Phi(G')=1$ and $G/G_3\cong M_p(n,m,1)$, where $n>1$ for $p=2$ and
$n\ge m$. Then $G$ is one of the following non-isomorphic groups:

{\rm (F)} one of $3$-groups of maximal class of order $3^4$; %(A4)

{\rm (G1)}  $\langle a, b, c \di a^{p^{m+1}}=b^{p^m}=c^p=1,
[a,b]=c,[c,a]=1,[c,b]=a^{\nu p^m} \rangle$, where \mbox{\hskip0.6in}$m>1$ for $p\le 3$,
$\nu=1$ or a fixed quadratic non-residue modular
$p$; %(B2)

{\rm (G2)}  $\langle a, b, c,d \di a^{p^{m}}=b^{p^{m}}=c^p=d^p=1,
[a,b]=c,[c,a]=d,[c,b]=1,[d,a]=\mbox{\hskip0.6in}[d,b]=1 \rangle$,
where $m>1$ for $p\le 3$; %(B5)

{\rm (G3)}  $\langle a, b, c \di a^{p^{m+1}}=b^{p^m}=c^p=1,
[a,b]=c,[c,a]=a^{p^m},[c,b]=1 \rangle$, where
\mbox{\hskip0.6in}$m>1$ for $p\le 3$;%(B1)

{\rm (H1)}  $\langle a, b, c,d \di a^{4}=b^{2}=c^2=d^2=1,
[a,b]=c,[c,a]=d,[c,b]=[d,a]=[d,b]=\mbox{\hskip0.65in}1 \rangle$;%(A1)

{\rm (H2)}  $\langle a, b, c \di a^{8}=b^{2}=c^2=1,
[a,b]=c,[c,a]=a^{4},[c,b]=1 \rangle$; %(A2)

{\rm (H3)}  $\langle a, b, c \di a^{8}=c^{2}=1, b^2=a^4,
[a,b]=c,[c,a]=b^2,[c,b]=1 \rangle$; %(A3)

{\rm (I1)}  $\langle a, b, c,d \di a^{2^n}=b^{2}=c^2=d^2=1,
[a,b]=c,[c,a]=d,[c,b]=[d,a]=[d,b]=\mbox{\hskip0.65in}1 \rangle$, where $n\ge 3$;%(B5) for $p=2$ and $m=1$, $n\ge 3$;

{\rm (I2)}  $\langle a, b, c \di a^{2^{n+1}}=b^{2}=c^2=1,
[a,b]=c,[c,a]=a^{2^n},[c,b]=1 \rangle$, where $n\ge 3$;%(B1) for $p=2$ and $m=1$, $n\ge 3$;

{\rm (I3)}  $\langle a, b, c \di a^{2^n}=b^4=c^{2}=1,[a,b]=c,[c,a]=b^2,[c,b]=1 \rangle$, where $n\ge 3$;%(B4) for $p=2$ and $m=1$, $n\ge 3$;

{\rm (J1)}  $\langle a, b, c \di a^{p^{n+1}}=b^{p^m}=c^p=1,
[a,b]=c,[c,a]=1,[c,b]=a^{\nu p^n} \rangle$, where \mbox{\hskip0.55in}$m>1$
for $p=2$ and $n> m$,
$\nu=1$ or a fixed quadratic non-residue modular
$p$; %(B2) &(A11) & (A12)

{\rm (J2)}  $\langle a, b, c \di a^{p^{n+1}}=b^{p^m}=c^p=1,
[a,b]=c,[c,a]=a^{p^n},[c,b]=1 \rangle$, where
\mbox{\hskip0.55in}$m>1$ for $p=2$ and $n>m$;%(B1) & (A8)

{\rm (J3)}  $\langle a, b, c \di a^{p^{n}}=b^{p^{m+1}}=c^p=1,
[a,b]=c,[c,a]=1,[c,b]=b^{p^m} \rangle$, where $m>1$
\mbox{\hskip0.55in}for $p=2$ and $n>m$; %(B3) &(A10)

{\rm (J4)}  $\langle a, b, c \di a^{p^{n}}=b^{p^{m+1}}=c^p=1,
[a,b]=c,[c,a]=b^{\nu p^m},[c,b]=1 \rangle$, where
\mbox{\hskip0.55in}$m>1$ for $p=2$ and $n>m$, $\nu=1$ or
a fixed quadratic non-residue modular $p$; %(B4) &(A6) &(A7)

{\rm (J5)}  $\langle a, b, c,d \di a^{p^{n}}=b^{p^{m}}=c^p=d^p=1,
[a,b]=c,[c,a]=1,[c,b]=d, [d,a]=\mbox{\hskip0.55in}[d,b]=1 \rangle$,
where $m>1$ for $p=2$ and $n>m$; %(B6) &(A9)

{\rm (J6)}  $\langle a, b, c,d \di a^{p^{n}}=b^{p^{m}}=c^p=d^p=1,
[a,b]=c,[c,a]=d,[c,b]=1,[d,a]=\mbox{\hskip0.55in}[d,b]=1 \rangle$,
where $m>1$
for $p=2$ and $n>m$. %(B5) &(A5)

\end{thm}
\demo
{\bf Case 1:} $n=m$.

If $n=m=1$, then $p>2$. If $p=3$, then $|G|=3^4$ and hence $G$ is the group of type (F).
If $m>1$ or $p>3$, then, by Theorem \ref{isomorphic-4.1} and \ref{isomorphic-4.2},
$G\cong \bar{G}$ if and only if there exists $X={\left(
 \begin{array}{cc}
 x_{11}& x_{12}\\
 x_{21}& x_{22}
 \end{array}
 \right)}$, an invertible matrix over $F_p$, and $\lambda\in F_p^*$ such that
$$\left(
 \begin{array}{c}
 \bar{w}_{11}\\
 \bar{w}_{21}
  \end{array}
 \right)=
{\lambda^{-1}}
X\left(
 \begin{array}{c}
 w_{11}\\
 w_{21}
 \end{array}
 \right)$$ and
 $$\left(
 \begin{array}{c}
 \bar{w}_{12}\\
 \bar{w}_{22}
  \end{array}
 \right)=\lambda^{-1} |X|X\left(
 \begin{array}{cc}
 w_{12}\\
  w_{22}
 \end{array}
 \right).$$
That is,
\begin{equation}\label{eq 4.7}
\left(
 \begin{array}{cc}
 \bar{w}_{11}&\bar{w}_{12}\\
 \bar{w}_{21}&\bar{w}_{22}
  \end{array}
 \right)=
\lambda^{-1} X\left(
 \begin{array}{cc}
 w_{11}&w_{12}\\
 w_{21}&w_{22}
 \end{array}
 \right)
\left(
 \begin{array}{cc}
 1&0\\
 0&|X|
 \end{array}
 \right).
  \end{equation}

{\bf Subcase 1.1.} $w(G)$ is invertible.

 Assume that $|w(G)|=\nu d^2$ where $\nu=1$ or a fixed quadratic non-residue modular $p$. Let $X=\diag(1,\nu d)(w(G))^{-1}$ and $\lambda=1$. Then $|X|=d^{-1}$ and $w(\bar{G})=\diag(1,\nu)$. Hence we get the group of Type (G1). By
Equation \ref{eq 4.7}, $|w(\bar{G})|=|w(G)|\lambda^{-2}|X|^2$. Hence different $\nu$ give non-isomorphic groups.

\medskip

{\bf Subcase 1.2.} $w(G)$ is not invertible.

Since $G_3\neq 1$, $[c,a]\neq 1$ or $[c,b]\neq 1$.
Without loss of generality, we may assume that $[c,a]\neq 1$. Furthermore, we may assume that $(w_{12},w_{22})=(1,0)$.
Since $w(G)$ is not invertible, $w_{21}=0$. If $w_{11}=0$, then we get the group of Type (G2). If $w_{11}\neq 0$, then, letting $X=\diag(w_{11}^{-1},w_{11}^2)$ and $\lambda=1$, we have $\bar{w}_{11}=1$. Hence we get the group of Type (G3).
It is easy to check that a group of Type (G2) is not isomorphic to that of Type (G3).

\medskip

{\bf Case 2.} $n>m=1$ and $p=2$.

Since $m=1$, $(w_{12},w_{22})=(1,0)$. If $n=2$, then $|G|=2^5$. By checking the list of groups of order
$2^5$, we get the groups of Type (H1)--(H3). If $n\ge 3$, then, by Theorem \ref{isomorphic-4.3},
$G\cong \bar{G}$ if and only if there exists $X={\left(
 \begin{array}{cc}
 1& 0\\
 x_{21}& 1
 \end{array}
 \right)}$ such that
$\left(
 \begin{array}{c}
 \bar{w}_{11}\\
 \bar{w}_{21}
  \end{array}
 \right)=X\left(
 \begin{array}{c}
 w_{11}\\
 w_{21}
 \end{array}
 \right).$ It is easy to get the groups of Type (I1)--(I3).

\medskip

{\bf Case 3.} $n>m$ where $m>1$ for $p=2$.

Suppose that $G$ and $\bar{G}$ are two groups described in the theorem. By Theorem \ref{isomorphic-4.1} and \ref{isomorphic-4.2},
$G\cong \bar{G}$ if and only if there exists $X={\left(
 \begin{array}{cc}
 x_{11}& x_{12}\\
 x_{21}p^{n-m}& x_{22}
 \end{array}
 \right)}$, an invertible matrix over $F_p$, and $\lambda\in F_p^*$ such that
 \begin{equation}\label{eq: 4.10}
\left(
 \begin{array}{c}
 \bar{w}_{11}\\
 \bar{w}_{21}
  \end{array}
 \right)=
 \lambda^{-1}
{\left(
 \begin{array}{cc}
 x_{11}& 0\\
 x_{21}& x_{22}
 \end{array}
 \right)}\left(
 \begin{array}{c}
 w_{11}\\
 w_{21}
 \end{array}
 \right)
 \end{equation} and
 \begin{equation}\label{eq: 4.11}
\left(
 \begin{array}{c}
 \bar{w}_{12}\\
 \bar{w}_{22}
  \end{array}
 \right)=\lambda^{-1}|X|{\left(
 \begin{array}{cc}
 x_{11}& x_{12}\\
 0 & x_{22}
 \end{array}
 \right)}\left(
 \begin{array}{cc}
 w_{12}\\
  w_{22}
 \end{array}
 \right).
 \end{equation}
 By suitably choosing $x_{21}$, that is, using an elementary row operation, we can simplify $\left(
 \begin{array}{cc}
 w_{11}\\
  w_{21}
 \end{array}
 \right)$ to be one of the following three types:
\begin{center}
(a1)  $\left(
 \begin{array}{c}
 w_{11}\\
 0
 \end{array}
 \right)$ where $w_{11}\neq 0$,
(b1)  $\left(
 \begin{array}{c}
 0\\
 w_{21}
 \end{array}
 \right)$ where $w_{21}\neq 0$,
(c1)  $\left(
 \begin{array}{c}
 0\\
 0
 \end{array}
 \right).$
\end{center}
In the following, we assume that both $\left(
 \begin{array}{cc}
 w_{11}\\
  w_{21}
 \end{array}
 \right)$ and
$\left(
 \begin{array}{cc}
 \bar{w}_{11}\\
  \bar{w}_{21}
 \end{array}
 \right)$ are such matrices. By Equation \ref{eq: 4.10} and \ref{eq: 4.11}, it is easy to check that (i)
 different types give non-isomorphic groups, (ii) $G\cong \bar{G}$ if and only if there exists $X={\left(
 \begin{array}{cc}
 x_{11}& x_{12}\\
 0 & x_{22}
 \end{array}
 \right)}$, an invertible matrix over $F_p$, and $\lambda\in F_p^*$ such that
 \begin{equation}\label{eq: 4.10'}
\left(
 \begin{array}{c}
 \bar{w}_{11}\\
 \bar{w}_{21}
  \end{array}
 \right)=
 \lambda^{-1}
{\left(
 \begin{array}{cc}
 x_{11}& 0\\
 0& x_{22}
 \end{array}
 \right)}\left(
 \begin{array}{c}
 w_{11}\\
 w_{21}
 \end{array}
 \right)
 \end{equation} and
 \begin{equation}\label{eq: 4.11'}
\left(
 \begin{array}{c}
 \bar{w}_{12}\\
 \bar{w}_{22}
  \end{array}
 \right)=\lambda^{-1}|X|X\left(
 \begin{array}{cc}
 w_{12}\\
  w_{22}
 \end{array}
 \right).
 \end{equation}
 Similarly, we can simplify $\left(
 \begin{array}{cc}
 w_{12}\\
  w_{22}
 \end{array}
 \right)$ to be one of the following two types (note that $(w_{12},w_{22})\neq (0,0)$ by $G_3\neq 1$):
\begin{center}
(a2)  $\left(
 \begin{array}{c}
 0\\
 w_{22}
 \end{array}
 \right)$ where $w_{22}\neq 0$,
(b2)  $\left(
 \begin{array}{c}
 w_{12}\\
 0
 \end{array}
 \right)$ where $w_{21}\neq 0$.
\end{center}
In the following, we assume that both $\left(
 \begin{array}{cc}
 w_{12}\\
  w_{22}
 \end{array}
 \right)$ and $\left(
 \begin{array}{cc}
 \bar{w}_{12}\\
  \bar{w}_{22}
 \end{array}
 \right)$ are such matrices. By Equation \ref{eq: 4.10'} and \ref{eq: 4.11'}, it is easy to check that (i)
different types give non-isomorphic groups, (ii) $G\cong \bar{G}$ if and only if there exists $X=\diag(x_{11},x_{22})$, an invertible matrix over $F_p$, and $\lambda\in F_p^*$ such that
 \begin{equation}\label{eq: 4.10''}
\left(
 \begin{array}{c}
 \bar{w}_{11}\\
 \bar{w}_{21}
  \end{array}
 \right)=
 \lambda^{-1}
X\left(
 \begin{array}{c}
 w_{11}\\
 w_{21}
 \end{array}
 \right)
\ \rm{and}\
 \left(
 \begin{array}{c}
 \bar{w}_{12}\\
 \bar{w}_{22}
  \end{array}
 \right)=\lambda^{-1}|X|X\left(
 \begin{array}{cc}
 w_{12}\\
  w_{22}
 \end{array}
 \right).
 \end{equation}
 By Equation \ref{eq: 4.10''}, $G\cong \bar{G}$ if and only if there exists $X=\diag(x_{11},x_{22})$, an invertible matrix over $F_p$, and $\lambda\in F_p^*$ such that $w(\bar{G})=\lambda^{-1}Xw(G)\diag(1,|X|)$.
\begin{table}[h]
  \centering
{\tiny
%\scriptsize
%\footnotesize %\small
\begin{tabular}{cccccccc}
\hline
$\left(
 \begin{array}{cc}
 w_{11}\\
  w_{21}
 \end{array}
 \right)$ &$\left(
 \begin{array}{cc}
 w_{12}\\
  w_{22}
 \end{array}
 \right)$ &$w(G)$&$\lambda$&$X$ &Remark & $w(\bar{G})$&Group\\
   \hline
(a1)&  (a2) & $\left(
 \begin{array}{cc}
 w_{11}& 0\\
 0& w_{22}
 \end{array}
 \right)$ &1& 
 $\diag(w_{11}^{-1},d^{-1})$
  &$w_{22}w_{11}^{-1}=\nu d^2$&$\left(
 \begin{array}{cc}
 1& 0\\
 0& \nu
 \end{array}
 \right)$& (J1)\\
 \hline
(a1)&  (b2) &$\left(
 \begin{array}{cc}
 w_{11}& w_{12}\\
 0&0
 \end{array}
 \right)$ &1&$\diag(w_{11}^{-1},w_{11}^2w_{12}^{-1})$&&$\left(
 \begin{array}{cc}
 1& 1\\
 0& 0
 \end{array}
 \right)$& (J2)\\
\hline
(b1)&  (a2) &$\left(
 \begin{array}{cc}
 0& 0\\
 w_{21}& w_{22}
 \end{array}
 \right)$& 1&$\diag(w_{21}^2w_{22}^{-1},w_{21}^{-1})$&&$\left(
 \begin{array}{cc}
 0& 0\\
 1& 1
 \end{array}
 \right)$& (J3)\\
\hline
(b1)&  (b2) &$\left(
 \begin{array}{cc}
 0& w_{12}\\
 w_{21}&0
 \end{array}
 \right)$ &1&$\diag(d^{-1},w_{21}^{-1})$&$w_{12}w_{21}^{-1}=\nu d^2$&$\left(
 \begin{array}{cc}
 0& \nu\\
 1& 0
 \end{array}
 \right)$& (J4)\\
\hline
(c1)&  (a2) &$\left(
 \begin{array}{cc}
 0& 0\\
 0& w_{22}
 \end{array}
 \right)$& 1&$\diag(w_{22}^{-1},1)$&&$\left(
 \begin{array}{cc}
 0& 0\\
 0& 1
 \end{array}
 \right)$& (J5)\\
\hline
(c1)&  (b2) &$\left(
 \begin{array}{cc}
 0& w_{12}\\
 0&0
 \end{array}
 \right)$&1 &$\diag(1, w_{12}^{-1})$&&$\left(
 \begin{array}{cc}
 0& 1\\
 0& 0
 \end{array}
 \right)$& (J6)\\
   \hline
\end{tabular}
\caption{Case 3 in Theorem \ref{c_p^2}}
\label{table 3}
}\end{table}
By Table \ref{table 3}, we get the groups of Type (J1)--(J6).\qed

\medskip

By Theorem \ref{property-4.1/2}, we have the following theorem.

\begin{thm}
\label{property-4.2} Let $G$ be a finite $p$-group listed in Theorem $\ref{c_p^2}$. Then the following conclusions hold:

{\rm (1)} $I_{\max}=1$ for the type {\rm (F)};

{\rm (2)} $I_{\max}=2$ for types {\rm (H1)--(H3), (I1)--(I3)};

{\rm (3)} $I_{\max}=m$  for types {\rm (G1)--(G3), (J2), (J4), (J6)};

{\rm (4)} $I_{\max}=n$  for types {\rm (J1), (J3), (J5)}.
\end{thm}

\section{The case  $G_3\cong C_p$ and $\Phi(G') G_3\cong C_p^2$}

Suppose that $G$ is a finite $p$-group with $G_3\cong C_p$, $\Phi(G') G_3\cong C_p^2$ and $G/\Phi(G')G_3\cong M_p(n,m,1)$, where $n>1$ for $p=2$ and
$n\ge m$. Let $$G/\Phi(G')G_3=\lg \bar{a},\bar{b},\bar{c}\mid \bar{a}^{p^n}=\bar{b}^{p^m}=\bar{c}^{p}=1,[\bar{c},\bar{a}]=[\bar{b},\bar{c}]=1\rg.$$
Then, without loss of generality, we may assume that $G=\lg a,b,c\rg$ where $[a,b]=c$. Since $\Phi(G')\cong C_p$, $c$ is of order $p^2$. Since $G_3\cong C_p$, $C_G(G')$ is maximal in $G$. If $[c,b]\neq 1$, then, without loss of generality, we may assume that $[c,a]=1$. In this case, the type of $C_G(G')/\Phi(G')G_3$ is $(p^n,p^{m-1},p)$. If $[c,b]=1$, then the type of $C_G(G')/\Phi(G')G_3$ is $(p^n,p^{m-1},p)$.

\subsection{The type of $C_G(G')/\Phi(G')G_3$ is $(p^n,p^{m-1},p)$}

Without loss of generality, we may assume that $[a,c]=1$ and $[b,c]\neq 1$. Let $[b,c]=x$. Then $G_3=\lg x\rg$. Since $a^{p^{n}}\in \Phi(G')G_3$, we may assume that $a^{p^{n}}=c^{w_{11}p}x^{w_{12}}$. By a similar reason, we may assume that $b^{p^m}=c^{w_{21}p}x^{w_{22}}$. Let $w(G)=(w_{ij})$. Then we get a $2\times 2$ matrix over $F_p$. $w(G)$ is called a characteristic matrix of $G$. Notice that $w(G)$ will be changed if we change the generators $a,b$. We also call $a,b$ a set of characteristic generators of $w(G)$.

\begin{thm}
\label{property-6.1.0}Suppose that $p$ is an odd prime, $G$ is a finite $p$-group such that
$G_3\cong C_p$, $\Phi(G')G_3=C_p^2$, $G/\Phi(G')G_3\cong M_p(n,m,1)$ and the type of $C_G(G')/\Phi(G')G_3$ is $(p^{n},p^{m-1},p)$ where $n\ge m\ge 2$. Let $w(G)=(w_{ij})$ be a characteristic matrix of $G$. Then the following conclusions hold:

{\rm (1)} $I_{\min}\ge 2$. That is, $G$ does not have an $\mathcal{A}_1$-subgroup of index $p$;

{\rm (2)} $I_{\max}\ge n$;

{\rm (3)} $I_{\max}=2$ if and only if $n=m=2$ and one of the following holds: {\rm (i)} $w_{11}=w_{12}=0$ and $w_{22}\neq 0$; {\rm (ii)} $w_{11}^2-4w_{12}$ is a quadratic non-residue modular $p$.
\end{thm}
\demo We get $I_{\max}$ and $I_{\min}$ by investigating all maximal subgroups of $G$. Let $N=\lg a,b^p,c,x\rg$ and $M_i=\lg a^ib, a^p,c,x\rg$. Then all maximal subgroups of $G$ are $N$ and $M_i$ where $0\le i\le p-1$.

(1) If $w_{12}\neq 0$ or $w_{22}\neq 0$, then $N=\lg a, b^p\rg\ast \lg c\rg\in\mathcal{A}_2$ by Lemma \ref{cor 2.4}.
If $w_{12}=w_{22}=0$, then $N=\lg a,b^p\rg\ast \lg c\rg\times \lg x\rg\in\mathcal{A}_3$ by Lemma \ref{cor 2.4}.
  By calculation, $M_i'=\lg c^p,x\rg=\Phi(G')G_3$ and hence $M_i\not\in\mathcal{A}_1$. Thus $G$ does not have an $\mathcal{A}_1$-subgroup of index $p$.

(2) Let $D=\lg c,b\rg$. Then $D\in\mathcal{A}_1$ and $|G:D|=p^n$. Hence $I_{\max}\ge n$.

(3) If $I_{\max}=2$, then $w_{12}\neq 0$ or $w_{22}\neq 0$ by the proof of (1), and $n=m=2$ by (2). 
 
 Let $A=\lg a^p,b^p,c,x\rg$, $B_r=\lg a^ibc^r,a^p,c^p,x\rg$ and $C_{st}=\lg a^iba^{sp},ca^{tp},c^p,x\rg$ where $0\le r,s,t\le p-1$. Then $A$, $B_r$ and $C_{st}$ are all maximal subgroups of $M_i$.
Since $G\in\mathcal{A}_3$, $M_i\in \mathcal{A}_2$. Hence $A$, $B_r$ and $C_{st}$ are abelian or minimal non-abelian.

First of all, $A$ is the unique abelian maximal subgroup of $M_i$. Since $w_{12}\neq 0$ or $w_{22}\neq 0$, $B_r=\lg a^ibc^r,a^p\rg\in \mathcal{A}_1$.
Since $[ca^{tp},a^iba^{sp}]=c^{tp}x^{-1}$ and $(a^iba^{sp})^{p^2}=c^{(iw_{11}+w_{21})p}x^{iw_{12}+w_{22}}$ and $(ca^{tp})^p=c^{(tw_{11}+1)p}x^{tw_{12}}$, $C_{st}\in \mathcal{A}_1$ if and only if the following equation set about $i$ and $t$ has no solution.
\begin{equation}
  \left\{
  \begin{aligned}
     iw_{11}+w_{21} &=-t(iw_{12}+w_{22}) &&(5.1.1)\\
     tw_{11}+1 &=-t^2w_{12} &&(5.1.2)\\
     \end{aligned}
     \right.
     \end{equation}
By (5.1), $i=tw_{21}+t^2w_{22}$. Hence $C_{st}\in\mathcal{A}_1$ if and only if Equation (5.1.2) has no solution. If $w_{12}=0$, then $w_{11}=0$. Hence (i) holds. If $w_{12}\neq 0$, then $w_{11}^2-4w_{12}$ is a quadratic non-residue modular $p$. Hence (ii) holds.

Conversely, if (i) or (ii) holds and $n=m=2$, then, by using the above argument, it is easy to see $I_{\max}=2$.  \qed

\begin{thm}
\label{property-6.1.00}Let $G$ be a finite $2$-group such that
$G_3\cong C_2$, $\Phi(G')G_3=C_2^2$, $G/\Phi(G')G_3\cong M_2(n,m,1)$ and the type of $C_G(G')/\Phi(G')G_3$ is $(2^{n},2^{m-1},2)$ where $n\ge m\ge 2$. Let $w(G)=(w_{ij})$ be a characteristic matrix of $G$. Then the following conclusions hold:

{\rm (1)} $I_{\min}\ge 2$. That is, $G$ does not have an $\mathcal{A}_1$-subgroup of index $2$;

{\rm (2)} $I_{\max}\ge n$;

{\rm (3)} $I_{\max}=2$ if and only if $n=m=2$ and one of the following holds: {\rm (i)} $w_{22}=1$ and $w_{21}=0$; {\rm (ii)} $w_{22}=0$ and $w_{21}=w_{12}=1$.
\end{thm}
\demo We get $I_{\max}$ and $I_{\min}$ by investigating maximal subgroups of $G$. Let $N=\lg a,b^2,c,x\rg$ and $M_i=\lg a^ib, a^2,c,x\rg$. Then all maximal subgroups of $G$ are $N$, $M_0$ and $M_1$.

(1) If $a^{2^n}=c^2$, then $N=\lg a,b^2\rg\times \lg ca^{2^{n-1}}\rg\in\mathcal{A}_2$ by Lemma \ref{cor 2.4}.
If $a^{2^n}=x$, then $N=\lg a,b^2\rg\ast \lg ca^{2^{n-1}}\rg\in\mathcal{A}_2$.
If $b^{2^m}=c^2$, then $N=\lg a,b^2\rg\times \lg cb^{2^{m-1}}\rg\in\mathcal{A}_2$.
If $b^{2^m}=x$, then $N=\lg a,b^2\rg\ast \lg cb^{2^{m-1}}\rg\in\mathcal{A}_2$.
If $a^{2^n}=c^{2i}x^i$ and $b^{2^m}=c^{2j}x^j$, then $N=\lg a^2,b\rg\times \lg c\rg\in\mathcal{A}_3$ by Lemma \ref{cor 2.4}.
By calculation, $M_i'=\lg c^2,y\rg=\Phi(G')G_3$ and hence $M_i\not\in\mathcal{A}_1$. To sum up, $G$ does not have an $\mathcal{A}_1$-subgroup of index $2$.

(2) Let $D=\lg c,b\rg$. Then $D\in\mathcal{A}_1$ and $|G:D|=2^n$. Hence $I_{\max}\ge n$.

(3) If $I_{\max}=2$, then $n=m=2$ by (2).
 
 Let $A=\lg a^2,b^2,c,x\rg$, $B_r=\lg a^ibc^r,a^2,c^2,x\rg$ and $C_{st}=\lg a^iba^{2s},ca^{2t},c^2,x\rg$ where $r,s,t=0,1$. Then $A$, $B_r$ and $C_{st}$ are all maximal subgroups of $M_i$.
Since $G\in\mathcal{A}_3$, $M_i\in \mathcal{A}_2$. Hence $A$, $B_r$ and $C_{st}$ are abelian or minimal non-abelian.

First of all, $A$ is the unique abelian maximal subgroup of $M_i$. Since $[a^ibc^r,a^2]=c^2$, $B_r\in \mathcal{A}_1$ if and only if $w_{12}\neq 0$ or $w_{22}\neq 0$.
Since $[a^iba^{2s},ca^{2t}]=c^{2t}x$ and $(a^iba^{2s})^{4}=c^{2(iw_{11}+w_{21}+i)}x^{iw_{12}+w_{22}}$ and $(ca^{2t})^2=c^{2(tw_{11}+1)}x^{tw_{12}}$, $C_{st}\in \mathcal{A}_1$ if and only if the following equation set about $i$ and $t$ has no solution.
\begin{equation}
  \left\{
  \begin{aligned}
     iw_{11}+w_{21}+i &=t(iw_{12}+w_{22}) &&(5.2.1)\\
     tw_{11}+1 &=t^2w_{12}=tw_{12} &&(5.2.2)\\
     \end{aligned}
     \right.
     \end{equation}
If Equation set (5.2) has a solution, then, by (5.2.2), $t=w_{11}+w_{12}=1$, and by (5.2.1), $w_{21}+w_{22}=0$. Hence $C_{st}\in \mathcal{A}_1$ if and only if $w_{21}+w_{22}=1$ or $w_{11}+w_{12}=0$. Since $w_{11}+w_{12}=1$ or $w_{21}+ w_{22}=1$ by the proof of (1), $w_{21}+w_{22}=1$. If $w_{22}=1$, then $w_{21}=0$ and (i) holds.
If $w_{22}=0$, then $w_{21}=w_{12}=1$ and (ii) holds.

Conversely, if (1) or (2) holds and $n=m=2$, then, by using the above argument, it is easy see $I_{\max}=2$.\qed

\begin{thm}\label{isomorphic-6.1.1}
Suppose that $G$ and $\bar{G}$ are finite $p$-groups such that $G_3\cong C_p$, $\Phi(G') G_3\cong C_p^2$, $G/\Phi(G')G_3\cong M_p(n,m,1)$, and the type of $C_G(G')/\Phi(G')G_3$ is $(p^n,p^{m-1},p)$ where $n\ge 3$ for $p=2$ and $n\ge m\ge 2$. Let two characteristic matrices of $G$ and $\bar{G}$ be $w(G)=(w_{ij})$ and $w(\bar{G})=(\bar{w}_{ij})$ respectively. Then $G\cong \bar{G}$ if and only if there exists an invertible matrix $X={\left(
 \begin{array}{cc}
 x_{11}& 0\\
 x_{21}& x_{22}
 \end{array}
 \right)}$ over $F_p$ such that $w(\bar{G})=Xw(G)\diag(x_{11}^{-1}x_{22}^{-1},x_{11}^{-1}x_{22}^{-2})$.
 \end{thm}
\demo Suppose that $a,b$ and
$\bar{a},\bar{b}$ are two set of characteristic generators of $w(G)$ and $w(\bar{G})$ respectively.
Let $\theta$ be an isomorphism from $\bar{G}$ onto $G$. We have $\Phi(\bar{G})^\theta=\Phi(G)$, $\Omega_m(\bar{G})^\theta=\Omega_m(G)$ and $(C_{\bar{G}}(\bar{G}'))^\theta=C_G(G')$ since these six subgroups are characteristic in $G$ or $\bar{G}$. Note that $C_{\bar{G}}(\bar{G}')=\lg \bar{a},\Phi(\bar{G})\rg$ and $C_G(G')=\lg a,\Phi(G)\rg$ respectively. So we may let
$$\bar{a}^\theta=a^{x_{11}}\phi_1,\ \bar{b}^\theta=a^{x_{21}p^{n-m}}b^{x_{22}}\phi_2$$
where $\phi_1\in\Phi(G)$, $\phi_2\in\Phi(G)\cap\Omega_m(G)$ and $x_{11}x_{22}\in F_p^*$. By calculation, $$\bar{c}^\theta=[\bar{a},\bar{b}]^\theta=[\bar{a}^\theta,\bar{b}^\theta]\equiv [a^{x_{11}},a^{x_{21}p^{n-m}}b^{x_{22}}]\equiv
c^{x_{11}x_{22}}\ (\mod G_3)$$
and $$\bar{x}^\theta=[\bar{b},\bar{c}]^\theta=[\bar{b}^\theta,\bar{c}^\theta]=[a^{x_{21}p^{n-m}}b^{x_{22}},c^{x_{11}x_{22}}]=x^{x_{11}x_{22}^2}.$$
 Transforming $\bar{c}^{\bar{w}_{11}p}\bar{x}^{\bar{w}_{12}}=\bar{a}^{p^{n}}$ by $\theta$, we have  ${c}^{\bar{w}_{11}x_{11}x_{22}p}{x}^{\bar{w}_{12}x_{11}x_{22}^2}=a^{x_{11}p^n}$. Hence
\begin{equation}\label{eq 6.1}
(\bar{w}_{11},\bar{w}_{12})\left(
 \begin{array}{cc}
 x_{11}x_{22}& 0\\
 0& x_{11}x_{22}^2
 \end{array}
 \right)=
 (x_{11},0)
 \left(
 \begin{array}{cc}
 w_{11}& w_{12}\\
 w_{21}& w_{22}
 \end{array}
 \right).
\end{equation}
 Transforming $\bar{c}^{\bar{w}_{21}p}\bar{x}^{\bar{w}_{22}}=\bar{b}^{p^{m}}$ by $\theta$, we have ${c}^{\bar{w}_{21}x_{11}x_{22}p}{x}^{\bar{w}_{22}x_{11}x_{22}^2}=a^{x_{21}p^n}b^{x_{22}p^m}$. Hence
\begin{equation}\label{eq 6.2}
(\bar{w}_{21},\bar{w}_{22})\left(
 \begin{array}{cc}
 x_{11}x_{22}& 0\\
 0& x_{11}x_{22}^2
 \end{array}
 \right)=
 (x_{21},x_{22})
 \left(
 \begin{array}{cc}
 w_{11}& w_{12}\\
 w_{21}& w_{22}
 \end{array}
 \right).
\end{equation}
By Equation \ref{eq 6.1} and \ref{eq 6.2},
\begin{equation}\label{eq 6.3}
\left(
 \begin{array}{cc}
 \bar{w}_{11}& \bar{w}_{12}\\
 \bar{w}_{21}& \bar{w}_{22}
  \end{array}
 \right)=
 \left(
 \begin{array}{cc}
 x_{11}&0\\
 x_{21}& x_{22}
 \end{array}
 \right)
 \left(
 \begin{array}{cc}
 w_{11}& w_{12}\\
 w_{21}& w_{22}
 \end{array}
 \right)\diag(x_{11}^{-1}x_{22}^{-1},x_{11}^{-1}x_{22}^{-2}).
\end{equation}

Conversely, if there exists an invertible matrix $X=\left(
 \begin{array}{cc}
 x_{11}& 0\\
 x_{21}& x_{22}
 \end{array}
 \right)$ over
$F_p$ such that Equation \ref{eq 6.3}, then, by using the above argument, it is easy to check that the map $\theta:\bar{a}\mapsto a^{x_{11}},\bar{b}\mapsto a^{x_{21}p^{n-m}}b^{x_{22}}$ is an isomorphism from
$\bar{G}$ onto $G$.
 \qed

\begin{thm}\label{th=5.4}
 Let $G$ be a finite $p$-group such that $G_3\cong C_p$, $\Phi(G') G_3\cong C_p^2$, $G/\Phi(G')G_3\cong M_p(n,m,1)$ where $n\ge m$, and the type of $C_G(G')/\Phi(G')G_3$ is $(p^n,p^{m-1},p)$. Then $G$ is one of the following non-isomorphic groups:

{\rm (K1)} $\langle a, b, c \di a^{8}=b^{4}=c^{4}=1,
[a,b]=c, [c,a]=1,[c,b]=c^{2}a^{4}\rangle$;%(??)

{\rm (K2)} $\langle a, b, c \di a^{8}=b^{8}=1,c^{2}=b^4,
[a,b]=c, [c,a]=1,[c,b]=a^{4}b^{4}\rangle$;%(I_{\max}=2)

{\rm (K3)} $\langle a, b, c,d \di a^{8}=b^{4}=d^2=1,c^{2}=a^4,
[a,b]=c, [c,a]=1,[c,b]=d,[d,a]=\mbox{\hskip0.65in}[d,b]=1\rangle$;%(??)

{\rm (K4)} $\langle a, b, c \di a^{8}=b^{8}=1,c^{2}=a^4,
[a,b]=c, [c,a]=1,[c,b]=b^4\rangle$;%(I_{\max}=2)

{\rm (K5)} $\langle a, b, c,d \di a^{8}=d^2=1,b^{4}=c^{2}=a^4,
[a,b]=c, [c,a]=1,[c,b]=d,[d,a]=\mbox{\hskip0.65in}[d,b]=1\rangle$;%(??)

{\rm (K6)} $\langle a, b, c \di a^{8}=b^{8}=1, c^{2}=a^4,
[a,b]=c, [c,a]=1,[c,b]=a^{4}b^{4}\rangle$;%(??)

{\rm (K7)} $\langle a, b, c \di a^{8}=b^{4}=c^{4}=1,
[a,b]=c, [c,a]=1,[c,b]=a^{4}\rangle$;%(??)

{\rm (K8)} $\langle a, b, c \di a^{8}=b^{8}=1,c^{2}=b^4,
[a,b]=c, [c,a]=1,[c,b]=a^{4}\rangle$;%(I_{\max}=2)

{\rm (K9)} $\langle a, b, c,d \di a^{4}=b^{4}=c^{4}=d^2=1,
[a,b]=c,[c,b]=d,[c,a]=[d,a]=[d,b]=\mbox{\hskip0.65in}1\rangle$;%(??)

{\rm (K10)} $\langle a, b, c \di a^{4}=b^{8}=c^{4}=1,
[a,b]=c, [c,a]=1,[c,b]=b^{4}\rangle$;%(I_{\max}=2)

 {\rm (L1)} $\langle a, b, c \di a^{p^{n+1}}=b^{p^{m+1}}=1,
[a,b]=c,c^p=a^{p^n}b^{sp^m}, [c,a]=1,[b,c]=b^{p^m}\rangle$, \mbox{\hskip0.65in}where $n\ge 3$ for $p=2$ and $n\ge m\ge 2$, $s\in F_p$;%(??)

 {\rm (L2)} $\langle a, b, c \di a^{p^{n+1}}=b^{p^{m}}=c^{p^2}=1,
[a,b]=c, [c,a]=1,[c,b]=c^{tp}a^{-tp^n}\rangle$, where \mbox{\hskip0.65in}$n\ge 3$ for $p=2$ and $n\ge m\ge 2$, $t\in F_p^*$;%(??)

{\rm (L3)} $\langle a, b, c,d \di a^{p^{n+1}}=b^{p^{m}}=d^p=1, c^p=a^{p^n},
[a,b]=c, [c,a]=1,[c,b]=\mbox{\hskip0.65in}d,[d,a]=[d,b]=1\rangle$, where $n\ge 3$ for $p=2$ and $n\ge m\ge 2$;%(??)

{\rm (L4)} $\langle a, b, c \di a^{p^{n+1}}=b^{p^{m+1}}=1,
[a,b]=c,c^p=b^{p^m}, [c,a]=1,[b,c]=a^{\nu p^n}\rangle$, where \mbox{\hskip0.65in}$n\ge 3$ for $p=2$ and $n\ge m\ge 2$, $\nu=1$ or a fixed quadratic non-residue \mbox{\hskip0.65in}modular $p$;%(??)

{\rm (L5)} $\langle a, b, c \di a^{p^{n+1}}=b^{p^{m}}=c^{p^2}=1,
[a,b]=c, [c,a]=1,[b,c]=a^{\nu p^n}\rangle$, where $n\ge 3$ \mbox{\hskip0.65in}for $p=2$ and $n\ge m\ge 2$, $\nu=1$ or a fixed quadratic non-residue modular $p$;%(??)

 {\rm (L6)} $\langle a, b, c \di a^{p^{n}}=b^{p^{m+1}}=c^{p^2}=1,
[a,b]=c, [c,a]=1,[c,b]=c^{p}b^{-p^m},[b^{p^m},a]=\mbox{\hskip0.65in}1\rangle$, where $n\ge 3$ for $p=2$ and $n\ge m\ge 2$;%(??)

{\rm (L7)} $\langle a, b, c,d \di a^{p^{n}}=b^{p^{m+1}}=d^p=1, c^p=b^{p^m}
[a,b]=c, [c,a]=1,[c,b]=d,[d,a]=\mbox{\hskip0.65in}[d,b]=1\rangle$, where $n\ge 3$ for $p=2$ and $n\ge m\ge 2$;%(??)

 {\rm (L8)} $\langle a, b, c \di a^{p^{n}}=b^{p^{m+1}}=c^{p^2}=1,
[a,b]=c, [c,a]=1,[b,c]=b^{p^m}\rangle$, where $n\ge 3$ \mbox{\hskip0.65in}for $p=2$ and $n\ge m\ge 2$;%(??)

{\rm (L9)} $\langle a, b, c,d \di a^{p^{n}}=b^{p^{m}}=c^{p^2}=d^p=1,
[a,b]=c, [c,a]=1,[c,b]=d,[d,a]=\mbox{\hskip0.65in}[d,b]=1\rangle$, where $n\ge 3$ for $p=2$ and $n\ge m\ge 2$.%(??)

 \end{thm}
\demo
If $m=1$, then $1=[a,b^{p}]=c^p[c,b]^{p\choose 2}$. It follows that $\Phi(G')\le G_3$, a contradiction. Hence $m\ge 2$. If $p=n=m=2$, then $|G|=2^7$. By checking the list of groups of order
$2^7$, we get the groups of Type (K1)--(K10). In the following, we may assume that $n\ge 3$ for $p=2$.

 Suppose that $G$ and $\bar{G}$ are two groups described in the theorem. By Theorem \ref{isomorphic-6.1.1}, $G\cong \bar{G}$ if and only if there exists $X={\left(
 \begin{array}{cc}
 x_{11}& 0\\
 x_{21}& x_{22}
 \end{array}
 \right)}$, an invertible matrix over $F_p$, such that
 \begin{equation}\label{eq: 5.6}
 w(\bar{G})=Xw(G)\diag(x_{11}^{-1}x_{22}^{-1},x_{11}^{-1}x_{22}^{-2}).
 \end{equation}
By suitably choosing $x_{21}$, that is, using an elementary row operation, we can simplify $w(G)$ to be one of the following types:
(a)  $\left(
 \begin{array}{cc}
 w_{11}& w_{12}\\
 0& w_{22}
 \end{array}
 \right)$ where $w_{11}\neq 0$,
(b)  $\left(
 \begin{array}{cc}
 0  & w_{12}\\
 w_{21}& 0
 \end{array}
 \right)$ where $w_{12}\neq 0$, and
(c)  $\left(
 \begin{array}{cc}
 0&0\\
 w_{21}&w_{22}
 \end{array}
 \right).$
In the following, we assume that both $w(G)$ and
$w(\bar{G})$ are such matrices. By Equation \ref{eq: 5.6}, it is easy to check that (i)
 different types give non-isomorphic groups, (ii) $G\cong \bar{G}$ if and only if there exists $X=\diag(x_{11},x_{22})$, an invertible matrix over $F_p$, such that
$ w(\bar{G})=Xw(G)\diag(x_{11}^{-1}x_{22}^{-1},x_{11}^{-1}x_{22}^{-2})$.
\begin{table}[h]
  \centering
%\caption{Properties of the groups of Type (A) in Theorem
%\ref{main}}\label{table: 8.1}
{
%\scriptsize
\tiny
%\footnotesize
%\small
\begin{tabular}{ccccccc}
$w(G)$& Cases &$X$ &Remark 1& $w(\bar{G})$& Group & Remark 2\\
   \hline
(a)& $w_{22}\neq 0$ & $\diag(w_{11}^{-1}w_{22},w_{11})$ &&${\left(
 \begin{array}{cc}
 1& w_{11}^{-2}w_{12}\\
 0& 1
 \end{array}
 \right)}$& (L1) & $s=-w_{11}^{-2}w_{12}$\\
(a) &$w_{22}=0$ & $\diag(1,w_{11})$ &&${\left(
 \begin{array}{cc}
 1& w_{11}^{-2}w_{12}\\
 0& 0
 \end{array}
 \right)}$& $\begin{array}{c} {\rm (L2)\ if\ }w_{12}\neq 0\\ {\rm (L3)\ if\ }w_{12}=0\end{array}$ &$\begin{array}{c} t=w_{11}^{2}w_{12}^{-1}\\
 \\ \end{array}$ \\
%\hline
(b) &$w_{21}\neq 0$&
 $\diag(w_{21},z)$& $w_{12}=\nu z^2$&${\left(
 \begin{array}{cc}
 0& \nu\\
 1& 0
 \end{array}
 \right)}$&(L4) & \\
(b) & $w_{21}=0$& $\diag(1,z)$ &$w_{12}=\nu z^2$&${\left(
 \begin{array}{cc}
 0& \nu\\
 0& 0
 \end{array}
 \right)}$&(L5) & \\
%\hline
(c)& $w_{21}\neq0,w_{22}\neq 0$& $\diag(w_{21},w_{21}^{-1}w_{22})$ &&${\left(
 \begin{array}{cc}
 0& 0\\
 1& 1
 \end{array}
 \right)}$&(L6) & \\
(c) & $w_{22}=0$, $w_{21}\neq 0$& $\diag(w_{21},1)$ &&${\left(
 \begin{array}{cc}
 0& 0\\
 1& 0
 \end{array}
 \right)}$&(L7) & \\
(c) &$w_{21}=0$, $w_{22}\neq 0$& $\diag(1,w_{22})$ & &${\left(
 \begin{array}{cc}
 0& 0\\
 0& 1
 \end{array}
 \right)}$&(L8) & \\
(c) & $w_{21}=w_{22}=0$ &  &&${\left(
 \begin{array}{cc}
 0& 0\\
 0& 0
 \end{array}
 \right)}$&(L9) & \\
   \hline
\end{tabular}
\caption{Theorem \ref{th=5.4} where $n\ge 3$ for $p=2$}}
\label{table 4}
\end{table}
By Table \ref{table 4}, we get the groups of type (L1)--(L9).\qed

\subsection{The type of $C_G(G')/\Phi(G')G_3$ is $(p^{n-1},p^{m},p)$}

If $n=m$, then the type of $C_G(G')/\Phi(G')G_3$ is also $(p^{n},p^{m-1},p)$. The problem is reduced to the case in last subsection. Hence we may assume that $n>m$ in this section.
In this case, $[a,c]\neq 1$ and $[b,c]=1$. Let $[a,c]=y$. Then $G_3=\lg y\rg$. Since $a^{p^{n}}\in \Phi(G')G_3$, we may assume that $a^{p^{n}}=c^{w_{11}p}y^{w_{12}}$. By a similar reason, we may assume that $b^{p^m}=c^{w_{21}p}y^{w_{22}}$. Let $w(G)=(w_{ij})$. Then we get a $2\times 2$ matrix over $F_p$. $w(G)$ is called a characteristic matrix of $G$. Notice that $w(G)$ will be changed if we change the generators $a,b$. We also call $a,b$ a set of characteristic generators of $w(G)$.

\begin{thm}
\label{property-6.2.0}Suppose that $p$ is an odd prime, $G$ is a finite $p$-group such that
$G_3\cong C_p$, $\Phi(G')G_3=C_p^2$, $G/\Phi(G')G_3\cong M_p(n,m,1)$ and the type of $C_G(G')/\Phi(G')G_3$ is $(p^{n-1},p^{m},p)$ where $n>m$. Let $w(G)=(w_{ij})$ be a characteristic matrix of $G$. Then the following conclusions hold:

{\rm (1)} $I_{\min}\ge 2$. That is, $G$ does not have an $\mathcal{A}_1$-subgroup of index $p$;

{\rm (2)} $I_{\max}\ge m$;

{\rm (3)} $I_{\max}=2$ if and only if $m=2$ and one of the following conditions holds: {\rm (i)} $w_{12}\neq 0$ and $w_{11}|w(G)|\neq w_{12}^2$; {\rm (ii)} $w_{12}=0$, $w_{11}\neq 0$ and $w_{22}\neq 0$; {\rm (iii)} $w_{21}^2+4w_{22}$ is a quadratic non-residue modular $p$.
\end{thm}
\demo We get $I_{\max}$ and $I_{\min}$ by investigating maximal subgroups of $G$. Let $N=\lg b,a^p,c,y\rg$ and $M_i=\lg ab^i, b^p,c,y\rg$. Then all maximal subgroups of $G$ are $N$ and $M_i$ where $0\le i\le p-1$.

(1) If $w_{12}\neq 0$ or $w_{22}\neq 0$, then $N=\lg a^p,b\rg\ast \lg c\rg\in\mathcal{A}_2$ by Lemma \ref{cor 2.4}.
If $w_{12}=w_{22}=0$, then $N=\lg a^p,b\rg\ast \lg c\rg\times \lg y\rg\in\mathcal{A}_3$ by Lemma \ref{cor 2.4}.
  By calculation, $M_i'=\lg c^p,y\rg=\Phi(G')G_3$ and hence $M_i\not\in\mathcal{A}_1$. To sum up, $G$ does not have an $\mathcal{A}_1$-subgroup of index $p$.

(2) Let $D=\lg c,a\rg$. Then $D\in\mathcal{A}_1$ and $|G:D|=p^m$. Hence $I_{\max}\ge m$.

(3) If $I_{\max}=2$, then $w_{12}\neq 0$ or $w_{22}\neq 0$ by the proof of (1), and $m=2$ by (2). Since $G\in\mathcal{A}_3$, $M_i\in \mathcal{A}_2$. Notice that $M_i\cong M_0$ since $n\ge 3$. We only need to investigate all maximal subgroups of $M_0$. Let $A=\lg a^p,b^p,c,y\rg$, $B_r=\lg ac^r,b^p,c^p,y\rg$ and $C_{st}=\lg ab^{sp},cb^{tp},c^p,y\rg$ where $0\le r,s,t\le p-1$. Then $A$, $B_r$ and $C_{st}$ are all maximal subgroups of $M_0$.

First of all, $A$ is the unique abelian maximal subgroup of $M_0$. Since $w_{12}\neq 0$ or $w_{22}\neq 0$, $B_r=\lg ac^r,b^p\rg\in \mathcal{A}_1$.
Since $[ab^{sp},cb^{tp}]=c^{tp}y$ and $(ab^{sp})^{p^n}=c^{w_{11}p}y^{w_{12}}$ and $(cb^{tp})^p=c^{(tw_{21}+1)p}y^{tw_{22}}$, $C_{st}\in \mathcal{A}_1$ if and only if the following equation set about $t$ has no solution.
\begin{equation}
  \left\{
  \begin{aligned}
     w_{11} &=tw_{12} &&(5.6.1)\\
     tw_{21}+1 &=t^2w_{22} &&(5.6.2)\\
     \end{aligned}
     \right.
     \end{equation}
If $w_{12}\neq 0$, then the solution of Equation (5.6.1) is $t=w_{12}^{-1}w_{11}$. In this case, $C_{st}\in\mathcal{A}_1$ if and only if $t=w_{12}^{-1}w_{11}$ is not a solution of Equation (5.6.2). Hence  $w_{11}|w(G)|\neq w_{12}^2$. (i) holds.
If $w_{12}=0$, then $w_{22}\neq 0$. In this case, Equation (5.6.1) has no solution if and only if $w_{11}\neq 0$. Hence (ii) holds. Equation (5.6.2) has no solution if and only if $w_{21}^2+4w_{22}$ is a quadratic non-residue modular $p$. Hence (iii) holds.

Conversely, if (i), (ii) or (iii) holds and $m=2$, then, by using the above argument, it is easy to see $I_{\max}=2$.  \qed

\begin{thm}
\label{property-6.2.00}Suppose that $G$ is a finite $2$-group such that
$G_3\cong C_2$, $\Phi(G')G_3=C_2^2$, $G/\Phi(G')G_3\cong M_2(n,m,1)$ and the type of $C_G(G')/\Phi(G')G_3$ is $(2^{n-1},2^{m},2)$ where $n>m$. Let $w(G)=(w_{ij})$ be a characteristic matrix of $G$. Then the following conclusions hold:

{\rm (1)} $I_{\min}\ge 2$. That is, $G$ does not have an $\mathcal{A}_1$-subgroup of index $2$;

{\rm (2)} $I_{\max}\ge m$;

{\rm (3)} $I_{\max}=2$ if and only if $m=2$ and one of the following holds: {\rm (i)} $w_{11}=0$ and $w_{12}=1$; {\rm (ii)} $w_{12}=0$ and $w_{11}=w_{22}=1$.
\end{thm}
\demo We get $I_{\max}$ and $I_{\min}$ by investigating maximal subgroups of $G$. Let $N=\lg b,a^2,c,y\rg$ and $M_i=\lg ab^i, b^2,c,y\rg$. Then all maximal subgroups of $G$ are $N$, $M_0$ and $M_1$.

(1) If $a^{2^n}=c^2$, then $N=\lg a^2,b\rg\times \lg ca^{2^{n-1}}\rg\in\mathcal{A}_2$ by Lemma \ref{cor 2.4}.
If $a^{2^n}=y$, then $N=\lg a^2,b\rg\ast \lg ca^{2^{n-1}}\rg\in\mathcal{A}_2$.
If $b^{2^m}=c^2$, then $N=\lg a^2,b\rg\times \lg cb^{2^{m-1}}\rg\in\mathcal{A}_2$.
If $b^{2^m}=y$, then $N=\lg a^2,b\rg\ast \lg cb^{2^{m-1}}\rg\in\mathcal{A}_2$.
If $a^{2^n}=c^{2i}y^i$ and $b^{2^m}=c^{2j}y^j$, then $N=\lg a^2,b\rg\times \lg c\rg\in\mathcal{A}_3$ by Lemma \ref{cor 2.4}.
By calculation, $M_i'=\lg c^2,y\rg=\Phi(G')G_3$ and hence $M_i\not\in\mathcal{A}_1$. To sum up, $G$ does not have an $\mathcal{A}_1$-subgroup of index $2$.

(2) Let $D=\lg c,a\rg$. Then $D\in\mathcal{A}_1$ and $|G:D|=2^m$. Hence $I_{\max}\ge m$.

(3) If $I_{\max}=2$, then $m=2$ by (2). Since $G\in\mathcal{A}_3$, $M_i\in \mathcal{A}_2$. Notice that $M_i\cong M_0$ since $n\ge 3$. We only need to investigate all maximal subgroups of $M_0$. Let $A=\lg a^2,b^2,c,y\rg$, $B_r=\lg ac^r,b^2,c^2,y\rg$ and $C_{st}=\lg ab^{2s},cb^{2t},c^2,y\rg$ where $r,s,t=0,1$. Then $A$, $B_r$ and $C_{st}$ are all maximal subgroups of $M_0$.

First of all, $A$ is the unique abelian maximal subgroup of $M_0$. Since $[ac^r,b^2]=c^2$, $B_r\in \mathcal{A}_1$ if and only if $w_{12}\neq 0$ or $w_{22}\neq 0$.
Since $[ab^{2s},cb^{2t}]=c^{2t}y$ and $(ab^{2s})^{2^n}=c^{2w_{11}}y^{w_{12}}$ and $(cb^{2t})^2=c^{tw_{21}+1}y^{tw_{22}}$, $C_{st}\in \mathcal{A}_1$ if and only if the following equation set about $t$ has no solution.
\begin{equation}
  \left\{
  \begin{aligned}
     w_{11} &=tw_{12} &&(5.7.1)\\
     tw_{21}+1 &=t^2w_{22} &&(5.7.2)\\
     \end{aligned}
     \right.
     \end{equation}
If Equation set (5.7) has a solution, then, by (5.7.2), $t=w_{21}+w_{22}=1$, and by (5.7.1), $w_{11}=w_{12}$. Hence $C_{st}\in \mathcal{A}_1$ if and only if $w_{11}+w_{12}=1$ or $w_{21}+w_{22}=0$. Since $w_{11}+w_{12}=1$ or $w_{21}+ w_{22}=1$ by the proof of (1), $w_{11}+w_{12}=1$.
If $w_{12}=1$, then $w_{11}=0$ and (i) holds.
If $w_{12}=0$, then $w_{11}=w_{22}=1$ and (ii) holds.

Conversely, if (1) or (2) holds and $m=2$, then, by using the above argument, it is easy to see $I_{\max}=2$.\qed

\begin{thm}\label{isomorphic-6.2.1}
Suppose that $G$ and $\bar{G}$ are finite $p$-groups such that $G_3\cong C_p$, $\Phi(G') G_3\cong C_p^2$, $G/\Phi(G')G_3\cong M_p(n,m,1)$ and the type of $C_G(G')/\Phi(G')G_3$ is $(p^{n-1},p^{m},p)$ where $n>m\ge 2$.
Let two characteristic matrices of $G$ and $\bar{G}$ be $w(G)=(w_{ij})$ and $w(\bar{G})=(\bar{w}_{ij})$ respectively.
Then $G\cong \bar{G}$ if and only if there exists $X={\left(
 \begin{array}{cc}
 x_{11}& 0\\
 x_{21}& x_{22}
 \end{array}
 \right)}$, an invertible matrix over $F_p$, such that $
w(\bar{G})=Xw(G)\diag( x_{11}^{-1}x_{22}^{-1},x_{11}^{-2}x_{22}^{-1})$.
 \end{thm}
\demo Suppose that $a,b$ and
$\bar{a},\bar{b}$ are two set of characteristic generators of $w(G)$ and $w(\bar{G})$ respectively.
Let $\theta$ be an isomorphism from $\bar{G}$ onto $G$. We have $\Phi(\bar{G})^\theta=\Phi(G)$, $\Omega_m(\bar{G})^\theta=\Omega_m(G)$ and $(C_{\bar{G}}(\bar{G}'))^\theta=C_G(G')$ since these six subgroups are characteristic in $G$ or $\bar{G}$. Note that $C_{\bar{G}}(\bar{G}')=\lg \bar{a},\Phi(\bar{G})\rg$ and $C_G(G')=\lg a,\Phi(G)\rg$ respectively. So we may let
$$\bar{a}^\theta=a^{x_{11}}b^{x_{12}}\phi_1,\ \bar{b}^\theta=a^{x_{21}p^{n-m}}b^{x_{22}}\phi_2$$
where $\phi_1\in\Phi(G)$, $\phi_2\in\Phi(G)\cap\Omega_m(G)$ and $x_{11}x_{22}\in F_p^*$. By calculation, 
we have
$$\bar{c}^\theta=[\bar{a},\bar{b}]^\theta=[\bar{a}^\theta,\bar{b}^\theta]\equiv [a^{x_{11}},a^{x_{21}p^{n-m}}b^{x_{22}}]\equiv
c^{x_{11}x_{22}}\ (\mod G_3)$$
and $$\bar{y}^\theta=[\bar{a},\bar{c}]^\theta=[\bar{a}^\theta,\bar{c}^\theta]=[a^{x_{11}}b^{x_{12}},c^{x_{11}x_{22}}]=y^{x_{11}^2x_{22}}.$$
 Transforming $\bar{c}^{\bar{w}_{11}p}\bar{y}^{\bar{w}_{12}}=\bar{a}^{p^{n}}$ by $\theta$, we have ${c}^{\bar{w}_{11}x_{11}x_{22}p}{y}^{\bar{w}_{12}x_{11}^2x_{22}}=a^{x_{11}p^n}$. It follows that
\begin{equation}\label{eq 6.4}
(\bar{w}_{11},\bar{w}_{12})\left(
 \begin{array}{cc}
 x_{11}x_{22}& 0\\
 0& x_{11}^2x_{22}
 \end{array}
 \right)=
 (x_{11},0)
 \left(
 \begin{array}{cc}
 w_{11}& w_{12}\\
 w_{21}& w_{22}
 \end{array}
 \right).
\end{equation}
Transforming $\bar{c}^{\bar{w}_{21}p}\bar{y}^{\bar{w}_{22}}=\bar{b}^{p^{m}}$ by $\theta$, we have ${c}^{\bar{w}_{21}x_{11}x_{22}p}{y}^{\bar{w}_{22}x_{11}^2x_{22}}=a^{x_{21}p^n}b^{x_{22}p^m}$. It follows that
\begin{equation}\label{eq 6.5}
(\bar{w}_{21},\bar{w}_{22})\left(
 \begin{array}{cc}
 x_{11}x_{22}& 0\\
 0& x_{11}^2x_{22}
 \end{array}
 \right)=
 (x_{21},x_{22})
 \left(
 \begin{array}{cc}
 w_{11}& w_{12}\\
 w_{21}& w_{22}
 \end{array}
 \right).
\end{equation}
By Equation \ref{eq 6.4} and \ref{eq 6.5},
\begin{equation}\label{eq 6.6}
\left(
 \begin{array}{cc}
 \bar{w}_{11}& \bar{w}_{12}\\
 \bar{w}_{21}& \bar{w}_{22}
  \end{array}
 \right)=
 \left(
 \begin{array}{cc}
 x_{11}&0\\
 x_{21}& x_{22}
 \end{array}
 \right)
 \left(
 \begin{array}{cc}
 w_{11}& w_{12}\\
 w_{21}& w_{22}
 \end{array}
 \right)\diag( x_{11}^{-1}x_{22}^{-1},x_{11}^{-2}x_{22}^{-1}).
\end{equation}

Conversely, if there exists an invertible matrix $X=\left(
 \begin{array}{cc}
 x_{11}& 0\\
 x_{21}& x_{22}
 \end{array}
 \right)$ over
$F_p$ such that Equation \ref{eq 6.6}, then, by using the above argument, it is easy to check that the map $\theta:\bar{a}\mapsto a^{x_{11}},\bar{b}\mapsto a^{x_{21}p^{n-m}}b^{x_{22}}$ is an isomorphism from
$\bar{G}$ onto $G$.
 \qed

\begin{thm}\label{th=5.8}
 Let $G$ be a finite $p$-group such that $G_3\cong C_p$, $\Phi(G') G_3\cong C_p^2$, $G/\Phi(G')G_3\cong M_p(n,m,1)$ where $n>m$, and the type of $C_G(G')/\Phi(G')G_3$ is $(p^{n-1},p^{m},p)$. Then $G$ is one of the following non-isomorphic groups:

 {\rm (M1)} $\langle a,b,c\di a^{p^{n+1}}=b^{p^{m+1}}=1,
[a,b]=c,c^p=a^{p^n}b^{s\nu p^m},[a,c]=b^{\nu p^m},[b,c]=\mbox{\hskip0.65in}1\rangle$, where $n>m\ge 2$, $\nu=1$ or a fixed quadratic non-residue modular $p$, \mbox{\hskip0.65in}$s=0,1,\dots,\frac{p-1}{2}$;%(??);%(??)

 {\rm (M2)} $\langle a, b, c \di a^{p^{n+1}}=b^{p^{m}}=c^{p^2}=1,
[a,b]=c, [c,b]=1,[c,a]=c^{p}a^{-p^n},[a^{p^n},b]=\mbox{\hskip0.65in}1\rangle$, where $n>m\ge 2$;%(??)

 {\rm (M3)} $\langle a, b, c,d \di a^{p^{n+1}}=b^{p^{m}}=d^p=1,
[a,b]=c, c^p=a^{p^n}, [c,a]=d,[c,b]=\mbox{\hskip0.65in}1, [d,a]=[d,b]=1\rangle$, where $n>m\ge 2$;%(??)

{\rm (M4)} $\langle a, b, c \di a^{p^{n+1}}=b^{p^{m+1}}=1,
[a,b]=c,c^p=b^{p^m}, [b,c]=1,[a,c]=a^{p^n}\rangle$, where \mbox{\hskip0.65in}$n>m\ge 2$;%(??)

{\rm (M5)} $\langle a, b, c \di a^{p^{n+1}}=b^{p^{m}}=c^{p^2}=1,
[a,b]=c, [b,c]=1,[a,c]=a^{p^n}\rangle$, where \mbox{\hskip0.65in}$n>m\ge 2$;%(??)

 {\rm (M6)} $\langle a, b, c \di a^{p^{n}}=b^{p^{m+1}}=c^{p^2}=1,
[a,b]=c, [c,b]=1,[c,a]=c^{tp}b^{-tp^m},[b^{p^m},a]=\mbox{\hskip0.65in}1\rangle$, where $n>m\ge 2$ and $t\in F_p^*$;%(??)

{\rm (M7)} $\langle a, b, c,d \di a^{p^{n}}=b^{p^{m+1}}=d^p=1, c^p=b^{p^m},
[a,b]=c, [c,b]=1,[c,a]=\mbox{\hskip0.65in}d,[d,a]=[d,b]=1\rangle$, where $n>m\ge 2$;%(??)

 {\rm (M8)} $\langle a, b, c \di a^{p^{n}}=b^{p^{m+1}}=c^{p^2}=1,
[a,b]=c, [c,b]=1,[a,c]=b^{\nu p^m}\rangle$, where \mbox{\hskip0.65in}$n>m\ge 2$, $\nu=1$ or a fixed quadratic non-residue modular $p$;%(??)

{\rm (M9)} $\langle a, b, c,d \di a^{p^{n}}=b^{p^{m}}=c^{p^2}=d^p=1,
[a,b]=c, [c,b]=1,[c,a]=d,[d,a]=\mbox{\hskip0.65in}[d,b]=1\rangle$, where $n>m\ge 2$.%(??)

 \end{thm}
\demo
If $m=1$, then $1=[a,b^{p}]=c^p$. It follows that $\Phi(G')=1$, a contradiction. Hence $m\ge 2$.
 Suppose that $G$ and $\bar{G}$ are two groups described in the theorem. By Theorem \ref{isomorphic-6.2.1}, $G\cong \bar{G}$ if and only if there exists $X={\left(
 \begin{array}{cc}
 x_{11}& 0\\
 x_{21}& x_{22}
 \end{array}
 \right)}$, an invertible matrix over $F_p$, such that
 \begin{equation}\label{eq: 5.12}
 w(\bar{G})=Xw(G)\diag( x_{11}^{-1}x_{22}^{-1},x_{11}^{-2}x_{22}^{-1}).
 \end{equation}
 By suitably choosing $x_{21}$, that is, using an elementary row operation, we can simplify $w(G)$ to be one of the following types:
(a)  $\left(
 \begin{array}{cc}
 w_{11}& w_{12}\\
 0& w_{22}
 \end{array}
 \right)$ where $w_{11}\neq 0$,
(b)  $\left(
 \begin{array}{cc}
 0  & w_{12}\\
 w_{21}& 0
 \end{array}
 \right)$ where $w_{12}\neq 0$, and
(c)  $\left(
 \begin{array}{cc}
 0&0\\
 w_{21}&w_{22}
 \end{array}
 \right).$
In the following, we assume that both $w(G)$ and
$w(\bar{G})$ are such matrices. By Equation \ref{eq: 5.12}, it is easy to check that (i)
 different types give non-isomorphic groups, (ii) $G\cong \bar{G}$ if and only if there exists $X=\diag(x_{11},x_{22})$, an invertible matrix over $F_p$, such that
$ w(\bar{G})=Xw(G)\diag(x_{11}^{-1}x_{22}^{-1},x_{11}^{-2}x_{22}^{-1})$.
\begin{table}[h]
  \centering
%\caption{Properties of the groups of Type (A) in Theorem
%\ref{main}}\label{table: 8.1}
{
\tiny
%\scriptsize
\begin{tabular}{cccccc}
w(G)&Case &$X$ & $w(\bar{G})$&Group&Remark\\
   \hline
(a)&$w_{22}\neq 0$ & $
 \begin{array}{c}
 \diag(z,w_{11})\\
 \rm{ where}\ w_{22}=\nu z^2\\
 \end{array}$ &${\left(
 \begin{array}{cc}
 1& w_{12}w_{11}^{-1}z^{-1}\\
 0& \nu
 \end{array}
 \right)}$&(M1)& $s=-w_{12}w_{11}^{-1}z^{-1}$\\

(a) & $w_{22}=0$, $w_{12}\neq 0$ & $\diag(w_{11}^{-1}w_{12},w_{11})$ &${\left(
 \begin{array}{cc}
 1& 1\\
 0& 0
 \end{array}
 \right)}$ &(M2)&\\

(a) & $w_{22}=w_{12}=0$ & $\diag(1,w_{11})$ &${\left(
 \begin{array}{cc}
 1& 0\\
 0& 0
 \end{array}
 \right)}$ &(M3)&\\
\hline
(b) &$w_{21}\neq 0$ & $\diag(w_{21},w_{21}^{-1}w_{12})$ &${\left(
 \begin{array}{cc}
 0& 1\\
 1& 0
 \end{array}
 \right)}$ &(M4)&\\
(b) &$w_{21}=0$ & $\diag(1,w_{12})$ &${\left(
 \begin{array}{cc}
 0& 1\\
 0& 0
 \end{array}
 \right)}$&(M5)&\\
\hline
(c)& $w_{21}\neq 0$ & $\diag(w_{21},1)$ &${\left(
 \begin{array}{cc}
 0& 0\\
 1& w_{22}w_{21}^{-2}
 \end{array}
 \right)}$ & $\begin{array}{c} {\rm (M6)\ if\ }w_{22}\neq 0\\ {\rm (M7)\ if\ }w_{22}=0\end{array}$ &$\begin{array}{c} t=w_{22}^{-1}w_{21}^2\\ \\ \end{array}$ \\
(c)&$w_{21}=0$, $w_{22}\neq 0$ & $
 \begin{array}{c}
 \diag(1,z)\\
 \rm{ where}\ w_{22}=\nu z^2\\
 \end{array}$ &${\left(
 \begin{array}{cc}
 0& 0\\
 0& \nu
 \end{array}
 \right)}$ &(M8) &\\
(c) &$w_{21}=w_{22}=0$ &  &${\left(
 \begin{array}{cc}
 0& 0\\
 0& 0
 \end{array}
 \right)}$&(M9)&\\
   \hline
\end{tabular}
\caption{The types in Theorem \ref{th=5.8}}}
\label{table 5}
\end{table}
By Table \ref{table 5}, we get the groups of Type (M1)--(M9).
\qed
\section{The case  $\Phi(G')\le G_3\cong C_p^2$}

Suppose that $G$ is a finite $p$-group with $\Phi(G')\le G_3\cong C_p^2$, $G_3\le Z(G)$ and $G/G_3\cong M_p(n,m,1)$, where $n>1$ for $p=2$ and
$n\ge m$. Let $$G/G_3=\lg \bar{a},\bar{b},\bar{c}\mid \bar{a}^{p^n}=\bar{b}^{p^m}=\bar{c}^{p}=1,[\bar{c},\bar{a}]=[\bar{b},\bar{c}]=1\rg.$$
Then, without loss of generality, we may assume that $G=\lg a,b,c\rg$ where $[a,b]=c$.
Let $x=[b,c], y=[c,a]$. Then $G_3=\lg x,y\rg$. Since $a^{p^{n}}\in G_3$, we may assume that $a^{p^{n}}=x^{w_{11}}y^{w_{12}}$. By similar reasons, we may assume that $b^{p^m}=x^{w_{21}}y^{w_{22}}$ and $c^p=x^{w_{31}}y^{w_{32}}$. Let $w(G)=\left(
 \begin{array}{cc}
 {w}_{11}& {w}_{12}\\
 {w}_{21}& {w}_{22}
  \end{array}
 \right)$ and $v(G)=\left(
 \begin{array}{c}
 {w}_{31}\\
 {w}_{32}
  \end{array}
 \right)$. Then we get two matrices over $F_p$. $w(G)$ is called a characteristic matrix of $G$, and $v(G)$ is called a characteristic vector of $G$. Notice that $w(G)$ and $v(G)$ will be changed if we change the generators $a,b$. We also call $a,b$ a set of characteristic generators of $w(G)$ and $v(G)$.

\begin{thm}
\label{property-5.0} Suppose that $G$ is a finite $p$-group such that
$\Phi(G')\le G_3\cong C_p^2$, $G_3\le Z(G)$ and $G/G_3\cong M_p(n,m,1)$. Let $w(G)=(w_{ij})$ be a characteristic matrix of $G$. Then the following conclusions hold:

{\rm (1)} If $I_{\min}=1$, then $m=1$;

{\rm (2)} If $I_{\max}=2$, then $n\le 2$;

{\rm (3)} If $p=2$ and $m=1$, then $I_{\min}=1$;

{\rm (4)} If $p>3$ and $n=m=1$, then $I_{\min}\neq 1$ if and only if $w_{11}=w_{22}=w_{12}+w_{21}=0$;

{\rm (5)} If $p>2$ and $n>m=1$, then $I_{\min}\neq 1$ if and only if $w_{11}=w_{12}=0$;

{\rm (6)} If $p>3$ and $n=m=1$, then $I_{\max}\neq 2$ if and only if $(w_{12}+w_{21})^2-4w_{11}w_{22}$ is a quadratic non-residue modular $p$;

{\rm (7)} If $p>2$, $m=1$ and $n=2$, then $I_{\max}=2$ if and only if $w_{22}\neq 0$;

{\rm (8)} If $p>2$, $n=m=2$ and $\Phi(G')=1$, then $I_{\max}=2$ if and only if $(w_{12}+w_{21})^2-4w_{11}w_{22}$ is a quadratic non-residue modular $p$;

{\rm (9)} If $p>2$, $n=m=2$ and $(w_{31},w_{32})=(0,1)$, then $I_{\max}=2$ if and only if $w_{11}\neq 0$ and one of the following holds: {\rm (i)} $(w_{12}+w_{21})^2-4w_{11}(w_{22}+1)$ is a quadratic non-residue modular $p$; {\rm (ii)} $w_{21}=w_{11}w_{22}+w_{11}$ and $(w_{12}+w_{21})^2=4w_{21}$.
\end{thm}
\demo Let $N=\lg b,a^p,c,x,y\rg$ and $M_i=\lg ab^i,b^p,c,x,y\rg$ where $0\le i\le p-1$. The $N$ and $M_i$ are all maximal subgroups of $G$.

(1) Assume that $H$ is an $\mathcal{A}_1$-subgroup of index $p$. Then $d(H/H')=2$. Since $G/G_3\in\mathcal{A}_1$, $H/G_3$ is abelian and hence $H'\le G_3$. It follows that $d(H/G_3)=2$. Since $\Phi(G)\le H$, $d(\Phi(G)/G_3)\le 2$.
Since $\Phi(G)/G_3=\lg \bar{a}^p,\bar{b}^p,\bar{c}\rg$ is of type $(p^{n-1},p^{m-1},p)$, we have $m=1$.

(2) Let $K=\lg b,c\rg$. Then $K\in\mathcal{A}_1$ and $|G:K|\ge |G: KG_3|=p^n$. Since $I_{\max}=2$, $n\le 2$.

(3) If $m=1$ and $p=2$, then $1=[a,b^2]=c^2x$. Hence $c^2=x$. In this case, $M_1=\lg ab,c\rg\in \mathcal{A}_1$. Thus $I_{\min}=1$.

%If $I_{\max}=2$, then $N\in \mathcal{A}_2$. It follows that $n=2$. Let $A=\lg a^2,c,x,y\rg$, $B_r=\lg bc^r,a^2,x,y\rg$ and $C_{st}=\lg ba^{2s},ca^{2t},x,y\rg$ where $r,s,t=0,1$. Then $A$, $B_r$ and $C_{st}$ are all maximal subgroups of $N$.

%First of all, $A$ is the unique abelian maximal subgroup of $N$.
%Since $[bc^r,a^2]=xy$, $B_r\in \mathcal{A}_1$. It follows that $w_{11}+w_{12}\neq 0$ or $w_{21}+w_{22}\neq 0$.
%Since $[ba^{2s},ca^{2t}]=x(xy)^{t}$, $(ba^{2s})^{2}=x^{sw_{11}+sw_{12}+w_{21}+w_{22}}(xy)^{sw_{12}+w_{22}+s}$ and $(ca^{2t})^2=x^{tw_{11}+tw_{12}+1}(xy)^{tw_{12}}$, $C_{st}\in \mathcal{A}_1$ if and only if the following equation set about $s$ and $t$ has no solution.
%\begin{equation}
 % \left\{
  %  sw_{12}+w_{22}+s &=t(sw_{11}+sw_{12}+w_{21}+w_{22}) &&(5.1.1)\\
   %  tw_{12} &=t(tw_{11}+tw_{12}+1) &&(5.1.2)\\
    % \end{aligned}
     %\right.
     %\end{equation}
%If $w_{12}=0$, then $s=w_{22}, t=0$ is a solution of Equation set (5.1). Hence $w_{12}=1$. If $w_{11}=1$, then $w_{11}+w_{12}=0$ and hence $w_{21}+w_{22}=1$. In this case, $s=t=w_{22}$ is a solution of Equation set (5.1). Hence $w_{11}=0$. In this case, Equation set (5.1) has no solution if and only if $w_{22}=1$.
%
%Conversely, if $w_{11}=0$ and $w_{12}=w_{22}=1$, then it is easy see that $G\in\mathcal{A}_3$ and $I_{\max}=2$.\qed
(4) If $n=m=1$ and $p>3$, then $|G|=p^5$ and $1=[a,b^p]=c^p$. If $I_{\min}\neq 1$, then $N\in\mathcal{A}_2$ and $M_i\in \mathcal{A}_2$. Since $y\not\in \lg b,c\rg$, $w_{22}=0$. By calculation, $[ab^i,c]=y^{-1}x^{i}$ and $(ab^i)^p=a^pb^{ip}=x^{w_{11}+iw_{21}}y^{w_{12}}$. Hence, $\forall i$, $w_{11}+iw_{21}=-iw_{12}$. It follows that $w_{11}=w_{12}+w_{21}=0$.

Conversely, if $w_{11}=w_{22}=w_{12}+w_{21}=0$, then, by using above argument, $N\in\mathcal{A}_2$ and $M_i\in \mathcal{A}_2$. Hence $I_{\min}\neq 1$.

(5) If $n>m=1$ and $p>2$, then $1=[a,b^p]=c^p$. If $I_{\min}\neq 1$, then $M_i\not\in \mathcal{A}_1$. By calculation, $[ab^i,c]=y^{-1}x^{i}$ and $(ab^i)^{p^n}=a^{p^n}=x^{w_{11}}y^{w_{12}}$. Hence, $\forall i$, $w_{11}=-iw_{12}$. It follows that $w_{11}=w_{12}=0$.

Conversely, if $w_{11}=w_{12}=0$, then, by using the above argument, $N\not\in\mathcal{A}_1$ and $M_i\not\in \mathcal{A}_1$. Hence $I_{\min}\neq 1$.

(6) If $n=m=1$ and $p>3$, then $1=[a,b^p]=c^p$. If $I_{\max}\neq 2$, then $N\in\mathcal{A}_1$ and $M_i\in \mathcal{A}_1$. Since $y\in \lg b,c\rg$, $w_{22}\neq 0$. By calculation, $[ab^i,c]=y^{-1}x^{i}$ and $(ab^i)^p=a^pb^{ip}=x^{w_{11}+iw_{21}}y^{w_{12}+iw_{22}}$. Hence, $\forall i$, $w_{11}+iw_{21}\neq-iw_{12}-i^2w_{22}$. It follows that $(w_{12}+w_{21})^2-4w_{11}w_{22}$ is a quadratic non-residue modular $p$.

Conversely, if $(w_{12}+w_{21})^2-4w_{11}w_{22}\not \in F_p^2$, then, by using the above argument, $N\in\mathcal{A}_1$ and $M_i\in \mathcal{A}_1$. Hence $I_{\max}\neq 2$.

(7)  If $p>2$, $n=2$ and $m=1$, then $1=[a,b^p]=c^p$. By calculation, $M_i\in \mathcal{A}_1$ or $\mathcal{A}_2$. If $I_{\max}=2$, then
$N\in\mathcal{A}_2$. It follows that $y\in\lg b,c\rg$. Hence $w_{22}\neq 0$.

Conversely, if $w_{22}\neq 0$, then $N\in\mathcal{A}_2$ and $M_i\in \mathcal{A}_1$ or $\mathcal{A}_2$. Hence $I_{\max}=2$.

(8) If $n=m=2$, then, by (1), $I_{\min}\ge 2$. If $I_{\max}=2$, then $N\in\mathcal{A}_2$ and $M_i\in \mathcal{A}_2$.
  Since $|G:\lg b,c\rg|=p^2$, we have $w_{22}\neq 0$. By calculation, $[c,ab^i]=x^{-i}y$ and $(ab^i)^{p^2}=x^{w_{11}+iw_{21}}y^{w_{12}+iw_{22}}$. Since $|G:\lg c,ab^i\rg|=p^2$, Equation ${w_{11}+iw_{21}}=-i({w_{12}+iw_{22}})$ about $i$ has no solution. Hence $(w_{12}+w_{21})^2-4w_{11}w_{22}$ is a quadratic non-residue modular $p$.

Conversely, if $(w_{12}+w_{21})^2-4w_{11}w_{22}$ is a quadratic non-residue modular $p$, then, by the above argument, $N\in\mathcal{A}_2$ and $M_i\not\in \mathcal{A}_2$. Hence $I_{\max}=2$.

(9) If $n=m=2$, then, by (1), $I_{\min}\ge 2$. If $I_{\max}=2$, then $N\in\mathcal{A}_2$ and $M_i\in \mathcal{A}_2$.
Let $N_j=\lg a^jb, a^p,c,x,y\rg$. Then $N_0=N$ and $N_j=M_{j^{-1}}$ for $j\neq 0$.
  Since $M_0\in \mathcal{A}_2$, $|G:\lg a,c\rg|=p^2$ and hence $w_{11}\neq 0$. Let $A=\lg a^p,b^p,c,x\rg$, $B_r=\lg a^jbc^r,a^p,c^p,x\rg$ and $C_{st}=\lg a^jba^{sp},ca^{tp},c^p,x\rg$ where $0\le r,s,t\le p-1$. Then $A$, $B_r$ and $C_{st}$ are all maximal subgroups of $N_j$.

First of all, $A$ is the unique abelian maximal subgroup of $N_j$. Since $w_{11}\neq 0$, $B_r\in \mathcal{A}_1$.
Since $[ca^{tp},a^jba^{sp}]=x^{-1}y^{t+j}$ and $(a^jba^{sp})^{p^2}=x^{jw_{11}+w_{21}}y^{jw_{12}+w_{22}}$ and $(ca^{tp})^p=x^{tw_{11}}y^{tw_{12}+1}$, $C_{st}\in \mathcal{A}_1$ if and only if the following equation set about $j$ and $t$ has no solution.
\begin{equation}
  \left\{
  \begin{aligned}
     jw_{12}+w_{22} &=-(j+t)(jw_{11}+w_{21}) &&(6.1.1)\\
     tw_{12}+1 &=-(j+t)tw_{11} &&(6.1.2)\\
     \end{aligned}
     \right.
     \end{equation}
 Adding up (6.1.1) and (6.1.2), we have
  \begin{equation*}
   w_{11}(j+t)^2+(w_{12}+w_{21})(j+t)+w_{22}+1=0 \eqno(6.1.3)
  \end{equation*}

  If Equation (6.1.3) about $j+t$ has no solution, then $(w_{12}+w_{21})^2-4w_{11}(w_{22}+1)$ is a quadratic non-residue modular $p$. Hence (i) holds.

  If Equation (6.1.3) about $j+t$ has a unique solution, then $(w_{12}+w_{21})^2-4w_{11}(w_{22}+1)=0$.
  In this case, $C_{st}\in\mathcal{A}_1$ if and only if Equation (6.1.2) about $t$ has no solution.
That is, $$j+t=-w_{11}^{-1}w_{12}\eqno(6.1.4)$$
By (6.1.3) and (6.1.4), $w_{21}=w_{11}(w_{22}+1)$ and hence (ii) holds.

  If Equation (6.1.3) has two solutions, then there exists a solution such that $w_{11}(j+t)+w_{12}\neq 0$. It follows that Equation set (6.1) has a solution, and hence there exists $j,t$ such that $C_{st}\not\in \mathcal{A}_1$.

 Conversely, if (i) or (ii) holds, then, by using the above argument, it is easy to see that Equation set (6.1) has no solution.
   Hence $C_{st}\in\mathcal{A}_1$ and $I_{\max}=2$.
\qed

\begin{thm}\label{isomorphic-5.1}
Suppose that $G$ and $\bar{G}$ are finite $p$-groups such that $\Phi(G')\le G_3\cong C_p^2$, $G_3\le Z(G)$ and $G/G_3\cong M_p(n,m,1)$, where $p>2$, $n\ge m\ge 2$. Let two characteristic matrices of $G$ and $\bar{G}$ be $w(G)=(w_{ij})$ and $w(\bar{G})=(\bar{w}_{ij})$ respectively.
Then $G\cong \bar{G}$ if and only if there exists $Y={\left(
 \begin{array}{cc}
 y_{11}& y_{12}\\
 y_{21}p^{n-m}& y_{22}
 \end{array}
 \right)}$, an invertible matrix over $F_p$, such that $
w(\bar{G})=Y_1w(G)Y^t$ and $v(\bar{G})
=Yv(G)$
, where
$Y_1={\left(
 \begin{array}{cc}
 y_{11}& y_{12}p^{n-m}\\
 y_{21}& y_{22}
 \end{array}
 \right)}$.
 \end{thm}
\demo Suppose that $a,b$ and
$\bar{a},\bar{b}$ are two set of characteristic generators of $w(G)$, $v(G)$ and $w(\bar{G})$, $v(\bar{G})$.
Let $\theta$ be an isomorphism from $\bar{G}$ onto $G$. We have $\Phi(\bar{G})^\theta=\Phi(G)$ and $\Omega_m(\bar{G})^\theta=\Omega_m(G)$ since these four subgroups are characteristic in $G$ or $\bar{G}$. So we may let
$$\bar{a}^\theta= a^{x_{11}}b^{x_{12}}\phi_1,\
\bar{b}^\theta=a^{x_{21}p^{n-m}}b^{x_{22}}\phi_2$$
where $\phi_1\in \Phi(G)$, $\phi_2\in \Phi(G)\cap \Omega_m(G)$ and $X:=\left(
 \begin{array}{cc}
 x_{11}& x_{12}\\
 x_{21}p^{n-m}& x_{22}
 \end{array}
 \right)$ is an invertible matrix over $F_p$.
By calculations,
 we have
$$\bar{c}^\theta=[\bar{a},\bar{b}]^\theta=[\bar{a}^\theta,\bar{b}^\theta]\equiv [a^{x_{11}}b^{x_{12}},a^{x_{21}p^{n-m}}b^{x_{22}}]\equiv
c^{|X|}\ (\mod G_3),$$ \begin{equation*}\bar{x}^\theta=[\bar{b},\bar{c}]^\theta=[\bar{b}^\theta,\bar{c}^\theta]=[a^{x_{21}p^{n-m}}b^{x_{22}},c^{|X|}]=x^{|X|x_{22}}y^{-|X|x_{21}p^{n-m}},
\end{equation*}
and
\begin{equation*}\bar{y}^\theta=[\bar{c},\bar{a}]^\theta=[\bar{c}^\theta,\bar{a}^\theta]=[c^{|X|},a^{x_{11}}b^{x_{12}}]=x^{-|X|x_{12}}y^{|X|x_{11}}.
\end{equation*}
 By transforming $\bar{x}^{\bar{w}_{11}}\bar{y}^{\bar{w}_{12}}=\bar{a}^{p^{n}}$ by $\theta$, we have
\begin{equation}\label{eq 5.1}
(\bar{w}_{11},\bar{w}_{12})\left(
 \begin{array}{cc}
 |X|x_{22}& -|X|x_{21}p^{n-m}\\
 -|X|x_{12}& |X|x_{11}
 \end{array}
 \right)=
 (x_{11},x_{12}p^{n-m})
 \left(
 \begin{array}{cc}
 w_{11}& w_{12}\\
 w_{21}& w_{22}
 \end{array}
 \right).
\end{equation}
 By transforming $\bar{x}^{\bar{w}_{21}}\bar{y}^{\bar{w}_{22}}=\bar{b}^{p^{m}}$ by $\theta$, we have
\begin{equation}\label{eq 5.2}
(\bar{w}_{21},\bar{w}_{22})\left(
 \begin{array}{cc}
 |X|x_{22}& -|X|x_{21}p^{n-m}\\
 -|X|x_{12}& |X|x_{11}
 \end{array}
 \right)=
 (x_{21},x_{22})
 \left(
 \begin{array}{cc}
 w_{11}& w_{12}\\
 w_{21}& w_{22}
 \end{array}
 \right).
\end{equation}
By Equation \ref{eq 5.1} and \ref{eq 5.2},
\begin{equation}\label{eq 5.3}
|X|\left(
 \begin{array}{cc}
 \bar{w}_{11}& \bar{w}_{12}\\
 \bar{w}_{21}& \bar{w}_{22}
  \end{array}
 \right)\left(
 \begin{array}{cc}
 x_{22}& -x_{21}p^{n-m}\\
 -x_{12}& x_{11}
 \end{array}
 \right)=
 \left(
 \begin{array}{cc}
 x_{11}& x_{12}p^{n-m}\\
 x_{21}& x_{22}
 \end{array}
 \right)
 \left(
 \begin{array}{cc}
 w_{11}& w_{12}\\
 w_{21}& w_{22}
 \end{array}
 \right)
\end{equation}
Let $$Y=|X|^{-1}X=|X|^{-1}\left(
 \begin{array}{cc}
 x_{11}& x_{12}\\
 x_{21}p^{n-m}& x_{22}
 \end{array}
 \right)=\left(
 \begin{array}{cc}
 y_{11}& y_{12}\\
 y_{21}p^{n-m}& y_{22}
 \end{array}
 \right)
$$ and $$Y_1=|X|^{-1}\left(
 \begin{array}{cc}
 x_{11}& x_{12}p^{n-m}\\
 x_{21}& x_{22}
 \end{array}
 \right)
 =\left(
 \begin{array}{cc}
 y_{11}& y_{12}p^{n-m}\\
 y_{21}& y_{22}
 \end{array}
 \right).
$$ 
Right multiplying $Y^{t}$, we have
\begin{equation}\label{eq 5.4}
\left(
 \begin{array}{cc}
 \bar{w}_{11}& \bar{w}_{12}\\
 \bar{w}_{21}& \bar{w}_{22}
  \end{array}
 \right)=Y_1\left(
 \begin{array}{cc}
 w_{11}& w_{12}\\
 w_{21}& w_{22}
 \end{array}
 \right)Y^t.
\end{equation}
 By transforming $\bar{x}^{\bar{w}_{31}}\bar{y}^{\bar{w}_{32}}=\bar{c}^{p}$ by $\theta$, we have \begin{equation}\label{eq 5.5}
\left(
 \begin{array}{cc}
 |X|x_{22}&  -|X|x_{12}\\
-|X|x_{21}p^{n-m}& |X|x_{11}
 \end{array}
 \right)\left(
 \begin{array}{c}
 \bar{w}_{31}\\
\bar{w}_{32}
 \end{array}\right)
=\left(
 \begin{array}{c}
 |X|w_{31}\\
|X|w_{32}
 \end{array}\right).
\end{equation}
Left multiplying $Y$, we have
\begin{equation}\label{eq 5.6}
\left(
 \begin{array}{c}
 \bar{w}_{31}\\
\bar{w}_{32}
 \end{array}\right)
=Y\left(
 \begin{array}{c}
 w_{31}\\
w_{32}
 \end{array}\right).
 \end{equation}

Conversely, if there exists an invertible matrix $Y=\left(
 \begin{array}{cc}
 y_{11}& y_{12}\\
 y_{21}p^{n-m}& y_{22}
 \end{array}
 \right)$ over
$F_p$ such that equation \ref{eq 5.4} and \ref{eq 5.6}, then, letting $X=|Y|^{-1}Y=\left(
 \begin{array}{cc}
 x_{11}& x_{12}\\
 x_{21}p^{n-m}& x_{22}
 \end{array}
 \right)$, by using the above argument, it is easy to check that the map $\theta:\bar{a}\mapsto a^{x_{11}}b^{x_{12}},\bar{b}\mapsto a^{x_{21}p^{n-m}}b^{x_{22}}$ is an isomorphism from
$\bar{G}$ onto $G$.
 \qed

\medskip

If $p>2$ and $m=1$, then $c^p=[a,b]^p=[a,b^p]=1$ and hence $\Phi(G')=1$.
In addition, if $p>3$ or $n>1$, then we also have Equation \ref{eq 5.4} and Equation \ref{eq 5.6}. Thus we get the following theorem.

\begin{thm}\label{isomorphic-5.2}
Suppose that $G$ and $\bar{G}$ are finite $p$-groups such that $\Phi(G')\le G_3\cong C_p^2$, $G_3\le Z(G)$ and $G/G_3\cong M_p(n,1,1)$, where $p>2$ and $n>1$ for $p=3$. Let two characteristic matrices of $G$ and $\bar{G}$ be $w(G)=(w_{ij})$ and $w(\bar{G})=(\bar{w}_{ij})$ respectively. Then $\Phi(G')=1$, and $G\cong \bar{G}$ if and only if there exists $Y={\left(
 \begin{array}{cc}
 y_{11}& y_{12}\\
 y_{21}p^{n-m}& y_{22}
 \end{array}
 \right)}$, an invertible matrix over $F_p$, such that $w(\bar{G})=Y_1w(G)Y^t$, where
$Y_1={\left(
 \begin{array}{cc}
 y_{11}& y_{12}p^{n-m}\\
 y_{21}& y_{22}
 \end{array}
 \right)}$.
 \end{thm}

\medskip

If $p=2$, $m\ge 2$ and $n\ge 3$, then we also have Equation \ref{eq 5.4} and \ref{eq 5.6}.
If $p=2$, $m=1$ and $n\ge 3$, then $b^2\in Z(G)$. It follows that $1=[a,b^2]=c^2[c,b]$. Hence $[c,b]=c^2$. That is, $v(G)=\left(
 \begin{array}{c}
 1\\
 0
 \end{array}
 \right)$. In this case, we may let $\bar{b}^\theta=a^{x_{21}2^{n-1}}bc^{x_{23}}$. By calculations, we have $(\bar{b}^{p^{m}})^\theta=(\bar{b}^{2})^\theta=(a^{x_{21}2^{n-1}}bc^{x_{23}})^2=a^{x_{21}2^n}b^2c^{2x_{23}}[b,c]^{2x_{23}}=a^{x_{21}2^n}b^2$. Hence we also have Equation \ref{eq 5.2}. It follows that Equation \ref{eq 5.4} and \ref{eq 5.6} hold.
Thus we get the following theorem.

\begin{thm}\label{isomorphic-5.3}
Suppose that $G$ and $\bar{G}$ are finite $2$-groups such that $\Phi(G')\le G_3\cong C_2^2$, $G_3\le Z(G)$ and $G/G_3\cong M_p(n,m,1)$, where $n\ge 3$ and $n\ge m$. Let two characteristic matrices of $G$ and $\bar{G}$ be $w(G)=(w_{ij})$ and $w(\bar{G})=(\bar{w}_{ij})$ respectively.
Then $G\cong \bar{G}$ if and only if there exists $Y={\left(
 \begin{array}{cc}
 y_{11}& y_{12}\\
 y_{21}2^{n-m}& y_{22}
 \end{array}
 \right)}$, an invertible matrix over $F_2$, such that $
w(\bar{G})=Y_1w(G)Y^t$ and $v(\bar{G})
=Yv(G)$, where
$Y_1={\left(
 \begin{array}{cc}
 y_{11}& y_{12}2^{n-m}\\
 y_{21}& y_{22}
 \end{array}
 \right)}$. In addition, if $m=1$, then $v(G)=\left(
 \begin{array}{c}
 1\\
 0
 \end{array}
 \right)$.
 \end{thm}

\begin{thm}
\label{5}Let $G$ be a finite $p$-group such that
$\Phi(G')\le G_3\cong C_p^2$, $G_3\le Z(G)$ and $G/G_3\cong M_p(n,m,1)$, where $n>1$ for $p=2$ and
$n\ge m$. Then $G$ is one of the following non-isomorphic groups:

{\rm (N1)} $\langle a, b, c \di a^{8}=b^{8}=c^{2}=1,
[a,b]=c,[c,a]=b^{4},[c,b]=a^{4},[a^4,b]=1\rangle$;%(F?)

{\rm (N2)} $\langle a, b, c \di a^{8}=b^{8}=c^{2}=1,
[a,b]=c,[c,a]=a^{4},[c,b]=b^{4}\rangle$;%(F1)

{\rm (N3)} $\langle a, b, c \di a^{8}=b^{8}=c^{2}=1,
[a,b]=c,[c,a]=a^{4}b^{4},[c,b]=a^{4},[a^4,b]=1\rangle$;%(F?)

{\rm (N4)} $\langle a, b, c,d \di a^{8}=b^{4}=c^{2}=d^2=1,
[a,b]=c,[c,a]=a^{4},[c,b]=d,[d,a]=\mbox{\hskip0.6in}[d,b]=1\rangle$;%(G4)

{\rm (N5)} $\langle a, b, c,d \di a^{4}=b^{8}=c^{2}=d^2=1,
[a,b]=c,[c,a]=b^{4},[c,b]=d,[d,a]=\mbox{\hskip0.6in}[d,b]=1\rangle$;%(G3)

{\rm (N6)} $\langle a, b, c,d,e \di a^{4}=b^{4}=c^{2}=d^2=e^2=1,
[a,b]=c,[c,a]=d,[c,b]=e,[d,a]=\mbox{\hskip0.6in}[d,b]=[e,a]=[e,b]=1\rangle$;%(H)

{\rm (N7)} $\langle a, b, c \di a^{8}=b^{4}=c^{4}=1,
[a,b]=c,[c,a]=c^2,[c,b]=a^4\rangle$;%??

{\rm (N8)} $\langle a, b, c \di a^{8}=b^{8}=1,
[a,b]=c,[c,a]=c^2=b^4,[c,b]=a^4\rangle$;% I_{\max}=2

{\rm (N9)} $\langle a, b, c \di a^{8}=c^4=1,
[a,b]=c,[c,a]=c^2,[c,b]=a^4=b^4\rangle$;% I_{\max}=2

{\rm (N10)} $\langle a, b, c \di a^{8}=b^{8}=1,
[a,b]=c,[c,a]=c^2=a^4b^4,[c,b]=a^4\rangle$;%

{\rm (N11)} $\langle a, b, c \di a^{8}=b^{8}=1,
[a,b]=c,[c,a]=c^2=a^4,[c,b]=b^4\rangle$;%

{\rm (N12)} $\langle a, b, c \di a^{4}=b^{8}=c^4=1,
[a,b]=c,[c,a]=c^2,[c,b]=b^4\rangle$;%

{\rm (N13)} $\langle a, b, c \di a^{4}=b^{8}=c^4=1,
[a,b]=c,[c,a]=c^2,[c,b]=c^2b^4,[c^2,b]=[b^4,a]=\mbox{\hskip0.65in}1\rangle$;%

{\rm (N14)} $\langle a, b, c,d \di a^{4}=b^{4}=c^4=d^2=1,
[a,b]=c,[c,a]=c^2,[c,b]=d,[d,a]=\mbox{\hskip0.65in}[d,b]=1\rangle$;%

{\rm (N15)} $\langle a, b, c,d \di a^{8}=b^{4}=d^2=1,
[a,b]=c,[c,a]=c^2=a^4,[c,b]=d,[d,a]=\mbox{\hskip0.65in}[d,b]=1\rangle$;%

{\rm (N16)} $\langle a, b, c,d \di a^{8}=d^2=1,
[a,b]=c,[c,a]=c^2=a^4=b^4,[c,b]=d,[d,a]=\mbox{\hskip0.65in}[d,b]=1\rangle$;%

{\rm (O1)}  $\langle a, b, c,d,e \di a^{3}=b^{3}=c^3=d^3=e^3=1,
[a,b]=c,[c,a]=d,[c,b]=e,[d,a]=\mbox{\hskip0.6in}[d,b]=[e,a]=[e,b]=1
\rangle$;%(E1)

{\rm (O2)}  $\langle a, b, c,d \di a^{3}=b^{9}=c^3=d^3=1,
[a,b]=c,[c,a]=d,[c,b]=b^3,[d,a]=[d,b]=\mbox{\hskip0.6in}1 \rangle$;%(E2)

{\rm (O3)}  $\langle a, b, c,d\di a^{9}=c^{3}=d^3=1, b^3=a^3,
[a,b]=c,[c,a]=d,[c,b]=a^{3},[d,a]=\mbox{\hskip0.6in}[d,b]=1
\rangle$;%(E3)

{\rm (O4)}  $\langle a, b, c,d \di a^{9}=b^{3}=c^3=d^3=1,
[a,b]=c,[c,a]=d,[c,b]=a^{-3},[d,a]=\mbox{\hskip0.6in}[d,b]=1
\rangle$;%(E4)

{\rm (O5)}  $\langle a, b, c \di a^{9}=b^9=c^3=1,
[a,b]=c,[c,a]=a^3,[c,b]=b^{3} \rangle$;%(E5)

{\rm (O6)}  $\langle a, b, c\di a^{9}=b^{9}=c^3=1,
[a,b]=c,[c,a]=b^3,[c,b]=a^3,[a^3,b]=1 \rangle$;%(E6)

{\rm (O7)}  $\langle a, b, c \di a^{9}=b^{9}=c^3=1,
[a,b]=c,[c,a]=b^{-3},[c,b]=a^3,[a^3,b]=1 \rangle$;%(E7)

{\rm (P1)} $\langle a, b, c \di a^{p^{n+1}}=b^{p^{n+1}}=c^{p}=1,
[a,b]=c,[c,a]=a^{p^{n}},[c,b]=b^{p^n}\rangle$, where \mbox{\hskip0.65in}$p>2$ and $n\ge 2$
for $p=3$;%(F1)

{\rm (P2)} $\langle a, b, c \di a^{p^{n+1}}=b^{p^{n+1}}=c^{p}=1,
[a,b]=c,[c,a]=a^{p^n}b^{\nu p^{n}},[c,b]=b^{p^n}\rangle$,
where \mbox{\hskip0.65in}$p>2$, $n\ge 2$ for $p=3$, $\nu=1$ or a
fixed quadratic non-residue modular $p$;%(F?)

{\rm (P3)} $\langle a, b, c \di a^{p^{n+1}}=b^{p^{n+1}}=c^{p}=1,
[a,b]=c,[c,a]=b^{\nu p^{n}},[c,b]=a^{-p^n},[a^{p^n},b]=\mbox{\hskip0.65in}1\rangle$, where
$p>2$, $n\ge 2$ for $p=3$, $\nu=1$ or a
fixed quadratic non-residue \mbox{\hskip0.65in}modular $p$;%(F?)

{\rm (P4)} $\langle a, b, c \di a^{p^{n+1}}=b^{p^{n+1}}=c^{p}=1,
[a,b]=c,[c,a]^{1+r}=a^{p^n}b^{p^{n}},[c,b]^{1+r}=$ \mbox{\hskip0.65in}$a^{-r p^n}b^{
p^n},[a^{p^n},b]=1\rangle$, where $p>2$, $n\ge 2$ for $p=3$, $r=1,2,\dots,p-2$;%(F?)

{\rm (P5)} $\langle a, b, c \di a^{2^{n+1}}=b^{2^{n+1}}=c^{2}=1,
[a,b]=c,[c,a]=b^{2^{n}},[c,b]=a^{2^n},[a^{2^n},b]=1\rangle$, \mbox{\hskip0.65in}where
$n\ge 3$;%(F?)

{\rm (P6)} $\langle a, b, c \di a^{2^{n+1}}=b^{2^{n+1}}=c^{2}=1,
[a,b]=c,[c,a]=a^{2^{n}},[c,b]=b^{2^n}\rangle$, where \mbox{\hskip0.65in}$n\ge 3$;%(F1)

{\rm (P7)} $\langle a, b, c \di a^{2^{n+1}}=b^{2^{n+1}}=c^{2}=1,
[a,b]=c,[c,a]=a^{2^n}b^{2^{n}},[c,b]=a^{2^n},[a^{2^n},b]=\mbox{\hskip0.65in}1\rangle$,
where $n\ge 3$;%(F?)

{\rm (P8)} $\langle a, b, c,d \di a^{p^{n+1}}=b^{p^{n}}=c^{p}=d^p=1,
[a,b]=c,[c,a]=a^{p^{n}},[c,b]=d,[d,a]=\mbox{\hskip0.65in}[d,b]=1\rangle$, where $n\ge 3$ for $p=2$ and $n\ge 2$
for $p=3$;%(G4)

{\rm (P9)} $\langle a, b, c,d \di a^{p^{n}}=b^{p^{n+1}}=c^{p}=d^p=1,
[a,b]=c,[c,a]=b^{\nu p^n},[c,b]=d,[d,a]=\mbox{\hskip0.65in}[d,b]=1\rangle$,
where $n\ge 3$ for $p=2$ and $n\ge 2$
for $p=3$, $\nu=1$ or a
fixed \mbox{\hskip0.65in}quadratic non-residue modular $p$;%(G3)

{\rm (P10)} $\langle a, b, c,d,e \di a^{p^{n}}=b^{p^{n}}=c^{p}=d^p=e^p=1,
[a,b]=c,[c,a]=d,[c,b]=\mbox{\hskip0.65in}e,[d,a]=[d,b]=[e,a]=[e,b]=1\rangle$, where
$n\ge 3$ for $p=2$ and $n\ge 2$
for \mbox{\hskip0.65in}$p=3$;%(H)

{\rm (Q1)} $\langle a, b, c \di a^{p^{n+1}}=b^{p^{n+1}}=1,
[a,b]=c,[c,a]=c^p=b^{sp^{n}},[c,b]=a^{-\nu p^n}c^{\nu tp^n}\rangle$, \mbox{\hskip0.65in}where $n\ge 3$ for $p=2$ and $n\ge 2$, $\nu=1$ or a
fixed quadratic non-residue \mbox{\hskip0.65in}modular $p$, $s\in F_p^*$, $t=0,1,\dots,\frac{p-1}{2}$;%(??)

{\rm (Q2)} $\langle a, b, c \di a^{p^{n+1}}=b^{p^{n}}=c^{p^2}=1,
[a,b]=c,[c,a]=c^{p},[c,b]=a^{-\nu p^n}c^{t\nu p}\rangle$, where \mbox{\hskip0.65in}$n\ge 3$ for $p=2$ and $n\ge 2$, $\nu=1$ or a
fixed quadratic non-residue modular \mbox{\hskip0.65in}$p$, $t=0,1,\dots,\frac{p-1}{2}$;%(??)

{\rm (Q3)} $\langle a, b, c \di a^{p^{n+1}}=b^{p^{n+1}}=1,
[a,b]=c,[c,a]=c^{p}=a^{p^n},[c,b]=a^{s p^n}b^{p^n}\rangle$, where \mbox{\hskip0.65in}$n\ge 3$ for $p=2$ and $n\ge 2$, $s\in F_p$;%(??)

{\rm (Q4)} $\langle a, b, c \di a^{p^{n+1}}=b^{p^{n+1}}=1,
[a,b]=c,[c,a]=c^{p}=a^{p^n},[c,b]=b^{sp^n}\rangle$, where \mbox{\hskip0.65in}$n\ge 3$ for $p=2$ and $n\ge 2$, $s=2,3,\dots,p-1$;%(??)

{\rm (Q5)} $\langle a, b, c,d \di a^{p^{n+1}}=b^{p^{n}}=d^p=1,
[a,b]=c,[c,a]=c^p=a^{p^n},[c,b]=d,[d,a]=\mbox{\hskip0.6in}[d,b]=1\rangle$, where $n\ge 3$ for $p=2$ and $n\ge 2$;%(??)

{\rm (Q6)} $\langle a, b, c \di a^{p^{n}}=b^{p^{n+1}}=c^{p^2}=1,
[a,b]=c,[c,a]=c^p,[c,b]=b^{p^n}\rangle$, where $n\ge 3$ \mbox{\hskip0.65in}for $p=2$ and $n\ge 2$;%(??)

{\rm (Q7)} $\langle a, b, c,d \di a^{p^{n}}=b^{p^{n+1}}=d^p=1,
[a,b]=c,[c,a]=c^p=b^{sp^n},[c,b]=d,[d,a]=\mbox{\hskip0.6in}[d,b]=1\rangle$, where $n\ge 3$ for $p=2$ and $n\ge 2$, $s\in F_p^*$;%(??)

{\rm (Q8)} $\langle a, b, c,d \di a^{p^{n}}=b^{p^{n}}=c^{p^2}=d^p=1,
[a,b]=c,[c,a]=c^p,[c,b]=d,[d,a]=\mbox{\hskip0.6in}[d,b]=1\rangle$, where $n\ge 3$ for $p=2$ and $n\ge 2$;%(??)

{\rm (R1)} $\langle a, b, c,d \di a^4=b^2=c^4=d^2=1,
[a,b]=c,[c,a]=d,[c,b]=c^2,[d,a]=[d,b]=\mbox{\hskip0.6in}1\rangle$; %(I1)

{\rm (R2)} $\langle a, b, c \di a^4=b^4=c^4=1,
[a,b]=c,[c,a]=b^2,[c,b]=c^2\rangle$;%(I2)

{\rm (R3)} $\langle a, b, c \di a^8=b^2=c^4=1,
[a,b]=c,[c,a]=a^4,[c,b]=c^2\rangle$;%(I3)

{\rm (R4)} $\langle a, b, c \di a^8=c^4=1,b^2=a^4,
[a,b]=c,[c,a]=a^4,[c,b]=c^2\rangle$;%(I4)

{\rm (R5)} $\langle a, b, c \di a^8=b^2=c^4=1,
[a,b]=c,[c,a]=a^4c^2,[c,b]=c^2,[c^2,a]=1\rangle$;%(I5)

{\rm (R6)} $\langle a, b, c \di a^8=b^4=1,c^2=a^4b^2,
[a,b]=c,[c,a]=b^2,[c,b]=c^2\rangle$;%(I6)

{\rm (R7)} $\langle a, b, c\di a^8=b^4=1,c^2=b^2,
[a,b]=c,[c,a]=a^4b^2,[c,b]=c^2\rangle$;%(I7)

{\rm (S1)} $\langle a, b, c \di a^{p^{n+1}}=b^{p^{m+1}}=c^p=1,
[a,b]=c,[c,a]=a^{p^n},[c,b]=b^{sp^m}\rangle$, where \mbox{\hskip0.6in}$n\ge 3$ and $m\ge 2$ for $p=2$ and $n>m$, $s\in F_p^*$;%(??)

{\rm (S2)} $\langle a, b, c \di a^{p^{n+1}}=b^{p^{m+1}}=c^p=1,
[a,b]=c,[c,a]=b^{\nu_1p^m},[c,b]=a^{-\nu_2 p^n},[a,b^{p^m}]=\mbox{\hskip0.6in}1\rangle$, where $n\ge 3$ and $m\ge 2$ for $p=2$ and $n>m$, $\nu_1,\nu_2=1$ or a fixed \mbox{\hskip0.6in}quadratic non-residue modular $p$;%(??)

{\rm (S3)} $\langle a, b, c,d \di a^{p^{n+1}}=b^{p^{m}}=c^p=d^p=1,
[a,b]=c,[c,a]=d,[c,b]=a^{-\nu p^n},[d,a]=$ \mbox{\hskip0.6in}$[d,b]=1\rangle$, where $n\ge 3$ and $m\ge 2$ for $p=2$ and $n>m$, $\nu=1$ or a fixed \mbox{\hskip0.6in}quadratic non-residue modular $p$;%(??)

{\rm (S4)} $\langle a, b, c,d \di a^{p^{n}}=b^{p^{m+1}}=c^p=d^p=1,
[a,b]=c,[c,a]=b^{\nu p^m},[c,b]=d,[d,a]=\mbox{\hskip0.6in}[d,b]=1\rangle$, where $n\ge 3$ and $m\ge 2$ for $p=2$ and $n>m$, $\nu=1$ or a fixed \mbox{\hskip0.6in}quadratic non-residue modular $p$;%(??)

{\rm (S5)} $\langle a, b, c,d \di a^{p^{n}}=b^{p^{m+1}}=c^p=d^p=1,
[a,b]=c,[c,a]=d,[c,b]=b^{p^m},[d,a]=$ \mbox{\hskip0.6in}$[d,b]=1\rangle$, where $n\ge 3$ and $m\ge 2$ for $p=2$ and $n>m$;%(??)

{\rm (S6)} $\langle a, b, c,d \di a^{p^{n+1}}=b^{p^{m}}=c^p=d^p=1,
[a,b]=c,[c,a]=a^{p^n},[c,b]=d,[d,a]=$ \mbox{\hskip0.6in}$[d,b]=1\rangle$, where $n\ge 3$ and $m\ge 2$ for $p=2$ and $n>m$;%(??)

{\rm (S7)} $\langle a, b, c,d,e \di a^{p^{n}}=b^{p^{m}}=c^p=d^p=e^p=1,
[a,b]=c,[c,a]=d,[c,b]=e,$ \mbox{\hskip0.6in}$[d,a]=[d,b]=[e,a]=[e,b]=1\rangle$, where $n\ge 3$ and $m\ge 2$ for $p=2$ and \mbox{\hskip0.6in}$n>m$;%(??)

{\rm (T1)} $\langle a, b, c\di a^{2^{n+1}}=b^4=1,c^2=b^2,
[a,b]=c,[c,a]=a^{2^n},[c,b]=c^2\rangle$, where $n\ge 3$;%{\rm (J7)}

{\rm (T2)} $\langle a, b, c \di a^{2^{n+1}}=b^4=1, c^2=a^{2^n},
[a,b]=c,[c,a]=b^2,[c,b]=c^2\rangle$, where $n\ge 3$;%{\rm (J5)}

{\rm (T3)} $\langle a, b, c,d \di a^{2^{n+1}}=b^2=d^2=1,c^2=a^{2^n},
[a,b]=c,[c,a]=d,[c,b]=c^2,[d,a]=\mbox{\hskip0.65in}[d,b]=1\rangle$,
where $n\ge 3$;%{\rm (J4)}

{\rm (T4)} $\langle a, b, c \di a^{2^n}=b^4=c^4=1,
[a,b]=c,[c,a]=b^2,[c,b]=c^2\rangle$, where $n\ge 3$;%{\rm (J3)}

{\rm (T5)} $\langle a, b, c,d \di a^{2^n}=b^4=d^2=1, c^2=b^2,
[a,b]=c,[c,a]=d,[c,b]=c^2,[d,a]=\mbox{\hskip0.65in}[d,b]=1\rangle$,
where $n\ge 3$;%{\rm (J2)}

{\rm (T6)} $\langle a, b, c \di a^{2^{n+1}}=b^2=c^4=1,
[a,b]=c,[c,a]=a^{2^n},[c,b]=c^2\rangle$, where $n\ge 3$;%{\rm (J6)}

{\rm (T7)} $\langle a, b, c,d \di a^{2^n}=b^2=c^4=d^2=1,
[a,b]=c,[c,a]=d,[c,b]=c^2,[d,a]=\mbox{\hskip0.65in}[d,b]=1\rangle$,
where $n\ge 3$;%{\rm (J1)}

{\rm (U1)} $\langle a, b, c \di a^{p^{n+1}}=b^{p^{m+1}}=1,
[a,b]=c,[c,a]=a^{p^n},[c,b]=c^{-p}=b^{sp^m}\rangle$, where \mbox{\hskip0.65in}$n>m\ge 2$, $s\in F_p^*$;%(??)

{\rm (U2)} $\langle a, b, c \di a^{p^{n+1}}=b^{p^{m+1}}=1,
[a,b]=c,[c,a]=b^{\nu p^m},[c,b]=c^{-p}=a^{sp^n}\rangle$, where \mbox{\hskip0.65in}$n>m\ge 2$, $s\in F_p^*$, $\nu=1$ or a fixed quadratic non-residue modular $p$;%(??)

{\rm (U3)} $\langle a, b, c,d \di a^{p^{n+1}}=b^{p^{m}}=d^p=1,
[a,b]=c,[c,b]=c^{-p}=a^{s p^n},[c,a]=\mbox{\hskip0.65in}d,[d,a]=[d,b]=1\rangle$, where $n>m\ge 2$, $s\in F_p^*$, $\nu=1$ or a fixed quadratic \mbox{\hskip0.65in}non-residue modular $p$;%(??)

{\rm (U4)} $\langle a, b, c \di a^{p^{n}}=b^{p^{m+1}}=c^{p^2}=1,
[a,b]=c,[c,a]=b^{\nu p^m},[c,b]=c^{-p}\rangle$, where \mbox{\hskip0.65in}$n>m\ge 2$, $\nu=1$ or a fixed quadratic non-residue modular $p$;%(??)

{\rm (U5)} $\langle a, b, c,d \di a^{p^{n}}=b^{p^{m+1}}=d^p=1,
[a,b]=c,[c,a]=d,[c,b]=c^{-p}=b^{p^m},[d,a]=$ \mbox{\hskip0.65in}$[d,b]=1\rangle$, where $n>m\ge 2$;%(??)

{\rm (U6)} $\langle a, b, c \di a^{p^{n+1}}=b^{p^{m}}=c^{p^2}=1,
[a,b]=c,[c,a]=a^{p^n},[c,b]=c^{-p}\rangle$, where \mbox{\hskip0.65in}$n>m\ge 2$;%(??)

{\rm (U7)} $\langle a, b, c,d \di a^{p^{n}}=b^{p^{m}}=c^{p^2}=d^p=1,
[a,b]=c,[c,a]=d,[c,b]=c^{-p},[d,a]=\mbox{\hskip0.65in}[d,b]=1\rangle$, where $n>m\ge 2$;%(??)

{\rm (V1)} $\langle a, b, c \di a^{p^{n+1}}=b^{p^{m+1}}=1,
[a,b]=c,[c,a]=c^p=b^{sp^{m}},[c,b]=a^{-\nu p^n}c^{\nu tp}\rangle$, \mbox{\hskip0.65in}where $n>m\ge 2$, $\nu=1$ or a
fixed quadratic non-residue modular $p$, $s\in F_p^*$, \mbox{\hskip0.65in}$t=0,1,\dots,\frac{p-1}{2}$;%(??)

{\rm (V2)} $\langle a, b, c \di a^{p^{n+1}}=b^{p^{m}}=c^{p^2}=1,
[a,b]=c,[c,a]=c^{p},[c,b]=a^{-\nu p^n}c^{t\nu p}\rangle$, \mbox{\hskip0.65in}where $n>m\ge 2$, $\nu=1$ or a
fixed quadratic non-residue modular $p$, $t=\mbox{\hskip0.65in}0,1,\dots,\frac{p-1}{2}$;%(??)

{\rm (V3)} $\langle a, b, c \di a^{p^{n+1}}=b^{p^{m+1}}=1,
[a,b]=c,[c,a]=c^{p}=a^{p^n},[c,b]=b^{sp^m}\rangle$, where \mbox{\hskip0.65in}$n>m\ge 2$, $s=F_p^*$;%(??)

{\rm (V4)} $\langle a, b, c,d \di a^{p^{n+1}}=b^{p^{m}}=d^p=1,
[a,b]=c,[c,a]=c^p=a^{p^n},[c,b]=d,[d,a]=\mbox{\hskip0.6in}[d,b]=1\rangle$, where $n>m\ge 2$;%(??)

{\rm (V5)} $\langle a, b, c \di a^{p^{n}}=b^{p^{m+1}}=c^{p^2}=1,
[a,b]=c,[c,a]=c^p,[c,b]=b^{p^m}c^{sp},[c^p,b]=1\rangle$, \mbox{\hskip0.6in}where $n>m\ge 2$, $s\in F_p$;%(??)

{\rm (V6)} $\langle a, b, c,d \di a^{p^{n}}=b^{p^{m+1}}=d^p=1,
[a,b]=c,[c,a]=c^p=b^{sp^m},[c,b]=d,[d,a]=\mbox{\hskip0.6in}[d,b]=1\rangle$, where $n>m\ge 2$, $s\in F_p^*$;%(??)

{\rm (V7)} $\langle a, b, c,d \di a^{p^{n}}=b^{p^{m}}=c^{p^2}=d^p=1,
[a,b]=c,[c,a]=c^p,[c,b]=d,[d,a]=\mbox{\hskip0.6in}[d,b]=1\rangle$, where $n>m\ge 2$.
\end{thm}
\demo {\bf Case 1.} $n=m$.

In this case, $m\ge 2$ for $p=2$. If $p=n=m=2$, then
$|G|=2^7$. By checking the list of groups of order
$2^7$, we get the groups of Type (N1)--(N16). If $p=3$ and $n=m=1$, then $|G|=3^5$. By checking the list of groups of order
$3^5$, we get the groups of Type (O1)--(O7). In the following, we may assume that $n>1$ for $p=3$ and $n>2$ for $p=2$.

\medskip

{\bf Subcase 1.1} $\Phi(G')=1$.

In this subcase, $v(G)=(0,0)^t$. Suppose that $G$ and $\bar{G}$ are two groups described in the theorem. By Theorem \ref{isomorphic-5.1} and \ref{isomorphic-5.3},
$G\cong \bar{G}$ if and only if there exists $Y={\left(
 \begin{array}{cc}
 y_{11}& y_{12}\\
 y_{21}& y_{22}
 \end{array}
 \right)}$, an invertible matrix over $F_p$, such that $
w(\bar{G})=
Yw(G)Y^t$. That is, $w(\bar{G})
$ and $w(G)$ are mutually congruent. By Lemma \ref{hetong-p}, \ref{hetong-2} and \ref{hetong-p-non-invertible}, we get the groups of Type (P1)--(P10).

 \medskip

{\bf Subcase 1.2.} $\Phi(G')\neq 1$.

If $n=m=1$ and $p>3$, then $[a,b^p]=[a,b]^p=1$ and hence $\Phi(G')=1$. Thus, in this subcase, we have $n\ge 2$. Since $v(G)=\left(
 \begin{array}{c}
 w_{31}\\
 w_{32}
 \end{array}
 \right)\neq \left(
 \begin{array}{c}
 0\\
 0
 \end{array}
 \right)$, there exist $w_{41},w_{42}$ such that $\left(
 \begin{array}{cc}
 w_{41}& w_{31}\\
 w_{42}& w_{32}
 \end{array}
 \right)$ is invertible. Let $Y=\left(
 \begin{array}{cc}
 w_{41}& w_{31}\\
 w_{42}& w_{32}
 \end{array}
 \right)^{-1}$. Then $Yv(G)=\left(
 \begin{array}{c}
 0\\
 1
 \end{array}
 \right)$.
 In the following, we assume that $G$ and $\bar{G}$ be two groups described in the theorem with $v(G)=v(\bar{G})=\left(
 \begin{array}{c}
 0\\
 1
 \end{array}
 \right)$.
By Theorem \ref{isomorphic-5.1} and \ref{isomorphic-5.3},
$G\cong \bar{G}$ if and only if there exists $Y={\left(
 \begin{array}{cc}
 y_{11}& 0\\
 y_{21}& 1
 \end{array}
 \right)}$, an invertible matrix over $F_p$, such that $
w(\bar{G})=
Yw(G)Y^t.$

By suitably choosing $y_{21}$,  we can simplify $w(G)$ to be one of the following types:
(a)  $\left(
 \begin{array}{cc}
 w_{11}& w_{12}\\
 0& w_{22}
 \end{array}
 \right)$ where $w_{11}\neq 0$,
(b)  $\left(
 \begin{array}{cc}
 0  & w_{12}\\
 -w_{12}& w_{22}
 \end{array}
 \right)$ where $w_{12}\neq 0$,
(c)  $\left(
 \begin{array}{cc}
 0&w_{12}\\
 w_{21}&0
 \end{array}
 \right)$ where $w_{12}\neq 0$ and $w_{21}\neq -w_{12}$,
(d)  $\left(
 \begin{array}{cc}
 0  & 0\\
 w_{21}& 0
 \end{array}
 \right)$ where $w_{21}\neq 0$ and
(e)  $\left(
 \begin{array}{cc}
 0&0\\
 0&w_{22}
 \end{array}
 \right)$.
In the following, we assume that both $w(G)$ and
$w(\bar{G})$ are such matrices. It is easy to check that (i)
 different types give non-isomorphic groups, (ii) $G\cong \bar{G}$ if and only if there exists $y_{11}\in F_p^*$ such that $ w(\bar{G})=Yw(G)Y^t$ where $Y=\diag(y_{11},1)$.
\begin{table}[h]
  \centering
%\caption{Properties of the groups of Type (A) in Theorem
%\ref{main}}\label{table: 8.1}
{\small
%\tiny
%\scriptsize
%\footnotesize
\begin{tabular}{c|c|c|c|c|c}
$w(G)$ &$y_{11}$ &Remak 1& $w(\bar{G})$&Group&Remark 2 \\
   \hline
(a)&$z^{-1}$ & $w_{11}=\nu z^2$ & ${\left(
 \begin{array}{cc}
 \nu & w_{12}z^{-1}\\
 0& w_{22}
 \end{array}
 \right)}$& $\begin{array}{c} {\rm (Q1)\ if\ }w_{22}\neq 0 \\  {\rm (Q2)\ if\ }w_{22}=0\end{array}$ & $\begin{array}{c} s=(w_{22})^{-1}\\ t=w_{12}z^{-1}\\ \end{array}$\\
\hline
(b)&$w_{12}^{-1}$ & &${\left(
 \begin{array}{cc}
 0 & 1\\
 -1& w_{22}
 \end{array}
 \right)}$&(Q3) &\\
\hline
(c)&$w_{12}^{-1}$ & &${\left(
 \begin{array}{cc}
 0& 1\\
 w_{21}w_{12}^{-1}& 0
 \end{array}
 \right)}$ &$\begin{array}{c} {\rm (Q4)\ if\ }w_{21}\neq 0\\ {\rm (Q5)\ if\ }w_{21}=0\end{array}$ & $\begin{array}{c}s=-(w_{21})^{-1}w_{12}\\ \\ \end{array}$\\
\hline
(d)& $-w_{21}^{-1}$&  &${\left(
 \begin{array}{cc}
 0& 0\\
 -1& 0
 \end{array}
 \right)}$ & (Q6) & \\
\hline
(e)&  &  &${\left(
 \begin{array}{cc}
 0& 0\\
 0& w_{22}
 \end{array}
 \right)}$ & $\begin{array}{c} {\rm (Q7)\ if\ }w_{22}\neq 0\\ {\rm (Q8)\ if\ }w_{22}=0\end{array}$ & $\begin{array}{c}s=w_{22}^{-1}\\ \\ \end{array}$\\
   \hline
\end{tabular}}
\caption{Subcase 1.2 in Theorem \ref{5}}
\label{table 6}
\end{table}
By Table \ref{table 6}, we get the groups of Type (Q1)--(Q8).

{\bf Case 2.} $n>m$.

If $p=n=2$ and $m=1$, then
$|G|=2^6$. By checking the list of groups of order
$2^6$, we get the groups of Type (R1)--(R7). In the following, we may assume that $n\ge 3$ for $p=2$.

\medskip

{\bf Subcase 2.1.} $\Phi(G')=1$.

By Theorem \ref{isomorphic-5.3}, $m\ge 2$ for $p=2$. Suppose that $G$ and $\bar{G}$ are two groups described in the theorem. By Theorem \ref{isomorphic-5.1}, \ref{isomorphic-5.2} and \ref{isomorphic-5.3},
$G\cong \bar{G}$ if and only if there exists $Y={\left(
 \begin{array}{cc}
 y_{11}& y_{12}\\
 y_{21}p^{n-m}& y_{22}
 \end{array}
 \right)}$, an invertible matrix over $F_p$, such that $
w(\bar{G})=
Y_1w(G)Y^t$, where
$Y_1={\left(
 \begin{array}{cc}
 y_{11}& 0\\
 y_{21}& y_{22}
 \end{array}
 \right)}$.

By suitably choosing $y_{21}$ and $y_{12}$, that is, using an elementary operation and an elementary column operation, we can simplify $w(G)$ to be such a matrix, in which every column and every row have at most one non-zero entry.
In the following, we assume that both $w(G)=(w_{ij})$ and $w(\bar{G})=(\bar{w}_{ij})$ are such matrices.
It is easy to check that (i) for all possible subscripts $i,j$, $\bar{w}_{ij}\neq 0$ if and only if ${w}_{ij}\neq 0$; (ii) $G\cong \bar{G}$ if and only if there exists $Y=\diag(y_{11},y_{22})$, an invertible matrix over $F_p$, such that $
w(\bar{G})=Yw(G)Y$.

If $w(G)={\left(
 \begin{array}{cc}
 0&w_{12}\\
 w_{21}&0
 \end{array}
 \right)}$ where $w_{12}w_{21}\neq 0$, then, letting $Y=\diag(w_{12}^{-1},1)$, we have $w(\bar{G})=Yw(G)Y={\left(
 \begin{array}{cc}
 0& 1\\
 w_{21}w_{12}^{-1}& 0
 \end{array}
 \right)}$. Hence we get the groups of Type (S1) where $s=-w_{21}^{-1}w_{12}$.
It is easy to see that different $s$ give non-isomorphic groups.

If $w(G)={\left(
 \begin{array}{cc}
 w_{11}&0\\
 0& w_{22}
 \end{array}
 \right)}$ where $w_{11}w_{22}\neq 0$, then, letting $Y=\diag(y_{11},y_{22})$, we have $w(\bar{G})=Yw(G)Y={\left(
 \begin{array}{cc}
 w_{11}y_{11}^2& 0\\
 0&w_{22}y_{22}^{2}
 \end{array}
 \right)}$. Hence we can simplify $w(G)$ to be $\left(
 \begin{array}{cc}
 \nu_1& 0\\
 0& \nu_2
 \end{array}
 \right)$ where $\nu_1, \nu_2=1$ or a fixed quadratic non-residue modular $p$.
Thus we get the groups of Type (S2).
It is easy to see that different $\nu_1$ or $\nu_2$ give non-isomorphic groups.

If $w(G)$ is invertible, then $w(G)$ is one of the above types. If $w(G)$ is of rank 1, then $w(G)$ is one of the following types:
\begin{center}
{\rm (a)} $\left(
 \begin{array}{cc}
 w_{11}& 0\\
 0& 0
 \end{array}
 \right)$,
{\rm (b)} $\left(
 \begin{array}{cc}
 0& 0\\
 0& w_{22}
 \end{array}
 \right)$,
  {\rm (c)} $\left(
 \begin{array}{cc}
 0& 0\\
 w_{21}& 0
 \end{array}
 \right)$,
{\rm (d)} $\left(
 \begin{array}{cc}
 0& w_{12}\\
 0& 0
 \end{array}
 \right)$.
 \end{center}
By similar arguments as above, we get the groups of Type (S3)--(S6) respectively. If $w(G)=0$, then $G$ is the group of Type
(S7).

\medskip

 {\bf Subcase 2.2.} $\Phi(G')\neq 1$.

\medskip

 {\bf Subcase 2.2.1.} $m=1$.

If $p>2$, then $c^p=[a,b]^p=[a,b^p]=1$ and hence $\Phi(G')=1$, a contradiction.
Thus $p=2$. By Theorem \ref{isomorphic-5.3}, $v(G)=(1,0)^t$ and
$G\cong \bar{G}$ if and only if there exists $Y={\left(
 \begin{array}{cc}
 1& y_{12}\\
 y_{21}2^{n-1}& 1
 \end{array}
 \right)}$, an invertible matrix over $F_p$, such that $
w(\bar{G})=
Y_1w(G)Y^t$, where
$Y_1={\left(
 \begin{array}{cc}
 1& 0\\
 y_{21}& 1
 \end{array}
 \right)}$.
By suitably choosing $y_{21}$ and $y_{12}$, that is, using an elementary operation and an elementary column operation, we can simplify $w(G)$ to be such a matrix, in which every column and every row have at most one non-zero entry.
Then we may assume that $w(G)$ is one of the following matrices:

{\rm (a)} $\left(
 \begin{array}{cc}
 0& 1\\
 1& 0
 \end{array}
 \right)$,\ \
{\rm (b)} $\left(
 \begin{array}{cc}
 1& 0\\
 0& 1
 \end{array}
 \right)$,\ \
{\rm (c)} $\left(
 \begin{array}{cc}
 1& 0\\
 0& 0
 \end{array}
 \right)$,\ \
{\rm (d)} $\left(
 \begin{array}{cc}
 0& 0\\
 0& 1
 \end{array}
 \right)$,\ \
 {\rm (e)} $\left(
 \begin{array}{cc}
 0& 0\\
 1& 0
 \end{array}
 \right)$,

{\rm (f)} $\left(
 \begin{array}{cc}
 0& 1\\
 0& 0
 \end{array}
 \right)$,\ \
{\rm (g)} $\left(
 \begin{array}{cc}
 0& 0\\
 0& 0
 \end{array}
 \right)$.

\noindent It is easy to check that different types give non-isomorphic groups. Hence we get the groups of Type (T1)--(T7).

\medskip

{\bf Subcase 2.2.2.} $m\ge 2$.

If $w_{32}\neq 0$, then, letting $Y=\left(
 \begin{array}{cc}
 w_{32}& -w_{31}\\
 0& w_{32}^{-1}
 \end{array}
 \right)$, we have $Yv(G)=\left(
 \begin{array}{c}
 0\\
 1
 \end{array}
 \right)$.
If $w_{32}=0$, then $w_{31}\neq 0$. Let $Y=\left(
 \begin{array}{cc}
 w_{31}^{-1}& 0\\
0& 1
 \end{array}
 \right)$. Then $Yv(G)=\left(
 \begin{array}{c}
 1\\
 0
 \end{array}
 \right)$. By Theorem \ref{isomorphic-5.1} and \ref{isomorphic-5.3}, groups with $v(G)=\left(
 \begin{array}{c}
 1\\
 0
 \end{array}
 \right)$ and $\left(
 \begin{array}{c}
 0\\
 1
 \end{array}
 \right)$ respectively are mutually non-isomorphic.

\medskip

{\bf Subcase 2.2.2.1.} $v(G)=\left(
 \begin{array}{c}
 1\\
 0
 \end{array}
 \right)$.

By calculation, $\left(
 \begin{array}{cc}
 y_{11}&  y_{12}\\
 0& y_{22}
 \end{array}
 \right)\left(
 \begin{array}{c}
 1\\
0
 \end{array}\right)
=\left(
 \begin{array}{c}
 1\\
0
 \end{array}\right)$
 if and only if $y_{11}=1$.
Suppose that $G$ and $\bar{G}$ are two groups described in the theorem with $v(G)=\left(
 \begin{array}{c}
 1\\
 0
 \end{array}
 \right)$. By Theorem \ref{isomorphic-5.1} and \ref{isomorphic-5.3},
$G\cong \bar{G}$ if and only if there exists $Y={\left(
 \begin{array}{cc}
 1& y_{12}\\
 0& y_{22}
 \end{array}
 \right)}$ and $Y_1={\left(
 \begin{array}{cc}
 1& 0\\
 y_{21}& y_{22}
 \end{array}
 \right)}$ where $y_{22}\neq 0$ such that $
w(\bar{G})=
Y_1w(G)Y^t$.

By suitably choosing $y_{21}$ and $y_{12}$, that is, using an elementary operation and an column operation, we can simplify $w(G)$ to be such a matrix, in which every column and every row have at most one non-zero entry.
In the following, we assume that both $w(G)=(w_{ij})$ and $w(\bar{G})=(\bar{w}_{ij})$ are such matrices.
It is easy to check that (i) for all possible subscripts $i,j$, $\bar{w}_{ij}\neq 0$ if and only if ${w}_{ij}\neq 0$; (ii) $G\cong \bar{G}$ if and only if there exists $Y=\diag(1,y_{22})$, an invertible matrix over $F_p$, such that $
w(\bar{G})=Yw(G)Y$.

If $w(G)={\left(
 \begin{array}{cc}
 0&w_{12}\\
 w_{21}&0
 \end{array}
 \right)}$ where $w_{12}w_{21}\neq 0$, then, letting $Y=\diag(1,w_{12}^{-1})$, we have $w(\bar{G})=Yw(G)Y={\left(
 \begin{array}{cc}
 0& 1\\
 w_{21}w_{12}^{-1}& 0
 \end{array}
 \right)}$. Hence we get the groups of Type (U1) where $s=-w_{21}^{-1}w_{12}$.
It is easy to see that different $s$ give non-isomorphic groups.

If $w(G)={\left(
 \begin{array}{cc}
 w_{11}&0\\
 0& w_{22}
 \end{array}
 \right)}$ where $w_{11}w_{22}\neq 0$, then, letting $Y=\diag(1,y_{22})$, we have $w(\bar{G})=Yw(G)Y={\left(
 \begin{array}{cc}
 w_{11}& 0\\
 0&w_{22}y_{22}^{2}
 \end{array}
 \right)}$. Hence we can simplify $w(G)$ to be $\left(
 \begin{array}{cc}
 w_{11}& 0\\
 0& \nu
 \end{array}
 \right)$ where $\nu=1$ or a fixed quadratic non-residue modular $p$.
Thus we get the groups of Type (U2) where $s=-w_{11}^{-1}$.
It is easy to see that different $s$ or $\nu$ give non-isomorphic groups.

If $w(G)$ is invertible, then $w(G)$ is one of the above types. If $w(G)$ is of rank 1, then $w(G)$ is one of the following types:
\begin{center}
{\rm (a)} $\left(
 \begin{array}{cc}
 w_{11}& 0\\
 0& 0
 \end{array}
 \right)$,
{\rm (b)} $\left(
 \begin{array}{cc}
 0& 0\\
 0& w_{22}
 \end{array}
 \right)$,
  {\rm (c)} $\left(
 \begin{array}{cc}
 0& 0\\
 w_{21}& 0
 \end{array}
 \right)$,
{\rm (d)} $\left(
 \begin{array}{cc}
 0& w_{12}\\
 0& 0
 \end{array}
 \right)$.
 \end{center}
By similar arguments as above, we get the groups of Type (U3)--(U6) respectively. If $w(G)=0$, then $G$ is the group of Type
(U7).

\medskip

{\bf Subcase 2.2.2.2.} $v(G)=\left(
 \begin{array}{c}
 0\\
 1
 \end{array}
 \right)$.

By calculation, $\left(
 \begin{array}{cc}
 y_{11}&  y_{12}\\
 0& y_{22}
 \end{array}
 \right)\left(
 \begin{array}{c}
 0\\
1
 \end{array}\right)
=\left(
 \begin{array}{c}
 0\\
1
 \end{array}\right)$
 if and only if $y_{12}=0$ and $y_{22}=1$.
Suppose that $G$ and $\bar{G}$ are two groups described in the theorem with $v(G)=\left(
 \begin{array}{c}
 0\\
1
 \end{array}\right)$. By Theorem \ref{isomorphic-5.1} and \ref{isomorphic-5.3},
$G\cong \bar{G}$ if and only if there exist $Y={\left(
 \begin{array}{cc}
 y_{11}& 0\\
 0& 1
 \end{array}
 \right)}$ and $Y_1={\left(
 \begin{array}{cc}
 y_{11}& 0\\
 y_{21}& 1
 \end{array}
 \right)}$ where $y_{11}\neq 0$ such that $
w(\bar{G})=
Y_1w(G)Y$.

By suitably choosing $y_{21}$, that is, using an elementary operation, we can simplify $w(G)$ to be one of the following types: (a) $\left(
 \begin{array}{cc}
 w_{11}&  w_{12}\\
 0 & w_{22}
 \end{array}
 \right)$ where $w_{11}\neq 0$, (b) $\left(
 \begin{array}{cc}
 0&  w_{12}\\
 w_{21} & 0
 \end{array}
 \right)$ where $w_{12}\neq 0$, (c) $\left(
 \begin{array}{cc}
 0&  0\\
 w_{21} & w_{22}
 \end{array}
 \right)$ where $w_{21}\neq 0$, (d) $\left(
 \begin{array}{cc}
 0&  0\\
 0 & w_{22}
 \end{array}
 \right)$.
In the following, we assume that both $w(G)=(w_{ij})$ and $w(\bar{G})=(\bar{w}_{ij})$ are such matrices.
It is easy to check that (i) different types give non-isomorphic groups; (ii) $G\cong \bar{G}$ if and only if there exists $Y=\diag(y_{11},1)$, an invertible matrix over $F_p$, such that $
w(\bar{G})=Yw(G)Y$.
\begin{table}[h]
  \centering
%\caption{Properties of the groups of Type (A) in Theorem
%\ref{main}}\label{table: 8.1}
{\small
%\tiny
%\scriptsize
%\footnotesize
\begin{tabular}{c|c|c|c|c|c}
$w(G)$ &$y_{11}$ &Remak 1& $w(\bar{G})$&Group&Remark 2 \\
   \hline
(a)&$z^{-1}$ & $w_{11}=\nu z^2$ & ${\left(
 \begin{array}{cc}
 \nu & w_{12}z^{-1}\\
 0& w_{22}
 \end{array}
 \right)}$& $\begin{array}{c} {\rm (V1)\ if\ }w_{22}\neq 0 \\  {\rm (V2)\ if\ }w_{22}=0\end{array}$ & $\begin{array}{c} s=(w_{22})^{-1}\\ t=w_{12}z^{-1}\\ \end{array}$\\
\hline
(b)&$w_{12}^{-1}$ & &${\left(
 \begin{array}{cc}
 0& 1\\
 w_{21}w_{12}^{-1}& 0
 \end{array}
 \right)}$ &$\begin{array}{c} {\rm (V3)\ if\ }w_{21}\neq 0\\ {\rm (V4)\ if\ }w_{21}=0\end{array}$ & $\begin{array}{c}s=-(w_{21})^{-1}w_{12}\\ \\ \end{array}$\\
\hline
(c)& $-w_{21}^{-1}$&  &${\left(
 \begin{array}{cc}
 0& 0\\
 -1& w_{22}
 \end{array}
 \right)}$ & (V5) & $s=-w_{22}$ \\
\hline
(e)&  &  &${\left(
 \begin{array}{cc}
 0& 0\\
 0& w_{22}
 \end{array}
 \right)}$ & $\begin{array}{c} {\rm (V6)\ if\ }w_{22}\neq 0\\ {\rm (V7)\ if\ }w_{22}=0\end{array}$ & $\begin{array}{c}s=w_{22}^{-1}\\ \\ \end{array}$\\
   \hline
\end{tabular}}
\caption{Subcase 2.2.2.1 in Theorem \ref{5}}
\label{table 7}
\end{table}
By Table \ref{table 7}, we get the groups of Type (V1)--(V7).
\qed
\section{The application of the classification}

Before applying the above results, we want to ensure that the classification is accurate. A useful check is to compare determination of the groups whose order is as small as possible. Determination may be the database of Magma or GAP or a group-list given in a paper.
For example, if $p=2$, then the smallest order of $G$ is $2^6-2^9$ respectively. Fortunately, Magma and GAP provide group lists for 2-groups of order less than $2^{10}$. We checked the database of Magma, and no error was found.

\begin{table}[h]
  \centering
%\caption{Properties of the groups of Type (A) in Theorem
%\ref{main}}\label{table: 8.1}
{%\scriptsize
%\footnotesize \small
\begin{tabular}{cc||cc}
\hline
Groups &the minimal possible order &Groups & the minimal possible  order\\
\hline

(B1)--(B3)&  $2^6$ &(M1)--(M9)&$p^8$\\
 (D1)--(D5) & $2^8$, $p^6$ $(p>2)$&(P1)--(P10)&$2^9$, $3^7$, $p^5$ $(p>3)$\\
 (E1)--(E9) & $p^7$ &(Q1)--(Q8)&$2^9$, $p^7$ $(p>2)$\\
  (G1)--(G3)& $2^6$, $3^6$, $p^4$ $(p>3)$ &(S1)--(S7)&$2^8$, $p^6$ $(p>2)$\\
  (I1)--(I3) & $2^6$ &(T1)--(T7)&$2^7$\\
 (J1)--(J6) &$2^7$, $p^5\ (p>2)$&(U1)--(U8)&$p^8$ \\
 (L1)--(L9)&$2^8$, $p^7\ (p>2)$&(V1)--(V7)&$p^8$\\
   \hline
\end{tabular}}
\caption{the minimal possible order of groups obtained in this paper}
\label{table 8}
\end{table}
\begin{thm}\label{d_0}Suppose that $G$ is not metacyclic, $d(G)=2$ and $G$ has at least two
$\mathcal{A}_1$-subgroups of index $p$. Then $G$ is one of following non-isomorphic groups:

\begin{itemize}
\rr{I} $\Phi(G')G_3\cong C_p$

{\rm (1)} one of $3$-groups of maximal class of order $3^4$; %(F)

{\rm (2)}  $\langle a, b, c \di a^{p^{2}}=b^{p}=c^p=1,
[a,b]=c,[c,a]=1,[c,b]=a^{\nu p} \rangle$, where $p> 3$,
\mbox{\hskip0.3in}$\nu=1$ or a fixed quadratic non-residue modular
$p$; %(G1)

{\rm (3)}  $\langle a, b, c,d \di a^{p}=b^{p}=c^p=d^p=1,
[a,b]=c,[c,a]=d,[c,b]=1,[d,a]=[d,b]=\mbox{\hskip0.3in}1 \rangle$,
where $p>3$; %(G2)

{\rm (4)}  $\langle a, b, c \di a^{p^{2}}=b^{p}=c^p=1,
[a,b]=c,[c,a]=a^{p},[c,b]=1 \rangle$, where $p>3$;%(G3)

{\rm (5)}  $\langle a, b, c,d \di a^{4}=b^{2}=c^2=d^2=1,
[a,b]=c,[c,a]=d,[c,b]=[d,a]=[d,b]=\mbox{\hskip0.3in}1 \rangle$;%(H1)

{\rm (6)}  $\langle a, b, c \di a^{8}=b^{2}=c^2=1,
[a,b]=c,[c,a]=a^{4},[c,b]=1 \rangle$; %(H2)

{\rm (7)}  $\langle a, b, c \di a^{8}=c^{2}=1, b^2=a^4,
[a,b]=c,[c,a]=b^2,[c,b]=1 \rangle$; %(H3)

{\rm (8)}  $\langle a, b, c,d \di a^{2^n}=b^{2}=c^2=d^2=1,
[a,b]=c,[c,a]=d,[c,b]=[d,a]=[d,b]=\mbox{\hskip0.3in}1 \rangle$, where $n\ge 3$;%(I1);

{\rm (9)}  $\langle a, b, c \di a^{2^{n+1}}=b^{2}=c^2=1,
[a,b]=c,[c,a]=a^{2^n},[c,b]=1 \rangle$, where $n\ge 3$;%(I2);

{\rm (10)}  $\langle a, b, c \di a^{2^n}=b^4=c^{2}=1,[a,b]=c,[c,a]=b^2,[c,b]=1 \rangle$, where $n\ge 3$;%(I3);

{\rm (11)}  $\langle a, b, c \di a^{p^{n+1}}=b^{p}=c^p=1,
[a,b]=c,[c,a]=1,[c,b]=a^{\nu p^n} \rangle$, where $p>2$ \mbox{\hskip0.3in}and $n>1$,
$\nu=1$ or a fixed quadratic non-residue modular
$p$; %(J1)

{\rm (12)}  $\langle a, b, c \di a^{p^{n+1}}=b^{p}=c^p=1,
[a,b]=c,[c,a]=a^{p^n},[c,b]=1 \rangle$, where
$p>2$ \mbox{\hskip0.3in}and $n>1$;%(J2)

{\rm (13)}  $\langle a, b, c \di a^{p^{n}}=b^{p^{2}}=c^p=1,
[a,b]=c,[c,a]=1,[c,b]=b^{p} \rangle$, where $p>2$ and \mbox{\hskip0.3in}$n>1$; %(J3)

{\rm (14)}  $\langle a, b, c \di a^{p^{n}}=b^{p^{2}}=c^p=1,
[a,b]=c,[c,a]=b^{\nu p},[c,b]=1 \rangle$, where $p>2$ and \mbox{\hskip0.3in}$n>1$, $\nu=1$ or
a fixed quadratic non-residue modular $p$; %(J4)

{\rm (15)}  $\langle a, b, c,d \di a^{p^{n}}=b^{p}=c^p=d^p=1,
[a,b]=c,[c,a]=1,[c,b]=d, [d,a]=\mbox{\hskip0.3in}[d,b]=1 \rangle$,
where $p>2$ and $n>1$; %(J5)

{\rm (16)}  $\langle a, b, c,d \di a^{p^{n}}=b^{p}=c^p=d^p=1,
[a,b]=c,[c,a]=d,[c,b]=1,[d,a]=\mbox{\hskip0.3in}[d,b]=1 \rangle$,
where $p>2$ and $n>1$. %(J6)

\rr{II} $\Phi(G')G_3\cong C_p^2$

{\rm (17)}  $\langle a, b, c,d,e \di a^{3}=b^{3}=c^3=d^3=e^3=1,
[a,b]=c,[c,a]=d,[c,b]=\mbox{\hskip0.3in}e,[d,a]=[d,b]=[e,a]=[e,b]=1
\rangle$;%(O1)

%{\rm (O2)}  $\langle a, b, c,d \di a^{3}=b^{9}=c^3=d^3=1,
%[a,b]=c,[c,a]=d,[c,b]=b^3,[d,a]=[d,b]=1 \rangle$;%(O2)只有一个$\mathcal{A}_1$子群

{\rm (18)}  $\langle a, b, c,d\di a^{9}=c^{3}=d^3=1, b^3=a^3,
[a,b]=c,[c,a]=d,[c,b]=a^{3},[d,a]=\mbox{\hskip0.3in}[d,b]=1
\rangle$;%(O3)

{\rm (19)}  $\langle a, b, c,d \di a^{9}=b^{3}=c^3=d^3=1,
[a,b]=c,[c,a]=d,[c,b]=a^{-3},[d,a]=\mbox{\hskip0.3in}[d,b]=1
\rangle$;%(O4)

{\rm (20)}  $\langle a, b, c \di a^{9}=b^9=c^3=1,
[a,b]=c,[c,a]=a^3,[c,b]=b^{3}\rangle$;%(O5)

{\rm (21)}  $\langle a, b, c\di a^{9}=b^{9}=c^3=1,
[a,b]=c,[c,a]=b^3,[c,b]=a^3,[a^3,b]=1 \rangle$;%(O6)

{\rm (22)}  $\langle a, b, c \di a^{9}=b^{9}=c^3=1,
[a,b]=c,[c,a]=b^{-3},[c,b]=a^3,[a^3,b]=1\rangle$;%(O7)

{\rm (23)} $\langle a, b, c \di a^{p^{2}}=b^{p^{2}}=c^{p}=1,
[a,b]=c,[c,a]=a^{p}b^{\nu p},[c,b]=b^{p}\rangle$,
where $p>3$, \mbox{\hskip0.3in}$\nu=1$ or a
fixed quadratic non-residue modular $p$;%(P2)

{\rm (24)} $\langle a, b, c \di a^{p^{2}}=b^{p^{2}}=c^{p}=1,
[a,b]=c,[c,a]=b^{\nu p},[c,b]=a^{-p}\rangle$, where
$p>3$, \mbox{\hskip0.3in}$\nu=1$ or a
fixed quadratic non-residue modular $p$;%(P3)

{\rm (25)} $\langle a, b, c \di a^{p^{2}}=b^{p^{2}}=c^{p}=1,
[a,b]=c,[c,a]^{1+r}=a^{p}b^{p},[c,b]^{1+r}=a^{-r p}b^{
p}\rangle$, \mbox{\hskip0.3in}where $p>3$, $r=1,2,\dots,p-2$;%(P4)

{\rm (26)} $\langle a, b, c,d \di a^{p^{2}}=b^{p}=c^{p}=d^p=1,
[a,b]=c,[c,a]=a^{p},[c,b]=d,[d,a]=\mbox{\hskip0.3in}[d,b]=1\rangle$, where $p>3$;%(P8)

{\rm (27)} $\langle a, b, c,d \di a^{p}=b^{p^{2}}=c^{p}=d^p=1,
[a,b]=c,[c,a]=b^{\nu p},[c,b]=d,[d,a]=\mbox{\hskip0.3in}[d,b]=1\rangle$,
where $p>3$, $\nu=1$ or a
fixed quadratic non-residue modular $p$;%(P9)

{\rm (28)} $\langle a, b, c,d \di a^4=b^2=c^4=d^2=1,
[a,b]=c,[c,a]=d,[c,b]=c^2,[d,a]=\mbox{\hskip0.3in}[d,b]=1\rangle$; %(R1)

{\rm (29)} $\langle a, b, c \di a^4=b^4=c^4=1,
[a,b]=c,[c,a]=b^2,[c,b]=c^2\rangle$;%(R2)

{\rm (30)} $\langle a, b, c \di a^8=b^2=c^4=1,
[a,b]=c,[c,a]=a^4,[c,b]=c^2\rangle$;%(R3)

{\rm (31)} $\langle a, b, c \di a^8=c^4=1,b^2=a^4,
[a,b]=c,[c,a]=a^4,[c,b]=c^2\rangle$;%(R4)

{\rm (32)} $\langle a, b, c \di a^8=b^2=c^4=1,
[a,b]=c,[c,a]=a^4c^2,[c,b]=c^2,[c^2,a]=1\rangle$;%(R5)

{\rm (33)} $\langle a, b, c \di a^8=b^4=1,c^2=a^4b^2,
[a,b]=c,[c,a]=b^2,[c,b]=c^2\rangle$;%(R6)

{\rm (34)} $\langle a, b, c\di a^8=b^4=1,c^2=b^2,
[a,b]=c,[c,a]=a^4b^2,[c,b]=c^2\rangle$;%(R7)

{\rm (35)} $\langle a, b, c \di a^{p^{n+1}}=b^{p^{2}}=c^p=1,
[a,b]=c,[c,a]=a^{p^n},[c,b]=b^{sp}\rangle$, where $p>2$ \mbox{\hskip0.3in}and $n>1$, $s\in F_p^*$;%(S1)

{\rm (36)} $\langle a, b, c \di a^{p^{n+1}}=b^{p^{2}}=c^p=1,
[a,b]=c,[c,a]=b^{\nu_1p},[c,b]=a^{-\nu_2 p^n}\rangle$, where \mbox{\hskip0.3in}$p>2$ and $n>1$, $\nu_1,\nu_2=1$ or a fixed quadratic non-residue modular $p$;%(S2)

{\rm (37)} $\langle a, b, c,d \di a^{p^{n+1}}=b^{p}=c^p=d^p=1,
[a,b]=c,[c,a]=d,[c,b]=a^{-\nu p^n},[d,a]=\mbox{\hskip0.3in}[d,b]=1\rangle$, where $p>2$ and $n>1$, $\nu=1$ or a fixed quadratic non-residue \mbox{\hskip0.3in}modular $p$;%(S3)

{\rm (38)} $\langle a, b, c,d \di a^{p^{n+1}}=b^{p}=c^p=d^p=1,
[a,b]=c,[c,a]=a^{p^n},[c,b]=d,[d,a]=\mbox{\hskip0.3in}[d,b]=1\rangle$, where $p>2$ and $n>1$;%(S6)

{\rm (39)} $\langle a, b, c\di a^{2^{n+1}}=b^4=1,c^2=b^2,
[a,b]=c,[c,a]=a^{2^n},[c,b]=c^2\rangle$, where \mbox{\hskip0.3in}$n\ge 3$;%{\rm (T1)}

{\rm (40)} $\langle a, b, c \di a^{2^{n+1}}=b^4=1, c^2=a^{2^n},
[a,b]=c,[c,a]=b^2,[c,b]=c^2\rangle$, where \mbox{\hskip0.3in}$n\ge 3$;%{\rm (T2)}

{\rm (41)} $\langle a, b, c,d \di a^{2^{n+1}}=b^2=d^2=1,c^2=a^{2^n},
[a,b]=c,[c,a]=d,[c,b]=\mbox{\hskip0.3in}c^2,[d,a]=[d,b]=1\rangle$,
where $n\ge 3$;%{\rm (T3)}

{\rm (42)} $\langle a, b, c \di a^{2^n}=b^4=c^4=1,
[a,b]=c,[c,a]=b^2,[c,b]=c^2\rangle$, where $n\ge 3$;%{\rm (T4)}

{\rm (43)} $\langle a, b, c,d \di a^{2^n}=b^4=d^2=1, c^2=b^2,
[a,b]=c,[c,a]=d,[c,b]=c^2,[d,a]=\mbox{\hskip0.3in}[d,b]=1\rangle$,
where $n\ge 3$;%{\rm (T5)}

{\rm (44)} $\langle a, b, c \di a^{2^{n+1}}=b^2=c^4=1,
[a,b]=c,[c,a]=a^{2^n},[c,b]=c^2\rangle$, where $n\ge 3$;%{\rm (T6)}

{\rm (45)} $\langle a, b, c,d \di a^{2^n}=b^2=c^4=d^2=1,
[a,b]=c,[c,a]=d,[c,b]=c^2,[d,a]=\mbox{\hskip0.3in}[d,b]=1\rangle$,
where $n\ge 3$;%{\rm (T7)}
\end{itemize}

\end{thm}
\demo By Lemma \ref{part1} (2), we have $\Phi(G')G_3\le Z(G)$, $\Phi(G')G_3\le C_p^2$ and $G/\Phi(G')G_3\in\mathcal{A}_1$. Since $G$ is not metacyclic, $G/\Phi(G')G_3\cong M_p(n,m,1)$ where $n>1$ for $p=2$ and $n\ge m$.

{\bf Case 1.} $\Phi(G')G_3\cong C_p$.

If $G_3\le \Phi(G')$, then, by Theorem \ref{property-3.1}, $G$ has either $0$ or $1$ $\mathcal{A}_1$-subgroup of index $p$. Hence $\Phi(G')=1$ and $G_3\cong C_p$. By Theorem \ref{property-4.1}, $m=1$. By Theorem \ref{c_p^2}, we get the groups (1)--(16).

{\bf Case 2.} $\Phi(G')G_3\cong C_p^2$.

If $G_3\cong C_p$, then, by Theorem \ref{property-6.1.0} and \ref{property-6.2.0}, $G$ has no $\mathcal{A}_1$-subgroup of index $p$. Hence $\Phi(G')\le G_3\cong C_p^2$. By Theorem \ref{property-5.0} and \ref{5}, we get the groups (17)--(44).
\qed

\begin{table}[h]
  \centering
{
%\scriptsize
%\footnotesize %\small
\tiny
\begin{tabular}{c|c||c|c||c|c||c|c}
\hline
(1)&(F)  & (13)& (J3) where $m=1$ & (25) &(P4) where $n=1$& (37) &(S3) where $m=1$\\
(2)&(G1) where $m=1$ & (14)& (J4) where $m=1$  & (26) &(P8) where $n=1$& (38) &(S6) where $m=1$\\
(3)&(G2) where $m=1$ & (15)& (J5) where $m=1$  & (27) &(P9) where $n=1$& (39) &(T1)\\
(4)&(G3) where $m=1$ & (16)& (J6) where $m=1$  & (28) &(R1)& (40) &(T2)\\
(5)&(H1) &(17)& (O1) & (29) &(R2)& (41) &(T3)\\
(6)&(H2)  &(18)& (O3) & (30) &(R3)& (42) &(T4)\\
(7)&(H3) & (19)& (O4) & (31) &(R4)& (43) &(T5)\\
(8)&(I1) & (20)& (O5) & (32) &(R5)& (44) &(T6)\\
(9)&(I2)  & (21)& (O6) & (33) &(R6)& (45) &(T7)\\
(10)&(I3) &(22)& (O7) & (34) &(R7)&&\\
(11)&(J1) where $m=1$   &(23)& (P2) where $n=1$ & (35) &(S1) where $m=1$&&\\
(12)&(J2) where $m=1$ & (24)& (P3) where $n=1$ & (36) &(S2) where $m=1$&&\\
\hline
\end{tabular}
\caption{Correspondence for Theorem \ref{d_0}}
}\label{table 7.2}
\end{table}

\end{document}